\documentclass[12pt]{amsart}

\usepackage{stmaryrd}
\usepackage{adjustbox}
\usepackage{cancel}
\usepackage{soul}
\usepackage{extarrows}

\usepackage{etoolbox}
\makeatletter
\patchcmd{\@maketitle}
  {\ifx\@empty\@dedicatory}
  {\ifx\@empty\@date \else {\vskip3ex \centering\footnotesize\@date\par\vskip1ex}\fi
   \ifx\@empty\@dedicatory}
  {}{}
\patchcmd{\@adminfootnotes}
  {\ifx\@empty\@date\else \@footnotetext{\@setdate}\fi}
  {}{}{}
\makeatother

\makeatletter
\newcommand{\linethrough}{\mathpalette\@thickbar}
\newcommand{\@thickbar}[2]{{#1\mkern0mu\vbox{
    \sbox\z@{$#1#2\mkern-1.5mu$}%
    \dimen@=\dimexpr\ht\tw@-\ht\z@+2\p@\relax 
    \hrule\@height0.5\p@ 
    \vskip\dimen@
    \box\z@}}
}
\makeatother

\usepackage{extarrows}
\usepackage[english]{babel}
\usepackage{blindtext}
\usepackage{everysel}

\EverySelectfont{%
\fontdimen2\font=0.25em
}

\usepackage{graphicx}

\allowdisplaybreaks

\usepackage[backend=bibtex,dateabbrev=false,style=alphabetic,doi=false,eprint=false,isbn=false,url=false]{biblatex}
\renewbibmacro{in:}{%
  \ifentrytype{article}{}{\printtext{\bibstring{in}\intitlepunct}}}

\usepackage{amsmath, amsthm, amssymb, mathrsfs}

\numberwithin{equation}{subsection}

\usepackage{fullpage}
\usepackage{xcolor}
\usepackage{enumerate}
\usepackage{tikz}

\usetikzlibrary{cd}

\tikzcdset{scale cd/.style={every label/.append style={scale=#1},
    cells={nodes={scale=#1}}}}

\usepackage[pdftex,plainpages=false,hypertexnames=false,pdfpagelabels]{hyperref}
\definecolor{medium-blue}{rgb}{0,0,.8}
\hypersetup{colorlinks, linkcolor={red}, citecolor={medium-blue}, urlcolor={gray}}

\DeclareMathOperator{\Aut}{Aut}

\DeclareMathOperator{\gr}{gr}
\DeclareMathOperator{\Gr}{Gr}
\DeclareMathOperator{\id}{id}

\DeclareMathOperator{\Ind}{Ind}

\DeclareMathOperator{\End}{End}
\DeclareMathOperator{\Ext}{Ext}
\DeclareMathOperator{\ovExt}{\overline{Ext}}

\DeclareMathSymbol{\mh}{\mathord}{operators}{`\-}

\DeclareMathOperator{\Rep}{Rep}

\DeclareMathOperator{\reg}{reg}
\DeclareMathOperator{\wreg}{wreg}
\DeclareMathOperator{\Sp}{Sp}
\DeclareMathOperator{\Span}{Span}

\DeclareMathOperator{\GL}{GL}
\DeclareMathOperator{\Hom}{Hom}

\DeclareMathOperator{\unr}{unr}

\DeclareMathOperator{\Res}{Res}

\newcommand{\set}[2]{\left\{#1 \middle| #2\right\}}

\newcommand{\heart}{\ensuremath\heartsuit}

\newcommand{\rec}{\mathrm{rec}}
\newcommand{\alg}{\mathrm{alg}}
\newcommand{\ab}{\mathrm{ab}}
\newcommand{\an}{\mathrm{an}}
\newcommand{\cris}{\mathrm{cris}}
\newcommand{\dR}{\mathrm{dR}}
\newcommand{\dif}{\mathrm{dif}}
\newcommand{\fs}{\mathrm{fs}}
\newcommand{\pdR}{\mathrm{pdR}}

\newcommand{\PS}{\mathrm{PS}}

\newcommand{\la}{\mathrm{la}}

\newcommand{\Fil}{{\rm Fil}}

\newcommand{\rank}{\mathrm{rank}}
\newcommand{\Rig}{\mathrm{Rig}}
\newcommand{\rig}{\mathrm{rig}}

\newcommand{\Sen}{\mathrm{Sen}}
\newcommand{\Aff}{\mathrm{Aff}}

\newcommand{\Sym}{\mathrm{Sym}}
\newcommand{\unif}{\mathrm{wu}}
\newcommand{\univ}{\mathrm{univ}}
\newcommand{\Nil}{\mathrm{Nil}}

\newcommand{\norm}{\mathrm{norm}}
\newcommand{\pr}{\mathrm{pr}}

\newcommand{\wt}{\mathrm{wt}}

\newcommand{\PGnc}{\Phi \Gamma_{\mathrm{nc}}(\phi, \mathbf{h})}

\newcommand{\DcrisKm}{D_{\cris}^{K_m}}

\newcommand{\frb}{{\mathfrak b}}

\newcommand{\frg}{{\mathfrak g}}

\newcommand{\frm}{{\mathfrak m}}

\newcommand{\frp}{{\mathfrak p}}

\newcommand{\frt}{{\mathfrak t}}

\newcommand{\frw}{{\mathfrak w}}

\newcommand{\frz}{{\mathfrak z}}

\newcommand{\vep}{{\varepsilon}}

\newcommand{\ovK}{{\overline K}}

\def\semicolon{;}
\def\applytolist#1{
    \expandafter\def\csname multi#1\endcsname##1{
        \def\multiack{##1}\ifx\multiack\semicolon
            \def\next{\relax}
        \else
            \csname #1\endcsname{##1}
            \def\next{\csname multi#1\endcsname}
        \fi
        \next}
    \csname multi#1\endcsname}

\def\calc#1{\expandafter\def\csname c#1\endcsname{{\mathcal #1}}}
\applytolist{calc}QWERTYUIOPLKJHGFDSAZXCVBNM;
\def\bbc#1{\expandafter\def\csname bb#1\endcsname{{\mathbb #1}}}
\applytolist{bbc}QWERTYUIOPLKJHGFDSAZXCVBNM;
\def\bfc#1{\expandafter\def\csname bf#1\endcsname{{\mathbf #1}}}
\applytolist{bfc}QWERTYUIOPLKJHGFDSAZXCVBNM;
\def\sfc#1{\expandafter\def\csname s#1\endcsname{{\sf #1}}}
\applytolist{sfc}QWERTYUIOPLKJHGFDSAZXCVBNM;
\def\fc#1{\expandafter\def\csname f#1\endcsname{{\mathfrak #1}}}
\applytolist{fc}QWERTYUIOPLKJHGFDSAZXCVBNM;
\def\scrc#1{\expandafter\def\csname scr#1\endcsname{{\mathscr #1}}}
\applytolist{scrc}QWERTYUIOPLKJHGFDSAZXCVBNM;

\theoremstyle{plain}
\newtheorem{thm}{Theorem}[section]
\newtheorem{cor}[thm]{Corollary}
\newtheorem{lem}[thm]{Lemma}

\newtheorem{prop}[thm]{Proposition}

\newtheorem*{claim*}{Claim}
\newtheorem*{thm*}{Theorem}
\theoremstyle{definition}
\newtheorem{defn}[thm]{Definition}

\newtheorem{ex}[thm]{Example}

\newtheorem{rem}[thm]{Remark}

\newtheorem*{exer*}{Exercise}
\newtheorem*{prob*}{Problem}
\newtheorem*{rems*}{Remarks}

\title{Change of weights operations for triangulated $(\varphi,\Gamma)$-modules}
\author{Zichuan Wang}

\address{Indiana University, Department of Mathematics, Rawles Hall, Bloomington, IN 47405, U.S.A.}
\email{wangzich@iu.edu}

\begin{document}
\maketitle

\begin{abstract}
Wu has shown the existence of ``change of weights'' operation on $(\varphi,\Gamma)$-modules in families, \cite[Prop. 3.16]{Wu_translation}.
We interpret it in the trianguline case as pullbacks with a discussion on related stacks.
Finally, we prove that it intertwines well with translation functors via a 1-1 correspondence defined by Ding \cite{Ding_crystabelline} in the non-critical crystabelline case.
\end{abstract}

\tableofcontents

\section{Introduction}
\subsection{Background and main results}

Under the locally analytic $p$-adic Langlands correspondence of $\GL_2(\bbQ_p)$ established in \cite{Colmez_localement}, the papers \cite{JLS_translation} and more generally \cite{Ding_weight} studied certain ``change of weights" operations on $2$-dimensional $p$-adic $G_{\bbQ_p}$-representations (or rather, ($\varphi,\Gamma$)-modules of rank $2$ over the Robba ring) that were shown to correspond to translation functors (\cite{Humphreys_CategoryO}) on the side of locally analytic representations of $\GL_2(\bbQ_p)$.
Then, such operations were generalized by Wu to families of $(\varphi,\Gamma)$-modules of rank $2$ over rigid spaces in \cite{Wu_translation}.
The paper \cite{Colmez_poids} can be seen as a precursor of these techniques and results (even though the connection to translation functors is not considered there).

\subsubsection{Change of weights for triangulated $(\varphi,\Gamma_K)$-modules}
It is natural to ask for a general theory of ``change of weights" for $(\varphi,\Gamma_K)$-modules.
For $\GL_2(\bbQ_p)$, using the basic constructions in the $p$-adic Langlands correspondence of Colmez, Ding equipped $(\varphi,\Gamma_{\bbQ_p})$-modules $D$ of rank $2$ over the Robba ring $\cR_{\bbQ_p,E}$ with $\mathfrak{gl}_2(\bbQ_p)$-module structures, equipped $D\otimes_E \Sym^k(E^2)$ with natural $\mathfrak{gl}_2(\bbQ_p)$-module and $(\varphi,\Gamma_{\bbQ_p})$-module structures, and developed a theory of change of weights in \cite{Ding_weight}.

For $n>2$ or a larger base field $K \neq \bbQ_p$, it is unclear to us how to equip $(\varphi,\Gamma_K)$-modules over $\cR_{K,E}$ with $\mathfrak{gl}_n(K)$-module structures, but we always have pullback and pushforward operations in the category of pairs $(D,\Fil^\bullet(D))$, where $D$ is a trianguline $(\varphi,\Gamma_{K})$-module and $\Fil^\bullet(D)$ is a triangulation on $D$, i.e., a filtration of $D$ with successive quotients of rank $1$.
The bijectivity of their induced maps between the cohomology groups in \cite[Theorem 2.22(i)]{Colmez_trianguline} is reminiscent of the fact in \cite[\S7.8]{Humphreys_CategoryO} that translation functors between dominant integral weights of the same regularity are equivalences of categories, which motivated us to consider pullback and pushforward operations as some kind of ``change of weights'' operations.

More precisely, given a triangulated $(\varphi,\Gamma_{K})$-module $(D,\Fil^\bullet(D))$ of rank $n$ over $\cR_{K,A}$ with
$$
\Fil^0(D) = 0 \subsetneq \Fil^1(D) \subsetneq \dots \subsetneq \Fil^n(D) = D,
$$
where $A$ is any affinoid algebra over a finite extension $E/\bbQ_p$ and $\Fil^i(D)$ are saturated $(\varphi,\Gamma_{K})$-submodules of $D$, 
for $0\leq j\leq n$ we consider the extension 
$$
0 \to \Fil^j(D) \to D \to D/\Fil^j(D) \to 0.
$$
Let $\Sigma_K$ be the set of embeddings from $K$ to $E$. We always assume that $|\Sigma_K|=[K:\bbQ_p]$.
For any ${\bf k} = (k_\sigma)_{\sigma} \in \bbN^{\Sigma_K}$, let $t^{\bf k} := \prod_{\sigma:K\to E} t_\sigma^{k_\sigma}\in \cR_{K,E}$, cf. \cite[Notation 6.2.7]{KPX_cohomology}.
Then,\footnote{As for the pushforward operations, starting from the extension $
0 \to \Fil^j(D) \to D \to D/\Fil^j(D) \to 0,
$
we may also consider its pushforward along the inclusion $\Fil^j(D) \subset t^{-{\bf k}}\Fil^j(D)$, resulting in a $(\varphi,\Gamma_K)$-submodule $\iota_{\bf k}(D)$ of $t^{-{\bf k}}D$ that is an extension of the form 
$$
0 \to t^{-{\bf k}}\Fil^j(D) \to \iota_{\bf k}(D) \to D/\Fil^j(D) \to 0.
$$
The pushforward $\iota_{\bf k}(D)$ becomes isomorphic to the pullback $p_{\bf k}(D)$ after a character twist (Lemma \ref{lem:push-pull}).
Hence, we mainly focus on the pullback operations $p_{\bf k}$ in the article.} we can pullback $D$ along $t^{\bf k}(D/\Fil^j(D))\subset D/\Fil^j(D)$ to get a subobject $p_{\bf k}(D,\Fil^j(D))$ of $D$:
$$
0 \to \Fil^j(D) \to p_{\bf k}(D,\Fil^j(D)) \to t^{\bf k}(D/\Fil^j(D)) \to 0.
$$
Note that if $(h_{i,\sigma})_{1\leq i \leq n, \sigma} \in (A^{n})^{\Sigma_K}$ are the Sen weights of $D$, with $(h_{i,\sigma})_\sigma \in A^{\Sigma_K}$ being the Sen weights of the $(\varphi,\Gamma_K)$-module $\Fil^i(D)/\Fil^{i-1}(D)$ of rank $1$ for $1\leq i\leq n$, then 
$$
h'_{i,\sigma} = \begin{cases}
h_{i,\sigma} &\text{ if } i\leq j, \\
h_{i,\sigma}+k_{\sigma} &\text{ if } i>j,
\end{cases}
$$
are the Sen weights of $p_{\bf k}(D,\Fil^j(D))$.

One can then ask whether these weight-shifting operations descend to the category of trianguline $(\varphi,\Gamma_K)$-modules (without fixing a filtration).
But it is easy to see that in general they depend on the choice of triangulation, as the following example shows.

When $K=\mathbb{Q}_p$, for an integer $i$, we write
$
x^i: \mathbb{Q}_p^{\times} \rightarrow E^{\times}, a \mapsto a^i,
$
and denote by
$
R_{\mathbb{Q}_p, E}(x^i)
$
the $(\varphi, \Gamma_{\mathbb{Q}_p})$-module of rank $1$ associated to this character. 
Explicitly, it is the free $R_{\mathbb{Q}_p, E}$-module $R_{\mathbb{Q}_p, E} e_i$ of rank $1$ with
$
\varphi(e_i)=p^i e_i
$
and
$\gamma(e_i)=\vep(\gamma)^i e_i
$
for $\gamma \in \Gamma_{\mathbb{Q}_p}$, where $\vep: \Gamma_{\mathbb{Q}_p} \xrightarrow{\sim} \mathbb{Z}_p^{\times}$ is the cyclotomic character. 
Then, $R_{\mathbb{Q}_p, E}(x^i)$ has Sen weight $i$.
Let $D = \cR_{\bbQ_p,E}(x)\oplus \cR_{\bbQ_p,E}(x^2)$ be split\footnote{We refer the reader to Remark \ref{rem:non-unique_triangulation} for a non-split example.} of rank $2$ and Hodge-Tate-Sen weights $\{1,2\}$.
Then we have two triangulations
\begin{align*}
&0 \to \cR_{\bbQ_p,E}(x) \to D \to \cR_{\bbQ_p,E}(x^2) \to 0 ,\\
&0 \to \cR_{\bbQ_p,E}(x^2) \to D \to \cR_{\bbQ_p,E}(x) \to 0.
\end{align*}
Applying pullbacks with $k>0$, we get $\cR_{\bbQ_p,E}(x)\oplus \cR_{\bbQ_p,E}(x^{k+2})$ and $\cR_{\bbQ_p,E}(x^2)\oplus \cR_{\bbQ_p,E}(x^{k+1})$ respectively, which are of different Hodge-Tate-Sen weights $\{k+2,1\}$ and $\{k+1,2\}$.

\subsubsection{Using Sen polynomials}
However, \cite[Proposition 3.16]{Wu_translation} implies that this weight issue is the only obstruction for the pullback to be independent of the triangulation. 
Let $K=\bbQ_p$ here for simplicity, and let $A$ be any affinoid $E$-algebra.
Wu showed that if the Sen polynomial $P_{\Sen,D}(T)\in A[T]$ of $D$ admits a factorization
$$
P_{\Sen,D}(T) = Q(T) S(T)
$$
with monic $Q(T),S(T) \in A[T]$, with $(Q,S)=1$, then there exists a unique $(\varphi,\Gamma)$-submodule $D'$ of $D$ satisfying
$$
tD \subset D' \subset D
$$
and having Sen polynomial 
$
Q(T-1)S(T).
$
The proof used Beauville-Laszlo gluing \cite[Proposition A.3]{Wu_translation} to show that, under the equivalence of $(\varphi,\Gamma_K)$-modules and $\Gamma_K$-equivariant vector bundles on the Fargues-Fontaine curve $X_{K_\infty,A}$ \cite[Theorem 5.1.5]{EGH_CatLLC}, with respect to the canonical decomposition of the Sen module 
$$
D_\Sen(D) = \ker(Q(\Theta_\Sen))\oplus \ker(S(\Theta_\Sen)),
$$
there is a unique modification of the $(\varphi,\Gamma)$-module $D$ at $\infty$ on the curve $X_{K_\infty,A}$ such that the resulting equivariant subbundle $D' \subset D$ has the Sen polynomial $P_{\Sen,D'}(T)=Q(T-1)S(T)$.

If $D$ is a trianguline $(\varphi,\Gamma_{K})$-module over $\cR_{K,A}$ and $\Fil^\bullet(D)$ is a triangulation on $D$, then for any $1\leq j\leq n$, the extension
$$
0 \to \Fil^j(D) \to D \to D/\Fil^j(D) \to 0
$$
induces a splitting of the Sen polynomial
\begin{equation}
\label{eqn:Sen_poly_Filj}
P_{\Sen,D}(T) = P_{\Sen,D/\Fil^j(D)}(T)\cdot P_{\Sen, \Fil^j(D)}(T),
\end{equation}
and $p_1(D,\Fil^j(D))$ is a $(\varphi,\Gamma)$-submodule of $D$ satisfying
$$
tD \subset p_1(D,\Fil^j(D)) \subset D
$$
and having Sen polynomial $P_{\Sen,D/\Fil^j(D)}(T-1)P_{\Sen, \Fil^j(D)}(T)$.
If \eqref{eqn:Sen_poly_Filj} is a comaximal factorization, then $p_1(D,\Fil^j(D))$ is the unique $(\varphi,\Gamma)$-submodule of $D$ containing $tD$ of the Sen polynomial $P_{\Sen,D/\Fil^j(D)}(T-1)P_{\Sen, \Fil^j(D)}(T)$.
Moreover, in this case, for any triangulation $\Fil^\bullet(D)'$ on $D$ inducing the same splitting as \eqref{eqn:Sen_poly_Filj}, by the uniqueness, we have that
$$
p_1(D,\Fil^j(D)') = p_1(D,\Fil^j(D)).
$$

Consider those $(\varphi,\Gamma_K)$-modules $D$ on $\Sp(A)$ that are Tate-fpqc locally trianguline.
Assume moreover that at each point $z\in \Sp(A)$, the fiber $D_z := D\otimes_A k(z)$ is ``weight-uniform trianguline" in the sense that all possible triangulations on $D_{z}\otimes_{k(z)}L$ induce the same ordering on the set of Sen weights of $D_z$, for all finite extensions $L/k(z)$.
Any local triangulation on $D$ locally gives us a factorization of the Sen polynomial. 
After imposing necessary restrictions on its Sen weights at geometric fibers, we can perform such ``modifications at $\infty$'' iteratively on $D$ in an invertible way that is independent of the choice of covers and local triangulations, whenever the resulting movement of our ordered Sen weights in the weight space does not meet any of the relevant walls, cf. Theorem \ref{thm:Wu_3.16} for our exposition on Wu's result and Theorem \ref{thm:comparing_Wu} for discussion in the weight-uniform trianguline case.

\subsubsection{The point of view of the analytic Emerton-Gee stacks}

Recall from \cite[\S5.3]{EGH_CatLLC} that over the category $\Rig_E$ of rigid analytic spaces over $E$ equipped with the Tate-fpqc topology, we have the moduli stack $\fX_n$ of rank $n$ $G_K$-equivariant vector bundles over the Fargues-Fontaine curve $X_{\overline{K}}$, and the stack $\fX_B$ of $G_K$-equivariant $B$-bundles on $X_{\overline{K}}$, where $B$ denotes the Borel subgroup of $G=\GL_n$ consisting of upper triangular invertible matrices.
Then, by the equivalence \cite[Theorem 5.1.5]{EGH_CatLLC}, a triangulated $(\varphi,\Gamma_{K})$-module $(D,\Fil^\bullet(D))$ of rank $n$ over $\cR_{K,A}$ defines a point in $\fX_B(A)$, while the $(\varphi,\Gamma_{K})$-module $D$ defines a point in $\fX_n(A)$.

We introduce in Definition \ref{defn:uniform_tri}  ``weight-uniform trianguline substack'' $\fX_n^{\unif}$ of $\fX_n$, and some additional substacks $\fX_n^{\sigma\mh\unif,i} \subset \fX_n^{\sigma\mh\unif}$ characterized by the property that, for any triangulation on $D$, the first $n-i$ $\sigma$-Sen weights are distinct from the last $i$ $\sigma$-Sen weights at any $z\in \Sp(A)$.
The pullback operator $p_{i,\sigma}$ (\S\ref{section:pullback_stack}) that increases the last $i$ $\sigma$-Sen weights by $1$ and leaves all other Sen weights invariant  descends to a map from $\fX_n^{\sigma\text{-}\unif,i}$ to $\fX_n$. 
Then, we deduce the following theorem, cf. Definition \ref{defn:Xn_SwuIk} (and Remark \ref{rem:SwuIk}) for the precise meaning of our notation.

\newtheorem*{theoremA}{Theorem A}
\begin{theoremA}
\textit{The pullback maps $p_{i,\sigma}:\fX_B\to \fX_B$ descend to canonical morphisms of stacks}
$$p_{i,\sigma} : \mathfrak{X}_n^{\sigma\mh\unif,i} \longrightarrow \mathfrak{X}_n$$
\textit{such that for $S\subset \Sigma_K$, $I= \prod_{\sigma\in S}I_\sigma \subset \{1,\dots,n\}^S$ and ${\bf k} = (k_{i,\sigma})_{\sigma\in S, i \in I_\sigma} \in \bbN^{I}$,}
$$
p_{\bf k}:= \prod_{i,\sigma} (p_{i,\sigma})^{k_{i,\sigma}} : \fX_n^{S\text{-}\unif, I, {\bf k}} \xlongrightarrow{\sim} \fX_n^{S\text{-}\unif, I, {\bf -k}} 
$$
\textit{are isomorphisms between these weight-uniform trianguline substacks, where the change of Sen weights does not change the regularity of the weights.}
\end{theoremA}

\subsubsection{}

We also mention that, for such directions of changing the weights, Wu obtained in \cite[\S3]{Wu_translation} general geometric results: the stack $\fX_n$ of rank $n$ $(\varphi,\Gamma_K)$-modules at integral Hodge-Tate weights ${\bf h} = (h_{i,\sigma})\in (\bbZ^{n})^{\Sigma_K}$ is described using a ``product formula" of the form (if $K=\bbQ_p$)
$$
(\fX_n)^\wedge_{\bf h} \cong (\fX)_0^{\wedge} \times_{\frg / \GL_n} \widetilde{\frg}_{P_{{\bf h}}} / \GL_n\,,
$$
which then induces change of weights maps 
$$f_{\bf h, \bf h'}:(\fX_n)^\wedge_{\bf h} \to (\fX_n)_{\bf h'}^\wedge\,,
$$ 
and this can be formulated for non-integral weights using local models developed in \cite[Ch. 5]{Wu_thesis}.
Moreover, change of weights maps $f_{\bf h, \bf h'}$ exist at arbitrary weights ${\bf h}$ when changing from ${\bf h}$ to ${\bf h'}$ does not increase the regularity, cf.
\cite[\S1.3]{Wu_translation} for a discussion and \cite[\S3.3]{Wu_translation} for details.

{\bf Question.} What can be said about the geometry of $\fX_n^{\unif}$?

\subsubsection{}
The change-of-weight operations are supposed to be the kind of operations on the side of $p$-adic Galois representations (or rather $(\varphi,\Gamma_K)$-modules $D$) that correspond to the effect of translation functor (as initially introduced in \cite{JLS_translation}) on the side of locally analytic $\GL_n(K)$-representations under a (hypothetical) $p$-adic locally analytic Langlands correspondence $D\mapsto \Pi^{\an}(D)$.
Roughly speaking, the change-of-weight operations and the translation functors are meant to be related by a formula of the type 
\begin{equation}
\label{eqn:1.1.5}
\Pi^{\an}(\text{(change-of-weight operation)}(D)) = \text{(translation functor)}(\Pi^{\an}(D)),
\end{equation}
which has been shown in certain cases by \cite{JLS_translation} and \cite{Ding_weight} for $\GL_2(\bbQ_p)$.
We also remark that there are analogous statements in the local Langlands correspondence for real reductive groups which relate operations on ``$L$-parameters'' (on the Galois side) with translation functors on the corresponding representations of the reductive group, cf. \cite[\mbox{}{16.6}]{ABV-RealLanglands}.\footnote{Note that in \cite[\mbox{}{16.6}]{ABV-RealLanglands} the objects on the Galois side are rather complete geometric parameters and on the automorphic side these are called ``final limit characters''.}

\subsubsection{Relation with translation functor under Ding's crystabelline correspondence}

A class of points of $\fX_n^{\unif}(E)$ consists of those non-critical crystabelline $(\varphi,\Gamma_K)$-modules over $\cR_{K,E}$.
Recently, for non-critical crystabelline $(\varphi,\Gamma_K)$-modules $D$ over $\cR_{K,E}$ of regular Hodge-Tate weights, Ding constructed in \cite{Ding_crystabelline} locally $\bbQ_p$-analytic representations $\pi_{\min}(D) \subset \pi_{\fs}(D)$ of $\GL_n(K)$, which are expected to occur as closed  subrepresentations of the hypothetical $\Pi^{\an}(D)$.
They are constructed as extensions of locally algebraic representations $\pi_{\alg}(\phi,{\bf h})$ by locally analytic representations $\pi(\phi,{\bf h})$ that only depend on $D[1/t]$ and the Hodge-Tate weights $\bf h$.
If $K=\bbQ_p$, these extension classes can ``recover the Hodge filtration'' and hence determine $D$, cf. \cite[Theorem 3.34]{Ding_crystabelline}.
The pullback $D'$ is a non-critical crystabelline $(\varphi,\Gamma_K)$-submodule of $D$ whenever this change of weights preserves the regularity.
It is then natural to expect that under $\pi :=\pi_{\min}$ or $\pi_{\fs}$, pulling back from $D$ to $D'$ corresponds to translating from $\pi(D)$ to $\pi(D')$.
We prove the expected intertwining \eqref{eqn:1.1.5} of the two kinds of weight-shifting operations:

\newtheorem*{theoremB}{Theorem B}
\begin{theoremB}
\label{thm:B}
\textit{Let $D$ be a non-critical crystabelline $(\varphi, \Gamma_{K})$-module with regular Sen weights ${\bf h}$.
Let $p_{\mathbf{k}}(D) = f_{{\bf h},{\bf h'}}(D)$ be the $(\varphi, \Gamma_{K})$-submodule of $D$ obtained by a sequence of pullback operators $p_{\bf k}=\prod_{i,\sigma} (p_{i,\sigma})^{k_{i,\sigma}}$ such that its weights ${\bf h'}$ are still regular.
Let $\lambda':={\bf h'}-\theta$ and $\lambda:={\bf h}-\theta$ be the ``automorphic weights'' with $\theta :=(0,-1,\dots, -(n-1))_{\sigma\in\Sigma_K} \in (\bbZ^n)^{\Sigma_K}$. Then,
$$
T_\lambda^{\lambda'}(\pi_\bullet(D)) = \pi_{\bullet}(p_{\bf k}(D)) = \pi_{\bullet}(f_{{\bf h},{\bf h'}}(D))
$$
as locally $\bbQ_p$-analytic representations of $\GL_n(K)$, for $\bullet \in \{{\min},{\fs}\}$.
}
\end{theoremB}

Since translation functors between regular weights are equivalences of categories, Theorem B follows from various results in \cite{Ding_crystabelline}, \cite{JLS_translation} and \cite{Wu_translation}, with some additional work.

\subsection{Structure of the paper}
In \S\ref{section:prelim}, we recall basics of $(\varphi,\Gamma_K)$-modules over affinoid algebras and their cohomology, generalizing a pointwise result \cite[Theorem 2.22(i)]{Colmez_trianguline} to the affinoid coefficients in Lemma \ref{lem:rank_one_affinoid}(i), leading to Theorem \ref{thm:p_k-invertible} for trianguline families.
In \S\ref{section:triangulations}, we review triangulations on generic crystabelline $(\varphi,\Gamma_K)$-modules over fields, and observe that ``generically, a trianguline $(\varphi,\Gamma_K)$-module has a unique non-split triangulation.''
In \S\ref{section:pullback_stack}, we discuss Wu's change of weights maps in general and in the weight-uniform trianguline families (cf. Theorem \ref{thm:Wu_3.16} = \cite[Proposition 3.16]{Wu_translation}, and Theorem \ref{thm:comparing_Wu} = Theorem A), which is followed by a discussion on when étaleness can be preserved up to twist for local Artinian $E$-algebras $A\in \cC_{E}$ of residue field $E$.
Finally, \S\ref{section:translation_functor} is devoted to Theorem B.
\vskip8pt

\subsection{Acknowledgements}
I heartily thank Yiwen Ding for helpful discussions regarding my questions on his work. Special thanks go to Zhixiang Wu for drawing my attention to his work on change of weights maps and for his valuable comments on an earlier draft of this paper. 
It is my great pleasure to thank my PhD advisor, Matthias Strauch, for his stimulating questions, persistent patience, and constant encouragement.
I also thank the anonymous referee for their careful reading of the manuscript and many helpful comments.

\section{Preliminaries}
\label{section:prelim}

\subsubsection{Notations and conventions}

Let $K$ be a finite extension of $\bbQ_p$ of ramification index $e$ and inertial degree $f$, with maximal unramified subfield $K_0$ and a chosen uniformizer $\pi_K$.
Let $G_K:=G(\overline{K}/K)$
be the absolute Galois group of $K$, $W_K$ be the Weil group of $K$, $
\rec_K: K^\times \xrightarrow{\sim} W_K^{\ab}
$ be the local reciprocity map normalized by sending $\pi_K$ to a geometric Frobenius element, and let $\Sigma_K$ be the set of all $\bbQ_p$-algebra embeddings of $K$ into $\overline{\bbQ_p}$.
Let $E$ be a finite extension of $\bbQ_p$ such that all the $\bbQ_p$-algebra embeddings of $K$ into $\overline{\bbQ_p}$ factor through $E$, which will be our coefficient field.
We allow $E$ to be enlarged at will.

Let $\Rig_E$ be the category of rigid $E$-analytic spaces, and $\Aff_E$ the category of affinoid $E$-algebras.
For $X\in \Rig_E$, let $\cR_{K,X}$ be the relative Robba ring of $K$ over $X$, and write $\cR_{K,A}:=\cR_{K,\Sp(A)}$ for affinoid $A\in \Aff_E$.
Let $\Gamma_K:=G(K(\mu_{p^\infty})/K)$.
For a review of (generalized) $(\varphi,\Gamma_K)$-modules $D$ over the Robba ring $\cR_{K,X}$, we refer the reader to \cite[\S2.1]{Bergdall}.
For a $(\varphi,\Gamma_K)$-module $D$ on $X$, and $z\in X$ a rigid point, we write $D_z := D\otimes_{\cO_X} k(z)$ for the fiber of $D$ at $z$. 
Similarly, if $\delta: K^\times \to \Gamma(X,\cO_{X})^\times$ is a locally $\bbQ_p$-analytic character, we write $\delta_z: K^\times \to \Gamma(X,\cO_{X})^\times \to k(z)^\times$ for the evaluation of $\delta$ at $z\in X$.
The notation $D = (Q - P)$ is used to indicate that $D$ is an extension of $P$ by $Q$, and the notation $D=(Q_1-\dots-Q_r)$ indicates that there exists an increasing filtration $(\Fil^i(D))_{0\leq i\leq r}$ on $D$ by subobjects such that $\Fil^0(D)=0$, $\Fil^r(D)=D$, and $\Fil^{i}(D)/\Fil^{i-1}(D) \cong Q_i$ for $1\leq i\leq r$. 

We set $\bbN := \bbZ_{\geq 0}$.
By local class field theory, the $p$-adic cyclotomic character $\vep:G_K\to \bbZ_p^\times$ is equivalently given by
$$
\varepsilon: K^\times \longrightarrow E^\times,\quad x\mapsto N_{K|\bbQ_p}(x)|N_{K|\bbQ_p}(x)|_p,
$$
which has $\sigma$-Sen weight $+1$ for all $\sigma\in \Sigma_K$.
For a $[K:\bbQ_p]$-tuple of integers ${\bf k} = (k_\sigma)_{\sigma\in \Sigma_K} \in \bbZ^{\Sigma_K}$, set $t^{\bf k}:=\prod_{\sigma\in \Sigma_K} t_{\sigma}^{k_\sigma}\in \cR_{K,E}$, where $t_\sigma\in \cR_{K,E}$ are the Lubin-Tate elements defined up to units in \cite[Notation 6.2.7]{KPX_cohomology}.
Let $x_\sigma$ be the embedding $\sigma: K^\times \to E^\times$ viewed as a character of $K^\times$, and set $x^{\bf k} := \prod_{\sigma\in \Sigma_K} x_{\sigma}^{k_\sigma}$.
Then, the $(\varphi,\Gamma_K)$-module $t^{\bf k}\cR_{K,A} \cong \cR_{K,A}(x^{\bf k})$ is free of rank 1 with $\sigma$-Sen weight $k_\sigma$ for each $\sigma\in \Sigma_K$.

\subsubsection{Extensions of $(\varphi,\Gamma_K)$-modules}

We recall certain operations on extensions of $(\varphi,\Gamma_K)$-modules and the induced linear maps on $(\varphi,\Gamma_K)$-cohomology.
\begin{defn}
For ${\bf k}\in \bbN^{\Sigma_K}$ and $(\varphi,\Gamma_K)$-modules $D_1, D_2$ over $\cR_{K,A}$, there are two natural $A$-linear maps between Yoneda Ext groups:
\begin{itemize}
\item 
pushing out along $D_1 \hookrightarrow t^{-{\bf k}}D_1$ defines a map
$$
\iota_{\bf k}: \Ext^1(D_2, D_1) \longrightarrow \Ext^1(D_2, t^{\bf -k} D_1),
$$
\item
pulling back along $t^{\bf k} D_2 \hookrightarrow D_2$ defines a map
$$
p_{\bf k}: \Ext^1(D_2, D_1) \longrightarrow \Ext^1(t^{\bf k} D_2, D_1).
$$
\end{itemize}
\end{defn}

The following lemma shows that these two maps are related by the ``twisting'' isomorphism  
$$
x^{\bf k}:
\Ext^1(D_2,t^{\bf -k}D_1) \overset{\simeq}{\longrightarrow} \Ext^1(t^{\bf k} D_2, D_1)
$$
induced by twisting by the character $x^{\bf k}:K^\times \to A^\times$, so that 
$
x^{\bf k} \circ \iota_{\bf k} = p_{\bf k}.
$

\begin{lem}
\label{lem:push-pull}
Let ${\bf k}=(k_\sigma)_\sigma \in \bbN^{\Sigma_K}$ and let
$$
0 \to D_1 \xrightarrow{i} D \xrightarrow{\pi} D_2 \to 0
$$
be an exact sequence of $(\varphi,\Gamma_K)$-modules over $\cR_{K,A}$.
Then, the pushout $\iota_{\bf k}(D)$ of $D$ along $D_1\hookrightarrow t^{\bf -k}D_1$ and the pullback $p_{\bf k}(D)$ of $D$ along $t^{\bf k} D_2\hookrightarrow D_2$ are related by a commutative diagram
$$
\begin{tikzcd}
0 \arrow[r] & t^{-{\bf k}} D_1 \arrow[r] & \iota_{\bf k}(D) \arrow[r] & D_2 \arrow[r] & 0 \\
0 \arrow[r] & D_1 \arrow[d, Rightarrow, no head] \arrow[u, hook] \arrow[r,"i"] & D  \arrow[u, dashed] \arrow[r,"\pi"] & D_2 \arrow[r] \arrow[u, Rightarrow, no head] & 0 \\
0 \arrow[r] & D_1 \arrow[r] & p_{\bf k}(D) \arrow[u, dashed] \arrow[r] & t^{\bf k} D_2 \arrow[r] \arrow[u, hook] & 0
\end{tikzcd}
$$
with exact rows and injective columns, from which we deduce 
$
p_{\bf k}(D) = t^{\bf k}\iota_{\bf k}(D) \cong \iota_{\bf k}(D)(x^{\bf k}).
$
\end{lem}
\begin{proof}
Explicitly, $p_{\bf k}(D) = \pi^{-1}(t^{\bf k}D_2)$, and $\iota_{\bf k}(D) = (t^{-{\bf k}}D_1) \oplus_{D_1} D$ is an amalgamated sum over $D_1$.
By the first row, we see that $t^{\bf k}\iota_{\bf k}(D)$ is a $(\varphi,\Gamma_K)$-submodule of $D[1/t]$ containing $D_1$ with the associated quotient being $t^{\bf k} D_2$, so it equals $p_{\bf k}(D)$ by the third row.
\end{proof}

\noindent
By Lemma \ref{lem:push-pull}, the pullback $p_{\bf k}$ is injective/surjective/zero if and only if the pushout $\iota_{\bf k}$ is.
We recall the following cohomological interpretation of the extensions.

\begin{lem}
\label{lem:H1-Ext1}
Let $D_1,D_2$ be $(\varphi,\Gamma_K)$-modules over $\cR_{K,A}$.
Then, 
$$
H^1(D_2^\vee \otimes D_1) \cong \mathrm{Ext}^1(D_2, D_1),
$$ 
where $H^1$ is computed using the Herr complex
$$
\cC^\bullet(D): \quad D \xrightarrow{(\varphi-1, \gamma-1)} D\oplus D \xrightarrow{(\gamma-1)\oplus(1-\varphi)} D
$$
for $D:= D_2^\vee \otimes D_1 \cong \Hom_{\cR_{K,A}}(D_2,D_1)$, where $\gamma\in \Gamma_K$ is a topological generator.
\end{lem}
\begin{proof}
This is essentially \cite[Lemma 5.1.2]{EG_modulistack}, but adapted to our convention and notation. 

Let $M$ be an extension of $D_2$ by $D_1$ as $(\varphi,\Gamma_K)$-modules over $\cR:=\cR_{K,A}$ given by
$$
0\to D_1 \xrightarrow{i} M \xrightarrow{\pi} D_2 \to 0.
$$
It splits on the level of $\cR$-modules.
We choose any $\cR$-linear section $s: D_2\to M$, which is unique up to an element $h \in D = \Hom_{\cR}(D_2,D_1)$.
Using the section $s$, we can write 
$$
\varphi_M = 
\begin{pmatrix}
\varphi_{D_1} & f\circ \varphi_{D_2} \\
 & \varphi_{D_2}
\end{pmatrix}
\quad \text{ and } \quad 
\gamma_M = 
\begin{pmatrix}
\gamma_{D_1} & g\circ \gamma_{D_2} \\
 & \gamma_{D_2}
\end{pmatrix}
$$
for uniquely determined $f,g \in \Hom_{\cR}(D_2,D_1) = D$ since $D_2$ has an $\cR$-basis in $\varphi(D_2)$.

That is, for any $x\in D_2$, we have
$$
f(\varphi_{D_2}(x)) = (\varphi_M- s\circ \varphi_{D_2}\circ \pi)(s(x)) = \varphi_M(s(x)) - s(\varphi_{D_2}(x)) \in \ker(\pi) = D_1.
$$
Similarly, since $\gamma$ acts invertibly on $\cR$, we have
$$
g = (\gamma_{M}- s\circ \gamma_{D_2}\circ \pi) \circ s \circ \gamma_{D_2}^{-1} = \gamma_M \circ s \circ \gamma_{D_2}^{-1} - s \in \Hom_{\cR}(D_2,D_1).
$$
The commutativity $\varphi_M \circ \gamma_M = \gamma_M \circ \varphi_M$ is then equivalent to the equality
$$
\varphi_{D_1} g \gamma_{D_2} + f\varphi_{D_2}\gamma_{D_2} = \gamma_{D_1} f \varphi_{D_2} + g\gamma_{D_2} \varphi_{D_2} \in \Hom_\cR(D_2,D_1),
$$
which, by precomposing with $\gamma_{D_2}^{-1}$ on $D_2$, is equivalent to
$$
\varphi_{D_1}g - g\varphi_{D_2} = \gamma_{D_1}f\gamma_{D_2}^{-1}\varphi_{D_2} - f\varphi_{D_2}.
$$
By the definition of $\Hom_\cR(D_2,D_1)$ as $(\varphi,\Gamma_K)$-module, this last displayed equation is the same as $(\varphi_D-1).g = (\gamma_D -1).f \in D$, hence $(f,g) \in \ker((\gamma-1)\oplus(1-\varphi))$ is a $1$-cocycle.

Modifying the section $s$ by $h\in D$ results in another 1-cocycle $(f',g')$ such that
$$
\begin{pmatrix}
\varphi_{D_1} & f'\circ \varphi_{D_2} \\
 & \varphi_{D_2}
\end{pmatrix}
=
\begin{pmatrix}
1 & -h \\
 & 1
\end{pmatrix}
\begin{pmatrix}
\varphi_{D_1} & f\circ \varphi_{D_2} \\
 & \varphi_{D_2}
\end{pmatrix}
\begin{pmatrix}
1 & h \\
 & 1
\end{pmatrix}
$$
(and a similar identity involving $g,g',\gamma_{D_1}$ and $\gamma_{D_2}$) which is equivalent to 
$$
f'\varphi_{D_2} = f \varphi_{D_2} -h\varphi_{D_2} + \varphi_{D_1}h \quad (\text{resp. }g' = g  -h + \gamma_{D_1}h\gamma_{D_2}^{-1}).
$$
Thus, $f' = f+(\varphi_D-1).h$, $g' = g + (\gamma_D-1).h$, and $(f',g')-(f,g)$ is a $1$-coboundary.
\end{proof}

Given Lemma \ref{lem:H1-Ext1}, we see that $\iota_{\bf k}$ and $p_{\bf k}$ induce the same map
$$
H^1(D) \longrightarrow H^1(t^{-\bf k}D)
$$
for $D := D_2^\vee \otimes D_1$, if we identify the codomains $H^1(D_2^\vee \otimes t^{-\bf k} D_1) \simeq H^1((t^{\bf k}D_2)^\vee \otimes D_1)$ via $x^{\bf k}$.

\begin{ex}
\label{ex:colmez}
When $D_2 = \cR_{\bbQ_p}$, Lemma \ref{lem:H1-Ext1} was proven by Colmez \cite[\S2.1]{Colmez_trianguline}: take $D_2=\cR$ and $D_1=\cR(\delta)$; then the isomorphism 
$$
\mathrm{Ext}^1(\cR, \cR(\delta)) \xrightarrow{\simeq } H^1(\delta)
$$
given in \cite[\S2.1]{Colmez_trianguline} is as follows: given an extension $M$, let $e\in M$ be a lift of the vector $1 \in \cR$, and we associate to $M$ (the class of) the 1-cocycle $[(\varphi_M-1)e, (\gamma_M-1)e] \in R(\delta)\oplus R(\delta)$.
This is the map constructed in Lemma \ref{lem:H1-Ext1}. 
Indeed, we have $s:1\mapsto e$, $\varphi_{D_2}=\gamma_{D_2}=\id$, so we deduce $f = \varphi_M(e) - e = (\varphi_M-1)e$.
Likewise, $g=(\gamma_M-1)e$, as desired.
\end{ex}

\subsubsection{}

We review some results on $(\varphi,\Gamma_K)$-modules of rank 1 and their cohomology.

\begin{thm}
\label{thm:rank-one_classification}
For any $(\varphi,\Gamma_K)$-module $D$ of rank $1$ over a rigid analytic space $X$, there exist a unique continuous character $\delta:K^\times \to \Gamma(X,\cO_X)^\times$ and, up to isomorphism, a unique line bundle $\cL$ on $X$ such that 
$$
\cR_{K,X}(\delta) \otimes_{\cO_X} \cL \cong D.
$$
In fact, the bundle $\cL$ is given by $H^0(D(\delta^{-1})) \cong  \Hom_{\varphi,\Gamma_K}(\cR_{K,X}(\delta),D)$, and the canonical map
$$
\cR_{K,X}(\delta) \otimes_{\cO_X} H^0(D(\delta^{-1})) \longrightarrow D
$$
is an isomorphism.
\end{thm}
\begin{proof}
This is {\cite[Theorem 6.2.14]{KPX_cohomology}}. For the definition of the $(\varphi,\Gamma_K)$-module $\cR_{K,X}(\delta)$ of rank $1$, see \cite[Construction 6.2.4]{KPX_cohomology}.
\end{proof}

\begin{defn}
Let $A$ be an affinoid $E$-algebra.
Let $\delta: K^\times \to A^\times$ be a continuous character, which is automatically locally $\bbQ_p$-analytic by \cite[Proposition 8.3]{Buzzard_eigenvar}.
Then, its derivative at $1$ defines a $\bbQ_p$-linear map $d\delta: K \to A$, and hence an $A$-linear map 
$
K\otimes_{\bbQ_p} A \to A.
$
Via the map 
\begin{equation}
\label{eqn:Galois_decoupling}
K\otimes_{\mathbb{Q}_p} A \xrightarrow{\sim} \prod_{\sigma\in \Sigma_K} A, \quad x\otimes y \mapsto (\sigma(x)y)_{\sigma\in \Sigma_K},
\end{equation}
we have an isomorphism 
$$
\Hom_{\bbQ_p}(K,A)
\cong 
\Hom_{A}(K\otimes_{\bbQ_p}A, A) \cong 
\Hom_{A}(\prod_{\sigma\in \Sigma_K} A, A) 
=A^{\Sigma_K}.
$$
The image of $d\delta$  is a $[K:\bbQ_p]$-tuple $\wt(\delta):=(\wt_\sigma(\delta))_{\sigma\in \Sigma_K} \in A^{\Sigma_K}$, which we call the {\bf weight} of the character $\delta$. 
We also call $\wt_\sigma(\delta)$ the {\bf $\sigma$-weight} of $\delta$.
By \cite[Lemma 6.2.12]{KPX_cohomology}, for any continuous character $\delta:K^\times \to A^\times$, the ($\sigma$-)Sen weight of $\cR_{K,A}(\delta)$ is the ($\sigma$-)weight of $\delta$.
\end{defn}

\begin{lem}
\label{lem:rank_one_coh}
Let $\delta: K^\times \to E^\times$ be a continuous character.
\begin{enumerate}[(i)]
\item 
\label{lem:Colmez-iota}
For ${\bf k}\in \bbN^{\Sigma_K}$, if $\wt_\sigma(\delta)\notin\{1,\dots, k_\sigma\}$ for each $\sigma\in \Sigma_K$, then 
$$
\iota_{\bf k}: H^1(\delta) \to H^1(x^{-{\bf k}}\delta)
$$ 
is an isomorphism.

\item
\label{lem:Liu-H0H2}
We have 
\begin{itemize}
\item 
$
\dim_E H^0(\delta) = \begin{cases}
1 &\text{ if } \delta = x^{-\bf k} \text{ for some } {\bf k} \in \bbN^{\Sigma_K}, \\
0 &\text{ otherwise}.
\end{cases}
$
\item 
$
\dim_E H^2(\delta) = \begin{cases}
1 &\text{ if } \delta = \left(N_{K|\bbQ_p}|N_{K|\bbQ_p}|_p\right)x^{\bf k} \text{ for some } {\bf k} \in \bbN^{\Sigma_K}, \\
0 &\text{ otherwise}.
\end{cases}
$
\item 
$
\dim_E H^1(\delta) = \begin{cases}
[K:\bbQ_p]+1 &\text{ if either $H^0$ or $H^2$ does not vanish}, \\
[K:\bbQ_p] &\text{ otherwise}.
\end{cases}
$
\end{itemize}

\item
Any nonzero $(\varphi,\Gamma_K)$-submodule of $\cR_{K,E}(\delta)$ must be of the form $t^{\bf k}\cR_{K,E}(\delta)$ for some ${\bf k}\in \bbN^{\Sigma_K}$.

\end{enumerate}
\end{lem}

\begin{proof}
The first statement is {\cite[Lemma 3.3.3]{BHS3}}, the second statement is \cite[Proposition 6.2.8]{KPX_cohomology}, and the third statement is {\cite[Corollary 6.2.9]{KPX_cohomology}}.
\end{proof}

\begin{defn}
Let $\cT$ denote the rigid analytic space of continuous characters of $K^\times$. 
Let $\cT_{\wreg}$ denote the open\footnote{Observe that any subset $\cZ$ of $\set{x^{-\bf k}, 
\left(N_{K|\bbQ_p}|N_{K|\bbQ_p}|_p\right)x^{\bf k}}{{\bf k} \in \bbN^{\Sigma_K}}$ is closed in $\cT$ since $\cZ$ is {\it locally finite}, i.e., for any affinoid subdomain $U=\Sp(A)$ of $\cT$, $U\cap \cZ$ is a finite set of rigid points, and hence closed. If $\cZ$ is finite, then there is nothing to show. Assume that $\cZ$ is infinite. Evaluation at $\pi_K\in K^\times$ of the {universal character} $\delta_A^{\univ}:K^\times \to \Gamma(U,\cO_U)^\times = A^\times$  gives rise to a unit $f:=\delta_A^{\univ}(\pi_K)\in A^\times$. By the maximum principle \cite[Theorem 3.1/15]{Bosch_lectures} applied to both $f^{-1}$ and $f$, there exists $C>0$ such that $C^{-1}\leq |f(z)|\leq C$ for every rigid point $z\in \Sp(A)$. 
However, since $\cZ$ is infinite, the $p$-adic valuations of its elements evaluated at $\pi_K$ are clearly unbounded.
Consequently, $U\cap \cZ$ is a finite subset of $\cZ$, and $\cZ$ is locally finite.} complement in $\cT$ to the points
$$
\set{
\left(N_{K|\bbQ_p}|N_{K|\bbQ_p}|_p\right)x^{\bf k}}{{\bf k} \in \bbN^{\Sigma_K}},
$$
and let $\cT_{\reg}$ denote the open complement in $\cT$ to the points 
$$
\set{x^{-\bf k}, 
\left(N_{K|\bbQ_p}|N_{K|\bbQ_p}|_p\right)x^{\bf k}}{{\bf k} \in \bbN^{\Sigma_K}}.
$$
For $n\geq 1$, we denote by $\cT^n_{\wreg}$ (resp. $\cT^n_{\reg}$) the open subspace of $\cT^n$ consisting of characters $\delta = (\delta_1,\dots,\delta_n)$ such that $\delta_i/\delta_j\in \cT_{\wreg}$ (resp. $\delta_i/\delta_j\in \cT_{\reg}$) for all $1\leq i\neq j\leq n$.
\end{defn}

\begin{lem}
\label{lem:rank_one_affinoid}
Let $\delta: K^\times \to A^\times$ be a continuous character.
For $z\in \Sp(A)$ corresponding to a maximal ideal $\frm_z \subset A$, let $\delta_z: K^\times \to A^\times \twoheadrightarrow k(z)^\times$ denote the mod-$\frm_z$ reduction of $\delta$.
\begin{enumerate}[(i)]
\item
For ${\bf k} \in \bbN^{\Sigma_K}$, if $\wt_{\sigma}(\delta_z) \notin \{1,\dots,k_\sigma\}$ for each $\sigma \in \Sigma_K$ and for all $z\in \Sp(A)$, then 
$$
\iota_{\bf k} : H^1(\delta) \to H^1(x^{-{\bf k}}\delta)
$$
is an isomorphism.
\item 
If $\delta \in \cT_{\wreg}(A)$, then $H^2(\delta)=0$, and $H^1(\delta)\otimes_A k(z) \cong H^1(\delta_z)$ for all $z\in \Sp(A)$.
\item 
If $\delta \in \cT_{\reg}(A)$, then $H^1(\delta)$ is locally free of rank $[K:\bbQ_p]$ over $A$ and $H^0(\delta)=H^2(\delta)=0$.
Moreover, $H^i(\delta)\otimes_A k(z) \cong H^i(\delta_z)$ for all $i\in \{0,1,2\}$ and $z\in \Sp(A)$.
\end{enumerate}
\end{lem}

\begin{proof}
By \cite[Corollary 6.3.3]{KPX_cohomology}, for any $(\varphi,\Gamma_K)$-module $D$ over $\cR_{K,A}$ and any $\sigma\in \Sigma_K$, $H^i(D)$ and $H^i(D/t_\sigma)$ are finitely generated $A$-modules; moreover, they are the cohomology of complexes of finite projective $A$-modules concentrated in degrees $[0,2]$.
Hence, as in \cite[\href{https://stacks.math.columbia.edu/tag/061Z}{Tag 061Z}]{stacks-project}, we have the base change spectral sequence
\begin{equation}
\label{eqn:spectral_sequence}
E_2^{j,-i}=\operatorname{Tor}_i^A\left(H^j(\heart), k(z)\right) \Rightarrow H^{j-i}(\heart \otimes_A k(z))
\end{equation}
for $\heart \in \{D,D/t_\sigma\}$.
As the $(\varphi,\Gamma_K)$-cohomology are concentrated in $[0,2]$, from \eqref{eqn:spectral_sequence} we see
\begin{equation}
\label{eqn:H2_base-change}
H^2(\heart) \otimes_A k(z) \xrightarrow{ \sim } H^2(\heart_z)
\end{equation}
is an isomorphism for every $z\in \Sp(A)$.
\begin{enumerate}[(i)]
\item
Since $\iota_{\bf k}$ is induced by $\cR_{K,A}(\delta) \hookrightarrow t^{-k}\cR_{K,A}(\delta) = \cR_{K,A}(x^{-\bf k}\delta)$, which factors as  
$$
\cR_{K,A} \subset t_{\sigma_1}^{-1}\cR_{K,A} \subset \dots \subset t_{\sigma_1}^{-k_{\sigma_1}}\cR_{K,A} 
\subset t_{\sigma_2}^{-1} t_{\sigma_1}^{-k_{\sigma_1}}\cR_{K,A} \subset \dots \subset  t_{\sigma_2}^{-k_{\sigma_2}}t_{\sigma_1}^{-k_{\sigma_1}}\cR_{K,A} 
\subset \dots \subset t^{-\bf k}\cR_{K,A},
$$
where we enumerate $\Sigma_K = \{\sigma_1,\sigma_2,\dots,\sigma_{[K:\bbQ_p]}\}$,
it suffices to prove that if $\wt_\sigma(\delta_z)\neq 1$ for all $z\in \Sp(A)$, then $\iota_\sigma: H^1(\delta) \to H^1(x_\sigma^{-1}\delta)$ is an isomorphism.

By the long exact sequence in cohomology attached to the exact sequence
$$
0\to \cR_{K,A}(\delta) \to t_\sigma^{-1}\cR_{K,A}(\delta) \to t_\sigma^{-1}\cR_{K,A}(\delta)/\cR_{K,A}(\delta) \to 0,
$$
it suffices to show the vanishing of $H^0$ and $H^1$ of $t_\sigma^{-1}\cR_{K,A}(\delta)/\cR_{K,A}(\delta)$.

For each $z\in \Sp(A)$, the fiber $t_\sigma^{-1}\cR_{K,k(z)}(\delta_z)/\cR_{K,k(z)}(\delta_z)$ has vanishing $H^2$ by \cite[Theorem 3.7(ii)]{Liu_cohomology} since it is a torsion $(\varphi,\Gamma_K)$-module.
By base change \eqref{eqn:H2_base-change} and that
$$
t_\sigma^{-1}\cR_{K,A}(\delta)/\cR_{K,A}(\delta) \cong \cR_{K,A}(\delta)/t_\sigma\cR_{K,A}(\delta)
$$ 
has coherent cohomology, we have $H^2(t_\sigma^{-1}\cR_{K,A}(\delta)/\cR_{K,A}(\delta)) = 0$ by Nakayama's lemma.
From the spectral sequence \eqref{eqn:spectral_sequence}, we see that for all $z\in \Sp(A)$,
$$
H^1(t_\sigma^{-1}\cR_{K,A}(\delta)/\cR_{K,A}(\delta)) \otimes_A k(z) \cong H^1(t_\sigma^{-1}\cR_{K,k(z)}(\delta_z)/\cR_{K,k(z)}(\delta_z)),
$$
which is zero since $\dim_{k(z)} H^0(t_\sigma^{-1}\cR_{K,k(z)}(\delta_z)/\cR_{K,k(z)}(\delta_z)) = 0$ by \cite[Lemma 2.16]{Nakamura_classification} and $\dim_{k(z)} H^0(t_\sigma^{-1}\cR_{K,k(z)}(\delta_z)/\cR_{K,k(z)}(\delta_z)) = \dim_{k(z)} H^1(t_\sigma^{-1}\cR_{K,k(z)}(\delta_z)/\cR_{K,k(z)}(\delta_z))$ by Euler-Poincare formula \cite[Theorem 4.3]{Liu_cohomology} and vanishing of $H^2$.
Again by Nakayama's lemma, we deduce the vanishing of $H^1$ and hence from the spectral sequence \eqref{eqn:spectral_sequence}
$$
H^0(t_\sigma^{-1}\cR_{K,A}(\delta)/\cR_{K,A}(\delta)) \otimes_A k(z) \cong H^0(t_\sigma^{-1}\cR_{K,k(z)}(\delta_z)/\cR_{K,k(z)}(\delta_z)),
$$
from which we deduce the vanishing of $H^0$ by Nakayama's lemma.
\item 
For $\delta\in \cT_{\wreg}(A)$, $\delta_z: K^\times \to k(z)^\times$ does not belong to $\set{
\left(N_{K|\bbQ_p}|N_{K|\bbQ_p}|_p\right)x^{\bf k}}{{\bf k} \in \bbN^{\Sigma_K}}$
for every $z\in \Sp(A)$, which implies that $H^2(\delta_z)=H^2(\delta)\otimes_A k(z) = 0$ for each $z\in \Sp(A)$ by \eqref{eqn:H2_base-change} and Lemma \ref{lem:rank_one_coh}(ii).
By Nakayama's lemma, $H^2(\delta)=0$.
Hence, \eqref{eqn:spectral_sequence} shows
$$
H^1(\cR_{K,A}(\delta)) \otimes_A k(z) \cong H^1(\cR_{K,k(z)}(\delta_z))
$$
for all $z\in \Sp(A)$.
\item 
This is \cite[Proposition 2.3]{HS_density}.
\qedhere
\end{enumerate}
\end{proof}

\subsubsection{}
Let $D$ be a $(\varphi,\Gamma_K)$-module of rank $n$ over $\cR_{K,A}$ for affinoid $A$ equipped with a filtration
$$
\Fil^\bullet(D) : \quad 
\Fil^0(D) = 0 \subsetneq \Fil^1(D) \subsetneq \dots \subsetneq \Fil^{n}(D) = D
$$
by saturated $(\varphi,\Gamma_K)$-submodules $\Fil^i(D)$ such that $\gr^i(\Fil^\bullet(D)):=\Fil^{i}(D)/\Fil^{i-1}(D)$ is locally free of rank $1$ over $\cR_{K,A}$ for $1\leq i\leq n$.
By Theorem \ref{thm:rank-one_classification}, there exists a unique continuous character $\delta_i: K^\times \to A^\times$ such that for the line bundle $\cL_i:= \Hom_{\varphi,\Gamma_K}(\cR_{K,A}(\delta_i), \gr^i(\Fil^\bullet(D))) $ the canonical map 
$$
\cR_{K,A}(\delta_i) \otimes_{A} \cL_i \longrightarrow \gr^i(\Fil^\bullet(D))
$$ 
is an isomorphism.
In this case, we say $D$ is {\bf trianguline}, the filtration $\Fil^\bullet(D)$ is a {\bf triangulation} of $D$, and $\delta = (\delta_i)_{1\leq i\leq n}$ is the {\bf parameter} of $(D,\Fil^\bullet(D))$.

If the triangulation on $D$ is clear, we write $D_i:= \Fil^i(D)$ and $D^i := D/\Fil^{n-i}(D)$ so that the ranks of the subobject $D_i$ and the quotient $D^i$ are always $i$, for $1\leq i\leq n$.

\subsubsection{}
We discuss various notions of ``non-split'' for $(\varphi,\Gamma_K)$-modules over $\cR_{K,A}$.

\begin{defn}
\label{defn:non-split}
\begin{enumerate}[(i)]
\item
Given a $(\varphi,\Gamma_K)$-module $D$ over $\cR_{K,A}$ with triangulation $\Fil^\bullet(D)$, if\footnote{In this definition, we do not need to write the line bundles $\cL_i\otimes_{A}k(z)$ because there are no nontrivial line bundles over the point $\Sp(k(z))$.}
\begin{equation}
\label{eqn:non-split}
0 \to \Fil^{i-1}(D_z) \to \Fil^i(D_z) \to \cR(\delta_{i,z}) \to 0
\end{equation}
are non-split as $(\varphi,\Gamma_K)$-modules for all $2\leq i\leq n$ and all rigid points $z$ of $\Sp(A)$, 
then $(D,\Fil^\bullet(D))$ is called {\bf non-split}.
\item 
Given a $(\varphi,\Gamma_K)$-module $D$ over $\cR_{K,A}$ with triangulation $\Fil^\bullet(D)$, if the exact sequences
\begin{equation}
\label{eqn:strong-non-split}
0 \to \cR(\delta_{i-1,z})\cong \Fil^{i-1}(D_z)/\Fil^{i-2}(D_z) \to \Fil^i(D_z)/\Fil^{i-2}(D_z) \to \cR(\delta_{i,z}) \to 0
\end{equation}
are non-split as $(\varphi,\Gamma_K)$-modules for all $2\leq i\leq n$ and all rigid points $z$ of $\Sp(A)$, then $(D,\Fil^\bullet(D))$ is called {\bf strongly non-split}.
\end{enumerate}
\end{defn}

If $(D,\Fil^\bullet(D))$ is strongly non-split, then it is non-split, cf. Remark \ref{rem:non-split}(i). 
But the converse is false when the rank of $D$ is at least 3, and a rank-3 counterexample is given in Remark \ref{rem:non-unique_triangulation}.

\begin{rem}
\label{rem:non-split}
\begin{enumerate}[(i)]
\item 
If $(D,\Fil^\bullet(D))$ is non-split \eqref{eqn:non-split} or strongly non-split  \eqref{eqn:strong-non-split}, then 
$$
0 \to \Fil^i(D) \to D \to D/\Fil^{i}(D) \to 0
$$
is a non-split extension of $(\varphi,\Gamma_K)$-modules, for all $1\leq i\leq n-1$.
\item 
If $D \in \Ext^1(D_2,D_1)$ is a non-split extension of trianguline $(\varphi,\Gamma_K)$-modules $(D_i,\Fil^\bullet(D_i))$ of rank $r_i$ for $i=1,2$, then $(D,\Fil^\bullet(D))$ need {\it not} be non-split in the sense of \eqref{eqn:non-split}, where $\Fil^\bullet(D)$ is the triangulation on $D$ induced by $\Fil^\bullet(D_1)$ and $\Fil^\bullet(D_2)$.
\end{enumerate}
\end{rem}

\begin{proof}
\begin{enumerate}[(i)]
\item 
The image of 
$$
0 \to \Fil^i(D) \to D \to D/\Fil^{i}(D) \to 0
$$ 
under the linear map
$$
\Ext^1(D/\Fil^{i}(D), \Fil^i(D)) \xrightarrow[\text{pullback}]{\Fil^{i+1}/\Fil^i(D)\subset D/\Fil^{i}(D)} \Ext^1(\delta_{i+1}, \Fil^i(D)) 
$$
is the extension $0 \to \Fil^{i}(D)\to \Fil^{i+1}(D) \to \cR(\delta_{i+1}) \to 0$ in \eqref{eqn:non-split}. 
Therefore, 
$$
0 \to \Fil^i(D) \to D \to D/\Fil^{i}(D) \to 0
$$
is non-split if $(D,\Fil^\bullet(D))$ is non-split.
Moreover, the image of this pullback extension
$$
0 \to \Fil^{i}(D)\to \Fil^{i+1}(D) \to \cR(\delta_{i+1}) \to 0
$$ under the linear map
$$
\Ext^1(\delta_{i+1}, \Fil^i(D))  \xrightarrow[\text{pushforward}]{\Fil^i(D)\twoheadrightarrow \Fil^i(D)/\Fil^{i-1}(D)} \Ext^1(\delta_{i+1},\delta_i)
$$
is the extension $0 \to \cR(\delta_{i}) \to \Fil^{i+1}(D)/\Fil^{i-1}(D) \to \cR(\delta_{i+1})\to 0$ in \eqref{eqn:strong-non-split}.
Consequently, ``strongly non-split" implies ``non-split", and
$$
0 \to \Fil^i(D) \to D \to D/\Fil^{i}(D) \to 0
$$
is non-split if $(D,\Fil^\bullet(D))$ is strongly non-split.
\item
We give a counterexample for $n=3$ and $\cR = \cR_{\bbQ_p,E}$.
Consider $D_1 = \cR(\delta_1)$ and $D_2 = (\cR(\delta_2) - \cR(\delta_3))$ the unique non-split extension with $\delta_1 = x, \delta_2 = |x|x$ and $\delta_3=1$.
From the short exact sequence 
$$
0 \to \cR(\delta_2) \to D_2 \to \cR(\delta_3) \to 0,
$$
we consider a segment of the associated long exact sequence 
$$
H^0(\delta_1\delta_2^{-1}) \to H^1(\delta_1\delta_3^{-1}) \to H^1(D_2^{\vee}(\delta_1)) \to H^1(\delta_1\delta_2^{-1}) \to H^2(\delta_1\delta_3^{-1}),
$$
and the maps between the $H^1$'s are pullbacks.
By Lemma \ref{lem:rank_one_coh}(ii), this sequence becomes $0 \to E \to E^2 \to E \to 0$, which in concrete terms means that if we take the non-split extension of $\cR(\delta_3)$ by $\cR(\delta_1)$, and pullback along $D_2\twoheadrightarrow \cR(\delta_3)$, then we get a non-split extension $D$ of $D_2$ by $D_1$ whose pullback along $\cR(\delta_2) \hookrightarrow D_2$ is a split extension of $\cR(\delta_2)$ by $\cR(\delta_1)$, which is $\Fil^2(D)$.
Hence, although $D$ is non-split as an extension of $D_2$ by $D_1$, it is not non-split in the sense of \eqref{eqn:non-split}.
\qedhere
\end{enumerate}
\end{proof}

\section{Triangulations on \text{$(\varphi,\Gamma_K)$}-modules}
\label{section:triangulations}
\begin{thm}
\label{thm:p_k-invertible}
Let $(D,\Fil^\bullet(D))$ be trianguline with parameters $(\delta_1,\dots,\delta_n):(K^\times)^n \to A^\times$ over $\cR_{K,A}$.
For any fixed ${\bf k} \in \bbN^{\Sigma_K}$ and $i\in \{1,\dots,n\}$, 
suppose that
$$
\wt_{\sigma}(\delta_{j,z}/\delta_{\ell,z}) \notin \{1,\dots,k_\sigma\}
$$ 
for all $1\leq j \leq n-i$, $n-i+1\leq \ell \leq n$, for all $\sigma \in \Sigma_K$ and for all $z\in \Sp(A)$.
Then
$$
p_{\bf k}: \Ext^1(D/\Fil^{n-i}(D), \Fil^{n-i}(D)) \longrightarrow \Ext^1(t^{\bf k}(D/\Fil^{n-i}(D)), \Fil^{n-i}(D))
$$
is an isomorphism.
\end{thm}
\begin{proof}
By induction on ${\bf k}$ it suffices to prove the statement for a fixed $\sigma\in \Sigma_K$ with $k_\sigma = 1$ and $k_\tau = 0$ for all $\tau\in \Sigma_K\setminus \{\sigma\}$, under the assumption that $\wt_{\sigma}(\delta_{j,z}/\delta_{\ell,z}) \neq 1$ for $j\leq n-i<\ell$.

Recall our notation that $D_{n-i}:=\Fil^{n-i}(D)$ and $D^i:= D/\Fil^{n-i}(D)$.
By Theorem \ref{thm:rank-one_classification}, $D$ is a successive extension of rank 1 $(\varphi,\Gamma_K)$-modules $\gr_m(\Fil^\bullet(D)) \cong \cR_{K,A}(\delta_m)\hat\otimes_A \cL_m$ for some line bundles $\cL_m$ over $A$. 
There is a finite admissible cover $\{\Sp(A_\nu)\}_{\nu=1}^r$ of $\Sp(A)$ trivializing them.
Since $H^i$ commutes with flat base change by \cite[Theorem 4.4.3(2)]{KPX_cohomology}, passing to $\Sp(A_\nu)$ we may assume that the line bundles $\cL_m$ are all trivial.

After the reductions above, we show that 
$$
p_{i,\sigma}: \Ext^1(D^i, D_{n-i}) \cong H^1({D^{i}}^\vee\otimes D_{n-i}) \to \Ext^1(t_\sigma D^i, D_{n-i}) \cong H^1(t_\sigma^{-1}({D^{i}}^\vee\otimes D_{n-i}))
$$
is bijective, where we have by assumption that 
$$
D_{n-i} = (\cR_{K,A}(\delta_1)-\dots-\cR_{K,A}(\delta_{n-i})) \quad \text{ and } \quad D^i = (\cR_{K,A}(\delta_{n-i+1})-\dots-\cR_{K,A}(\delta_{n})).
$$
Thus, putting $M:= {D^{i}}^\vee\otimes D_{n-i}$, we see that $M$ is trianguline over $\cR_{K,A}$ with parameters 
$$
\set{\delta_j\delta_\ell^{-1}}{1\leq j\leq n-i \text{ and } n-i+1\leq  \ell \leq n}.
$$
Note that $\wt_\sigma(\delta_j\delta_\ell^{-1})\neq 1$ for these parameters by the assumption.
By the short exact sequence
$$
0 \to M \to t_\sigma^{-1}M \to t_\sigma^{-1}M/M \to 0,
$$
it is enough to establish the vanishing of $H^i(t_\sigma^{-1}M/M)$ as coherent sheaves for $i=0,1,2$.

This follows from a dévissage argument on the above triangulation on $M$, using the proof of Lemma \ref{lem:rank_one_affinoid}(i).
Indeed, we know that $M = (\cR_{K,A}(\eta_1)-\dots-\cR_{K,A}(\eta_m))$ of parameters $\eta_s$ whose weights are not $1$ at all $z\in \Sp(A)$, so that
$$
t_\sigma^{-1}M/M \cong M/t_\sigma = (\cR_{K,A}(\eta_1)/t_\sigma -\dots-\cR_{K,A}(\eta_m)/t_\sigma).
$$
We induct on the rank of $M$ to show the vanishing of cohomologies of $M/t_\sigma$.
The base case where $m=1$ was shown in the proof of Lemma \ref{lem:rank_one_affinoid}(i).
For $m\geq 2$, we have
$$
\dots \to H^i(\cR_{K,A}(\eta_1)/t_\sigma) \to H^i(M/t_\sigma) \to H^i((M/\cR_{K,A}(\eta_1))/t_\sigma) \to \dots ,
$$
from which we know $H^i(M/t_\sigma)=0$ by the induction hypothesis.
\end{proof}

\subsection{Very generic case}

\begin{defn}
Let $(D,\Fil^\bullet(D))$ be a trianguline ($\varphi, \Gamma_K$)-module over $\cR_{K,E}$ with parameter $\delta = (\delta_1,\dots,\delta_n)$ and Sen weights $(h_{i,\sigma})_{i,\sigma}\in (E^n)^{\Sigma_K}$, where $h_{i,\sigma}:=\wt_\sigma(\delta_i)$.
We say that $D$ is {\bf very generic} if for each $\sigma\in \Sigma_K$, $h_{i,\sigma}-h_{j,\sigma}\notin \bbZ$ for all $i\neq j$.
\end{defn}

\begin{defn}
Let $\cT_{\circ}^n$ be the open subspace of $\cT^n$ such that for $A\in \Aff_E$,  $\cT_{\circ}^n(A)$ is the set of all the continuous $A^\times$-valued characters of $(K^\times)^n$
$$
\delta = (\delta_1,\dots, \delta_n):(K^\times)^n\to A^\times
$$
satisfying
$\wt_\sigma(\delta_{i,z}/\delta_{j,z}) \notin \bbZ_{\geq 1}$
for all $1\leq i<j \leq n$, all $\sigma\in \Sigma_K$, and all $z\in \Sp(A)$.
\end{defn}

\begin{prop}
\label{prop:unique_saturated_rank-j_subobject}
For a trianguline $(\varphi,\Gamma_K)$-module $(D,\Fil^\bullet(D))$ of parameters $\delta=(\delta_i)_{1\leq i\leq n} \in T^n_{\circ}(E)$ over $\cR_{K,E}$ and a fixed index $i\in \{1,\dots, n\}$, suppose 
\begin{enumerate}[(a).]
\item 
$\Fil^{n-i}(D)$ with the induced triangulation is strongly non-split, and 
\item 
$D/\Fil^{n-i-1}(D)$ with the induced triangulation is non-split.
\end{enumerate}
Then, for each $1\leq j\leq n-i$, $\Fil^j(D)$ is the unique saturated  $(\varphi,\Gamma_K)$-submodule of $D$ of rank $j$.
\end{prop}
\begin{proof}
Note that, similarly to Remark \ref{rem:non-split}, (a) and (b) imply that $(D,\Fil^\bullet(D))$ is non-split.

We proceed by induction.
Let $D'$ be a saturated ($\varphi,\Gamma_K$)-submodule of $D$.
We claim that $D'$ contains $\Fil^1(D)$.
Indeed, let $m$ be the largest integer such that $\Fil^{m-1}(D)\cap D'=0$.
Then 
$D'\cap \Fil^m(D)$ is an $\cR_{K,E}$-submodule of $D'$ stable under the $(\varphi,\Gamma_K)$-action such that the quotient $D'/(D'\cap \Fil^m(D))$ is torsion-free as it injects into $D/\Fil^m(D)$.
As 
$$
\cR_{K,E} \simeq \prod_{\tau:K_0\hookrightarrow E} \cR_{K,\bbQ_p} \otimes_{K_0,\tau} E
$$
is a product of B\'ezout domains on which $\varphi$ acts transitively, it follows from the torsion-freeness that $D'\cap \Fil^m(D)$ is free over $\cR_{K,E}$ (cf. \cite[\S2.6]{Bergdall}) and hence is a $(\varphi,\Gamma_K)$-module.

Now if $m> 1$, then $D'\cap \Fil^m(D)$ is a nonzero $(\varphi,\Gamma_K)$-submodule of $\Fil^m(D)$ injecting to the quotient $\Fil^m(D)/\Fil^{m-1}(D) \cong \cR_{K,E}(\delta_{m})$, which is of the form $t^{\bf k}\cR(\delta_m) = \cR_{K,E}(x^{\bf k}\delta_m)$ for some ${\bf k} \in \bbN^{\Sigma_K}$ by Lemma \ref{lem:rank_one_coh}(iii).
Then, $\Fil^m(D)$ contains the split $(\varphi,\Gamma_K)$-submodule 
$$
\Fil^{m-1}(D) \oplus (D'\cap \Fil^m(D)) \cong \Fil^{m-1}(D)\oplus t^{\bf k}\cR_{K,E}(\delta_{m}),
$$
which is the image of $\Fil^m(D)$ under the pullback map
$$
p_{\bf k}: \Ext^1(\delta_m, \Fil^{m-1}(D)) \to \Ext^1(t^{\bf k}\delta_m, \Fil^{m-1}(D)),
$$
which is bijective by Theorem \ref{thm:p_k-invertible}, contrary to our assumptions (a) and (b).
Hence, $m=1$. 
Since $D'\cap \Fil^1(D) = D'\cap \cR_{K,E}(\delta_1)$ is nonzero and saturated in $\cR_{K,E}(\delta_1)$, we conclude that $D'\cap \Fil^1(D) = \cR_{K,E}(\delta_1) = \Fil^1(D)$.

In particular, $\Fil^1(D)$ is the unique saturated $(\varphi,\Gamma_K)$-submodule of rank $1$ of $D$.
If $n-i>1$, we pass to the quotient by $\Fil^1(D)$ and claim that every saturated $(\varphi,\Gamma_K)$-submodule of $D/\Fil^1(D)$ contains $\Fil^2(D)/\Fil^1(D)$. 
By (a) and (b), $D/\Fil^1(D)$ with its induced filtration is non-split and satisfies the analogous (a) and (b), so the argument above works.
We can iterate this argument until $n-i$, hence the conclusion.
\end{proof}

\begin{cor}
\label{cor:generic-unique}
Let $(D,\Fil^\bullet(D))$ be strongly non-split of parameters $\delta\in \cT_\circ^n(E)$ over $\cR_{K,E}$.
Then, any saturated ($\varphi,\Gamma_K$)-submodule $D'$ of $D$ satisfies that $D'=\Fil^i(D)$ for some $0\leq i\leq n$.
In particular, the given triangulation $\Fil^\bullet(D)$ is the unique triangulation of $D$.
\end{cor}
\begin{proof}
This follows from Proposition \ref{prop:unique_saturated_rank-j_subobject} by taking $i=1$.
\end{proof}

\begin{rem}
\label{rem:non-unique_triangulation}
One cannot expect the uniqueness result for $(\varphi,\Gamma_K)$-module $(D,\Fil^\bullet(D))$ that is non-split with parameters in $\cT_\circ^n(E)$ without the strongly non-split assumption. 
Indeed, for regular characters $(\delta_1,\delta_2,\delta_3) \in \cT_{\reg}^3(E) \cap \cT_{\circ}^3(E)$, we have an exact sequence
$$
0 \to \Ext^1(\delta_3, \delta_1) = E \to \Ext^1(\delta_3,(\delta_1-\delta_2)) = E^2 \to \Ext^1(\delta_3,\delta_2) = E \to 0,
$$
where $(\delta_1-\delta_2)$ denotes the non-split extension of $\cR(\delta_2)$ by $\cR(\delta_1)$.
The image of any nonzero class of $\Ext^1(\delta_3, \delta_1)$ in $\Ext^1(\delta_3,(\delta_1-\delta_2))$ is a trianguline $(\varphi,\Gamma_K)$-module $(D,\Fil^\bullet(D),(\delta_1,\delta_2,\delta_3))$ that is non-split in the sense of \eqref{eqn:non-split}, but $D/\Fil^{1}(D)$ is a split extension of $\delta_2$ by $\delta_3$. 
Therefore, $D$ has another triangulation whose parameter is $(\delta_1,\delta_3,\delta_2)$.
\end{rem}

\subsection{Crystabelline case}

We discuss triangulations on a class of trianguline $(\varphi,\Gamma_K)$-modules that is closer to $p$-adic Hodge theory.
Recall that we have Fontaine's functors $D_{\cris},D_{\dR}$ defined for $(\varphi,\Gamma_K)$-modules as well, cf. \cite[Definition 2.5]{HS_density} or \cite[Definition 2.3]{Nakamura_B-pairs}.

\begin{defn}
\label{defn:cristabelline}
A $(\varphi,\Gamma_K)$-module $D$ over $\cR_{K, E}$ is {\bf crystabelline} if there exists a finite abelian extension $L/K$ such that $\cR_{L,E}\otimes_{\cR_{K,E}} D$ is crystalline, i.e., the $(L_0\otimes_{\bbQ_p}E)$-rank of 
$$
D_{\cris}^L (D) := (\cR_{L,E}\otimes_{\cR_{K,E}} D)[1/t]^{\Gamma_{L}}
$$
is equal to $\rank_{\cR_{K,E}}(D)$.
\end{defn}

As we describe below, one can understand this class of $(\varphi,\Gamma_K)$-modules in terms of simpler semilinear algebra and descent data.

\begin{defn}[{\cite[Definition 2.4]{Nakamura_B-pairs}}]
\label{defn:E-filtered_phi_GK_mod}
Let $L$ be a finite Galois extension of $K$ with Galois group $G(L/K)$. 
We say that $D$ is an {\bf $E$-filtered ($\varphi, G(L / K)$)-module over $K$} if
\begin{enumerate}[(i)]
\item 
$D$ is a finite free $(L_0\otimes_{\bbQ_p}E)$-module with a Frobenius semilinear operator $\varphi:D\xrightarrow{\simeq}D$, and a semilinear action by $G(L/K)$ that commutes with $\varphi$.
\item 
$D_L:= (L\otimes_{L_0}D)$ has a separated and exhaustive descending filtration $(\Fil^i(D_L))_{i\in\bbZ}$ by $G(L/K)$-stable $(L\otimes_{\bbQ_p}E)$-submodules $\Fil^i(D_L)$.
\end{enumerate}
An $(L_0\otimes_{\bbQ_p}E)$-module $D$ satisfying (i) will be called a {\bf ($\varphi, G(L / K)$)-module over $K$}.
\end{defn}

\begin{rem}
\label{rem:Hodge_fil}
\begin{enumerate}[(i)]
\item
For a ($\varphi, G(L / K)$)-module $D$ as in Definition \ref{defn:E-filtered_phi_GK_mod}, we have 
$$
L\otimes_{K} (D_L)^{G(L/K)} \xrightarrow{\sim} D_L
$$ 
by Galois descent, and we write $D_K:=(D_L)^{G(L/K)}$.
Then, a filtration by $G(L / K)$-stable $L \otimes_{\bbQ_p} E$-submodules on $D_L$ gives rise by taking invariants $(-)^{G(L/K)}$ to a filtration by $K \otimes_{\mathbb{Q}_p} E$-submodules on $D_K$. 
Conversely, a filtration by $K \otimes_{\mathbb{Q}_p} E$-submodules on $D_K$ gives rise by base change $L\otimes_K(-)$ to a filtration by $G(L / K)$-stable $L \otimes_{\bbQ_p} E$-submodules on $D_L$.
Via \eqref{eqn:Galois_decoupling}, $D_K \cong \prod_{\sigma \in \Sigma_K} D_{K,\sigma}$ decomposes into a product of $E$-vector spaces $D_{K,\sigma}$, and any filtration $\mathrm{Fil}^r D_K=\prod_{\sigma \in \Sigma_K} \mathrm{Fil}^i D_{\mathrm{dR}, \sigma}$ also decomposes into a product of filtrations by $E$-vector spaces on each $\sigma$-component $D_{K,\sigma}$.
\item
For a crystabelline $(\varphi,\Gamma_K)$-module $D$ over $\cR_{K,E}$ that becomes crystalline over $L$, $D_{\cris}^L(D)$ has the structure of an $E$-filtered ($\varphi, G(L / K)$)-module whose $G(L/K)$-stable filtration on $L\otimes_{L_0} D_{\cris}^L(D)$ is induced by the comparison isomorphism
$$
L\otimes_{L_0} D_{\cris}^L(D)\xrightarrow{\simeq} L\otimes_{K}D_{\dR}(D),
$$
where $(L\otimes_{L_0} D_{\cris}^L(D))^{G(L/K)}\cong D_{\dR}(D)$ has the Hodge filtration $\{\Fil^i(D_{\dR}(D))\}_{i\in \bbZ}$ by $K\otimes_{\bbQ_p} E$-submodules, cf. \cite[Definition 2.5]{HS_density}.
Via \eqref{eqn:Galois_decoupling}, we have 
$$
D_{\dR}(D) =\prod_{\sigma\in\Sigma_K}D_\dR(D)_\sigma
\quad \text{ and } \quad
\Fil^i(D_{\dR}(D)) = \prod_{\sigma\in \Sigma_K} \Fil^i(D_{\dR}(D)_{\sigma}),
$$
such that the jumps in the filtration $\prod_{\sigma\in \Sigma_K} \Fil^i(D_{\dR}(D)_{\sigma})$ by $E$-vector spaces are given by $(-1)\cdot(\text{$\sigma$-Sen weights of $D$})$, cf. \cite[Definition 6.2.5]{KPX_cohomology}.
\end{enumerate}  
\end{rem}

\begin{thm}
\label{thm:D_cris_WD_equivalence}
Let $L$ be a finite Galois extension of $K$.
\begin{enumerate}[(i)]
\item 
The functor $D\mapsto D_{\cris}^L(D)$ induces a $\otimes$-equivalence of categories from $(\varphi,\Gamma_K)$-modules over $\cR_{K,E}$ that become crystalline over $L$ to $E$-filtered $(\varphi,G(L/K))$-modules over $K$. 
The saturated submodules of $D$ correspond to $(\varphi, G(L/K))$-submodules of $D_{\cris}^L(D)$ with their filtrations induced by the {\it Hodge filtration} on $D_{\cris}^L(D)$.
\item 
Let $L_0$ be the maximal unramified subfield of $L$, and assume that $E$ contains all $\bbQ_p$-algebra embeddings of $L_0$ into $\overline{\bbQ_p}$.
For any free $(L_0\otimes_{\mathbb{Q}_p} E)$-module $D$ of rank $n$, we canonically have $D = \prod_{\sigma:L_0\hookrightarrow E} D_{\sigma}$, where each $D_\sigma$ is an $E$-vector space of rank $n$.
Fix an embedding $\sigma_0:L_0\hookrightarrow E$.
There is an equivalence of categories from the category of $(\varphi,G(L/K))$-modules $D$ over $K$ to the category of representations $(r,V)$ of the Weil group $W_K$ of $K$ on finite dimensional $E$-vector spaces such that $r$ is unramified when restricted to the Weil group $W_L$ of $L$.
On the underlying modules, this equivalence is given by $D\mapsto V:= D_{\sigma_0}$.
The action $r(w)$, for $w\in W_K$, on $V$ is the $\sigma_0$-component of the $(L_0\otimes_{\bbQ_p} E)$-linear map $r(w):=\overline{w}\circ \varphi^{-\alpha(w)}:D\to D$, where $\alpha(w)\in f\bbZ$ is the unique integer such that the image of $w$ in $G(\overline{\bbF_p}/\bbF_p)$ is given by $x\mapsto x^{p^{\alpha(w)}}$, and $\overline{w}\in G(L/K)$ is the restriction of $w$ to $L$.

\end{enumerate}
\end{thm}
\begin{proof}
(i) is {\cite[Theorem 2.5]{Nakamura_B-pairs}} and \cite[Proposition 3.3(a)]{Bergdall}. Also see \cite{Berger_filtres}.

(ii) is the special case of \cite[Proposition 4.1]{BS_functoriality} when the monodromy operator $N=0$.
\end{proof}

\begin{defn}
\label{defn:cris-generic}
\begin{enumerate}[(i)]
\item
A crystabelline $(\varphi,\Gamma_K)$-module $D$ is {\bf generic} if the $W_K$-representation associated to the $(\varphi,G(L/K))$-module $D_{\cris}^L(D)$ is a direct sum of characters $\phi_1,\dots, \phi_n$ of $W_{K}$ that are trivial on $I_L$ such that if $\phi_1,\dots, \phi_n$ are viewed as smooth characters of $K^\times$ via the local reciprocity map, then $\phi_i/\phi_j \notin \{1, |\cdot|_K^{-1}, |\cdot|_K\}$
for all $i\neq j$, cf. {\cite[\S2.1]{Ding_crystabelline}}.
\item For a crystabelline generic $(\varphi,\Gamma_K)$-module $D$ as in (i), a {\bf refinement} of $D$ is an ordering of the $n$ characters $\phi_1,\dots,\phi_n$ appearing in the associated Weil representation.
If we fix an ordering $\phi_1,\dots,\phi_n$, then all the refinements are indexed by elements of $S_n$, where $w\in S_n$ corresponds to the ordering $\phi_{w(1)},\dots, \phi_{w(n)}$.
\item
To each refinement $w\in S_n$ of $D$ as in (ii), we associate a triangulation $\Fil^\bullet_w(D)$ on $D$ by requiring that $\Fil^i_w(D)$ corresponds to the $(\varphi,G(L/K))$-submodule of $D_{\cris}^L(D)$, equipped with the induced filtration, whose associated Weil representation is the subrepresentation $\bigoplus_{1\leq j\leq i} \phi_{w(j)} $ of $\bigoplus_{j=1}^n \phi_j$ under the two equivalences recalled in Theorem \ref{thm:D_cris_WD_equivalence}, for all $i$.
\end{enumerate}
\end{defn}

\begin{defn}
\label{defn:ncr}
Let ${\bf h} = (h_{1,\sigma}\geq \dots \geq h_{n,\sigma})_{\sigma} \in (\bbZ^n)^{\Sigma_K}$ be the Sen weights of a crystabelline generic $(\varphi,\Gamma_K)$-module $D$ of rank $n$.
\begin{enumerate}[(i)]
\item 
We say ${\bf h}$ is  {\bf regular} if $h_{1,\sigma}>\dots>h_{n,\sigma}$ for each $\sigma\in \Sigma_K$.
\item 
A refinement of the crystabelline generic $(\varphi,\Gamma_K)$-module $D$
$$
w(\phi) = (\phi_{w(1)},\dots,\phi_{w(n)}), \text{ for $w\in S_n$},
$$ is {\bf non-critical} if the $\sigma$-Sen weights of $\Fil_w^i(D)$ are $(h_{1,\sigma}\geq \dots \geq h_{i,\sigma})$ for all $1\leq i\leq n$ and $\sigma\in \Sigma_K$.
We call a refinement of $D$ {\bf critical} if it is not non-critical.
Then, $D$ is called {\bf non-critical} if all its refinements are non-critical, and is called {\bf critical} if at least one of its refinements is critical.
\item 
Denote by $\PGnc$ the class of all crystabelline non-critical $(\varphi,\Gamma_K)$-modules of rank $n$ over $\cR_{K,E}$ with associated $W_K$-representation $\phi = \bigoplus_{i=1}^n \phi_i$ and Sen weights $\mathbf{h}$.
\end{enumerate}
\end{defn}

\begin{rem}
\label{rem:induced_fil}
In the situation of Definition \ref{defn:cris-generic}(iii) and Remark \ref{rem:Hodge_fil}, under the identifications $(L\otimes_{L_0} D_\cris^L(D))^{G(L/K)} = D_\dR(D)$ and $D_\dR(D)=\prod_{\sigma\in\Sigma_K} D_{\dR}(D)_\sigma$, we recall the relationship between refinements of $D_\cris^L(D)$ and the jumps in the induced Hodge filtration $\Fil^\bullet(D_{\dR}(\Fil_w^i(D))_\sigma)$ for $1\leq i\leq n$ and $\sigma\in \Sigma_K$.
Let $\pr_\sigma:D_\dR(D)\to D_{\dR}(D)_\sigma$ denote the projection onto the $\sigma$-component.
For each factor $\phi_j$ of the $W_K$-representation $\bigoplus_{j=1}^n \phi_j$, let $M_j$ be the corresponding $(\varphi,G(L/K))$-submodule of $D_\cris^L(D)$ of rank 1.
Via the construction $\left(\mathrm{pr}_\sigma\right) \circ\left(L \otimes_{L_0}(-)\right)^{G(L / K)}$, $M_j$ gives rise to a 1-dimensional $E$-subspace of $D_\dR(D)_\sigma$, and we choose a basis $e_{j,\sigma}$.
For each refinement $w\in S_n$, let $\cF^w_\bullet := (0\subset \cF^w_{1} \subset \dots \subset \cF^{w}_{n} = D_\cris^L(D))$ be the complete $(\varphi,G(L/K))$-stable flag on $D_\cris^L(D)$ corresponding to the triangulation $\Fil_w^\bullet(D)$.
Then, $\cF^w_{i}=\oplus_{j=1}^i M_{w(j)}$, and via the construction $\left(\mathrm{pr}_\sigma\right) \circ\left(L \otimes_{L_0}(-)\right)^{G(L / K)}$, $\cF^w_{i}$ gives rise to an $E$-subspace $\mathcal{F}^w_{i,\sigma} \subset D_\dR(D)_\sigma$ of dimension $i$ with a basis $\left(e_{w(1), \sigma}, \ldots, e_{w(i), \sigma}\right)$.
We also have the flag $\Fil^\bullet(D_{\dR}(D)_\sigma)$ coming from the Hodge filtration.
Assume moreover that ${\bf h}$ is regular.
The ``relative position'' of $\Fil^\bullet(D_{\dR}(D)_\sigma)$ and $\cF^w_{\bullet,\sigma}$ is given by the element $ \frw_\sigma \in S_n$ satisfying the condition that if the jumps of the $\sigma$-Hodge filtration on $\cF^w_{n,\sigma} = D_\dR(D)_\sigma$ are given by $\{-h_{1,\sigma} < \dots < -h_{n,\sigma}\}$, then the induced $\sigma$-Hodge filtration on $\cF^w_i$ has jumps at $\{-h_{\frw_\sigma^{-1}(j),\sigma}|1\leq j\leq i\}$, in which case the induced ordering of $\sigma$-Sen weights by the triangulation $\Fil_w^\bullet(D)$ is $\{h_{\frw_{\sigma}^{-1}(1),\sigma},\dots, h_{\frw_{\sigma}^{-1}(n),\sigma}\}$.
We say that the flags $\Fil^\bullet(D_{\dR}(D))$ and $\cF^w_{\bullet}$ are {\it of relative position $\frw=(\frw_\sigma)_\sigma\in S_n^{\Sigma_K}$}.
Then, $D$ is non-critical if and only if $\Fil^\bullet(D_{\dR}(D))$ and $\cF^w_{\bullet}$ are of relative position $\id \in S_n^{\Sigma_K}$ for all $w\in S_n$.
\end{rem}

\begin{prop}
\label{prop:crys_generic_triang}
Let $D$ be a crystabelline generic $(\varphi,\Gamma_K)$-module of rank $n$ over $\cR_{K,E}$ that is crystalline over a finite abelian extension $L/K$, with  Weil representation $\phi = \bigoplus_{i=1}^n \phi_i$ and Sen weights ${\bf h}= (h_{1,\sigma}\geq \dots \geq h_{n,\sigma})_{\sigma\in \Sigma_K}$.
Then, $D$ is trianguline and has $n!$ triangulations
$$
\Fil_w^\bullet(D): \quad 0 \subsetneq \Fil_w^1(D) \subsetneq \dots \subsetneq \Fil_w^n(D) = D,
$$
which are indexed by refinements $w\in S_n$, and are of parameters $\delta_w = (\delta_{w,1},\dots, \delta_{w,n})$ with 
$$
\delta_{w,i} = \left(\prod_{\sigma\in\Sigma_K} x_{\sigma}^{k_{i,\sigma}^w}\right)\phi_{w(i)}
$$ 
for some $k_{i,\sigma}^w \in \bbZ$ such that $\{h_{1,\sigma},\dots, h_{n,\sigma}\} = \{k_{1,\sigma}^w,\dots, k_{n,\sigma}^w\}$ as multisets, for all $\sigma$.
Moreover,  $D$ is non-critical if and only if $k_{i,\sigma}^w = h_{i,\sigma}$, for all $1\leq i\leq n$, $w\in S_n$, and $\sigma\in \Sigma_K$.
\end{prop}
\begin{proof}
The last assertion would follow from the previous assertions and the definition of non-critical.
The genericity of $D$ implies that there are $n!$ complete flags consisting of subobjects in the Weil representation $\bigoplus_{i=1}^n \phi_i$ attached to $D$.
By the equivalences in Theorem \ref{thm:D_cris_WD_equivalence}, the flags give rise to $n!$ distinct triangulations on $D$, indexed by permutations $w\in S_n$ of the set $\{\phi_1,\dots,\phi_n\}$.
It remains to show that the triangulation parameters are of the forms given in the statement.
Since the property of being crystabelline is stable under taking subquotients, by considering $\Fil_w^{i}(D)/\Fil_w^{i-1}(D)$ we are reduced to showing the following claim on crystabelline $(\varphi,\Gamma_K)$-modules of rank $1$: if $\cR_{K,E}(\delta)$ is a crystabelline $(\varphi,\Gamma_K)$-module of rank $1$ that becomes crystalline over $L$, then $\delta = \widetilde{\delta}\prod_{\sigma\in \Sigma_K}(x_\sigma^{k_\sigma})$ for a smooth character $\widetilde{\delta}:K^\times \to E^\times$ and $k_\sigma \in \bbZ$, and the $W_K$-representation $\phi$ attached to the $(\varphi,G(L/K))$-module $D_{\cris}^L(\cR_{K,E}(\delta))$ is $\widetilde{\delta}\circ \rec_K^{-1}$.

By construction, the $W_K$-representation $\phi$ on $E$-vector spaces is independent of the choice of $L$ over which $R_{K,E}(\delta)$ becomes crystalline.
Observe that any crystabelline $(\varphi,\Gamma_K)$-module $M$ becomes crystalline over a totally ramified abelian extension $L/K$.
For this observation, we can prove it for the $B$-pair $W = (W_e,W_{\dR}^+)$ associated to $M$ in \cite[\S1.3]{Nakamura_classification}, since $D_{\cris}^L(M) = D_{\cris}^L(W)$ by \cite[Proposition 2.3.4]{Berger_B-pairs}.
By local class field theory, there are $m\in \bbZ_{\geq 1}$ and a finite unramified extension $L'/K$ such that $L\subset K_m.L'$, where $K_m$ is the $m$-th Lubin-Tate extension of $K$.
We have the short exact sequence of Galois groups
$$
1 \to G_{K_m.L'} \to G_{K_m} \to G(K_m.L'/K_m) \to 1.
$$
Since $L \subset K_m.L'$, it follows that $G_{K_m.L'} \subset G_L$, and $M$ becomes crystalline over $K_m.L'$.
Then,
$$
D_{\cris}^{K_m.L'}(W) := (B_{\max} \otimes_{B_e} W_e)^{G_{K_m.L'}} \cong D_{\cris}^{K_m.L'}(M)
$$
is a semilinear $G(K_m.L'/K_m)$-module over $(K_{m}.L')_0 = L_0'$.
By Galois descent, we have
$$
D_{\cris}^{K_m}(W)\otimes_{(K_m)_0} (K_m.L')_0  = (B_{\max}\otimes_{B_e} W_e)^{G_{K_m}} \otimes_{K_0} L_0' \xrightarrow{\simeq} (B_{\max} \otimes_{B_e} W_e)^{G_{K_m.L'}} = D_{\cris}^{K_m.L'}(W),
$$
which shows that $D$ becomes crystalline over $K_m$.
Hence, from now on we may take $L=K_m$ for $m\gg 0$, which is totally ramified over $K$, in our proof of the claim.

Since $\cR_{K,E}(\delta)$ is crystabelline, $\delta = \widetilde{\delta}\prod_{\sigma\in \Sigma_K}(x_\sigma^{k_\sigma})$ for some smooth character $\widetilde{\delta}:K^\times \to E^\times$ and integers $k_\sigma \in \bbZ$ by \cite[Lemma 4.1]{Nakamura_classification}.
We need to show that the Weil representation $\phi$ attached to $D_{\cris}^{K_m}(R_{K,E}(\delta))$ is $\widetilde{\delta}\circ \rec_K^{-1}$.
Let $G_m:=G(K_m/K)$ and write $\widetilde{\delta}=\widetilde{\delta}^{\unr}.\widetilde{\delta}^{\wt}$ so that $\widetilde{\delta}^{\unr}(\pi_K)=\widetilde{\delta}(\pi_K)$ and $\widetilde{\delta}^{\wt}|_{\cO_K^\times} = \widetilde{\delta}|_{\cO_K^\times}$.
We first compute the $(\varphi, G_m)$-module $D_{\cris}^{K_m}(R_{K,E}(\delta))$ by Theorem \ref{thm:D_cris_WD_equivalence}(i), and then compute its Weil representation by Theorem \ref{thm:D_cris_WD_equivalence}(ii).
We have
\begin{align*}
D_{\cris}^{K_m}(\cR_{K,E}(\delta)) = D_{\cris}^{K_m}(\cR_{K,E}(\widetilde{\delta}^{\wt})) \otimes_{K_0\otimes_{\bbQ_p} E} D_{\cris}^{K_m}(\cR_{K,E}(\widetilde{\delta}^{\unr}\prod_{\sigma\in\Sigma_K} x_\sigma^{k_{\sigma}})).
\end{align*}
By $(K_m)_0=K_0$ and applying to $\cR_{K,E}(\widetilde{\delta}^{\unr}\prod_{\sigma\in\Sigma_K} x_\sigma^{k_{\sigma}})$  \cite[Example 6.2.6(3)]{KPX_cohomology}, we have
$$
D_{\cris}^{K_m}(\cR_{K,E}(\widetilde{\delta}^{\unr}\prod_{\sigma\in\Sigma_K} x_\sigma^{k_{\sigma}})) = (K_m)_0 \otimes_{K_0} D_{\cris}(\cR_{K,E}(\widetilde{\delta}^{\unr}\prod_{\sigma\in\Sigma_K} x_\sigma^{k_{\sigma}})) ,
$$
where $G(K_m/K)$ acts trivially, and as a $\varphi$-module, it is the free $(K_0\otimes_{\bbQ_p} E)$-module $D_{f,\widetilde{\delta}(\pi_K)}$ of rank 1 equipped with a Frobenius-linear endomorphism $\varphi$ such that $\varphi^f = 1\otimes\widetilde{\delta}(\pi_K)$, which is unique up to isomorphism by \cite[Lemma 6.2.3]{KPX_cohomology}.
As for $D_{\cris}^{K_m}(\cR_{K,E}(\widetilde{\delta}^{\wt}))$, since $\widetilde{\delta}^{\wt}(\pi_K)=1$, it extends by continuity to a continuous character of $G_K$ taking the value 1 at any geometric Frobenius.
Since $\widetilde{\delta}$ is smooth, we may also assume that $\widetilde{\delta}^{\wt}|_{1+\pi_K^m\cO_K} = 1$, i.e., the restriction of $\widetilde{\delta}^{\wt}$ to $G_{K_m}$ is the trivial character.
Then, $\widetilde{\delta}^{\wt}$ descends to a character of $G_m$, and 
$$
\DcrisKm(\widetilde{\delta}^\wt) = \left(B_\cris \otimes_{\bbQ_p} E(\widetilde{\delta}^\wt)\right)^{G_{K_m}} = K_0\otimes_{\bbQ_p}E(\widetilde{\delta}^\wt),
$$
which is the free $(K_0\otimes_{\bbQ_p}E)$-module $D_{f,1}$ of rank $1$ on which $G_m$ acts by $\widetilde{\delta}^{\wt}$.
Thus, we conclude that $D_{\cris}^{K_m}(\cR_{K,E}(\delta))$ is the $\varphi$-module $D_{f,\widetilde{\delta}(\pi_K)}$ free of rank $1$ over $(K_0\otimes_{\bbQ_p}E)$, where $G(K_m/K)$ acts by $\widetilde{\delta}^{\wt}$.
To show that the Weil representation $\phi$
attached to $D_{\cris}^{K_m}(\cR_{K,E}(\delta))$ is  $\widetilde{\delta}\circ \rec_K^{-1}$, we follow the recipe recalled in Theorem \ref{thm:D_cris_WD_equivalence}(ii).
For $w\in W_K$, let $\alpha(w) \in f\bbZ$ be the integer such that the image of $w$ in $G(\overline{\bbF_p}/\bbF_p)$ is the $\alpha(w)$-th power of the absolute arithmetic Frobenius.
Then the $W_K$-action is given by
$
\phi(w)={\overline{w}}\circ \varphi^{-\alpha(w)} = \widetilde{\delta}^{\wt}(\rec_K^{-1}(w))\cdot\varphi^{-\alpha(w)},
$
where $\overline{w}$ denotes the image of $w$ in $G_m$.
Since $\rec_K(\pi_K)$ is a geometric Frobenius element, 
$$
\phi(\rec_K(\pi_K)) = \widetilde{\delta}^{\wt}(\pi_K) \varphi^{-(-f)} = \widetilde{\delta}^{\wt}(\pi_K) \widetilde{\delta}(\pi_K) = \widetilde{\delta}(\pi_K).
$$
For $u\in \cO_K^\times = \rec^{-1}_{K}(I_K^{\ab})$, we have $\alpha(\rec_K(u))=0$ and hence
$$
\phi(\rec_K(u)) = \widetilde{\delta}^{\wt}(u)\varphi^{-(0)} = \widetilde{\delta}(u).
$$
We conclude that $\phi\circ \rec_K = \widetilde{\delta}$, and hence $\phi = \widetilde{\delta}\circ \rec_K^{-1}$, which proves the claim.
\end{proof}

\begin{prop}
\label{prop:nc_nonsplit}
For $D\in \Phi \Gamma_{\mathrm{nc}}(\phi, \mathbf{h})$ of regular Sen weight ($h_{1,\sigma}>\dots>h_{n,\sigma})$ for each $\sigma\in \Sigma_K$, $D$ is indecomposable.
In particular, all of its triangulations are strongly non-split.
\end{prop}
\begin{proof}
Any direct sum decomposition $D=D'\oplus D''$ of $D$ with two $(\varphi,\Gamma_K)$-submodules $D',D''$ is part of a triangulation $\Fil_w^\bullet(D)$ for some $w\in S_n$, i.e., $D' \cong \Fil_w^i(D)$ for $i=\rank(D')$ and $D''\cong D/\Fil_w^i(D)$.
This is because $D$ is crystabelline generic, and hence by Theorem \ref{thm:D_cris_WD_equivalence} and Definition \ref{defn:cris-generic}(i), its saturated subobjects correspond to the subobjects of the semisimple Weil representation $\bigoplus_i \phi_i$.
As in Remark \ref{rem:non-split}(i), the problem is reduced to the rank 2 case.
It suffices to show that for any non-critical crystabelline $(\varphi,\Gamma_K)$-module $D$ of rank $2$ of regular weights, if we have an extension $0\to \cR_{K,E}(\delta_1)\to D \to \cR_{K,E}(\delta_2) \to 0$, then this extension is non-split.
This is true since by the non-critical assumption, $\wt_\sigma(\delta_1) > \wt_\sigma(\delta_2)$ for all $\sigma\in \Sigma_K$. 
If the extension is split, we can write $D$ as an  extension $0\to \cR_{K,E}(\delta_2)\to D \to \cR_{K,E}(\delta_1) \to 0$, which is a critical triangulation of $D$.
\end{proof}

\begin{prop}
\label{prop:crit_split}
\begin{enumerate}[(i)]
\item
Let $D$ be a generic crystabelline $(\varphi,\Gamma_{\bbQ_p})$-module over $\cR_{\bbQ_p,E}$ of regular Sen weight ${\bf h}= (h_1 > \dots > h_n)$ and Weil representation $\bigoplus_{i=1}^n \phi_i$.
If $D$ has a critical refinement $w$, then $(D,\Fil_w^\bullet(D))$ is not strongly non-split in the sense of \eqref{eqn:strong-non-split}.
\item
For any nontrivial extension $K/\bbQ_p$, regular Sen weight ${\bf h}= (h_{1,\sigma} > \dots > h_{n,\sigma})_{\sigma\in \Sigma_K}$, and generic Weil representation $\bigoplus_{i=1}^n \phi_i$, there exists a generic crystabelline $(\varphi,\Gamma_{K})$-module $D$ over $\cR_{K,E}$ of regular Sen weight ${\bf h}$ and Weil representation $\bigoplus_{i=1}^n \phi_i$ that has a critical refinement $w$ such that $(D,\Fil_w^\bullet(D))$ is strongly non-split in the sense of \eqref{eqn:strong-non-split}.
\end{enumerate}
\end{prop}
\begin{proof}

\begin{enumerate}[(i)]
\item
Let $D$ be a generic crystabelline $(\varphi,\Gamma_{\bbQ_p})$-module over $\cR_{\bbQ_p,E}$ of regular Sen weight ${\bf h}$ and Weil representation $\bigoplus_{i=1}^n \phi_i$.
Suppose $w\in S_n$ is a critical refinement of $D$, which means that the induced ordering of weights satisfies $k_{i}^w<k_{i+1}^w$ for some $1\leq i\leq n-1$, where $\{k_{i}^w,k_{i+1}^w\}\subset \{h_1,\dots, h_n\}$, cf. Proposition \ref{prop:crys_generic_triang}.
We claim that the extension
$$
0 \to \cR_{\bbQ_p,E}(x^{k_i^w}\phi_{w(i)})  \to \Fil_w^{i+1}(D)/\Fil_w^{i-1}(D) \to \cR_{\bbQ_p,E}(x^{k_{i+1}^w}\phi_{w(i+1)}) \to 0
$$
is split, where the subobject is $\Fil_w^i(D)/\Fil_w^{i-1}(D)$ and the quotient is $\Fil_w^{i+1}(D)/\Fil_w^{i}(D)$.
In the notations of Remark \ref{rem:induced_fil}, since the induced jumps are $-\wt(\delta_{w,i})>-\wt(\delta_{w,i+1})$, the refinement $Ee_{w(i)} \subset Ee_{w(i)}\oplus Ee_{w(i+1)}$ on $D_\dR(\Fil_w^{i+1}(D)/\Fil_w^{i}(D))$ coincides with the Hodge filtration, which forces the other refinement $Ee_{w(i+1)} \subset Ee_{w(i)}\oplus Ee_{w(i+1)}$ to be of relative position $\id$ with the Hodge filtration.
This other refinement gives a triangulation 
$$
0\to \cR_{\bbQ_p,E}(x^{k_{i+1}^w}\phi_{w(i+1)})\to \Fil_w^{i+1}(D)/\Fil_w^{i-1}(D)) \to  \cR_{\bbQ_p,E}(x^{k_{i}^w}\phi_{w(i)}) \to 0.
$$
Since $D$ is generic, $\phi_{w(i)}\neq \phi_{w(i+1)}$.
Thus, by Lemma \ref{lem:rank_one_coh}(iii), the intersection of the subobjects $\cR_{\bbQ_p,E}(x^{k_{i}^w}\phi_{w(i)})$ and $\cR_{\bbQ_p,E}(x^{k_{i+1}^w}\phi_{w(i+1)})$ is zero.
The subobject $\cR_{\bbQ_p,E}(x^{k_{i+1}^w}\phi_{w(i+1)})$ reduces isomorphically to the quotient $\cR_{\bbQ_p,E}(x^{k_{i+1}^w}\phi_{w(i+1)})$ of $\Fil_w^{i+1}(D)/\Fil_w^{i}(D))$ in the first extension above, making it a split extension.
\item
The generic Weil representation $\bigoplus_{j=1}^n \phi_j$ determines an unfiltered $(\varphi,G(L/K))$-module $M$.
To get a crystabelline $(\varphi,\Gamma_K)$-module $D$ of Sen weight ${\bf h}$ with $D_{\cris}^L(D)=M$, we need to give filtrations $\Fil^i(D_{\dR}(D)_\sigma )$ on $D_{\dR}(D)_\sigma := (L\otimes_{L_0} M)_\sigma^{G(L/K)}$ whose jumps occur at $\{-h_{1,\sigma}<\dots<-h_{n,\sigma}\}$ for each $\sigma\in \Sigma_K$.
Since the property of being strongly non-split for $(D,\Fil_w^\bullet(D))$ can be checked at subquotients of rank $2$, it suffices to construct an example $D$ of rank $2$ with Sen weights $(h_{1,\sigma}>h_{2,\sigma})_{\sigma\in \Sigma_K}$ and generic Weil representation $\phi_1\oplus \phi_2$, which has a critical refinement $w\in S_2$ such that $\Fil^\bullet_{w}(D)$ is non-split.
As in Remark \ref{rem:induced_fil}, for the refinement $(\phi_1,\phi_2)$, we choose a basis $({e_{1,\sigma},e_{2,\sigma}})_{\sigma\in \Sigma_K}$ for $D_{\dR}(D) = \prod_\sigma D_{\dR}(D)_\sigma$.
Since $K\neq \bbQ_p$, there are at least two distinct embeddings $\sigma_1,\sigma_2\in \Sigma_K$.
Define
$$
D_{\dR}(D)_{\sigma_1} \supsetneq \Fil^{-h_{1,\sigma_1}+1}(D_{\dR}(D)_{\sigma_1}) := E(e_{1,\sigma_1}+e_{2,\sigma_1}) \supsetneq \Fil^{-h_{2,\sigma_1}+1}(D_{\dR}(D)_{\sigma_1}) = 0,
$$
$$
D_{\dR}(D)_{\sigma_2} \supsetneq \quad \quad \quad \Fil^{-h_{1,\sigma_2}+1}(D_{\dR}(D)_{\sigma_2}) := E(e_{1,\sigma_2}) \supsetneq \Fil^{-h_{2,\sigma_2}+1}(D_{\dR}(D)_{\sigma_2}) = 0,
$$
and define the $\tau$-Hodge filtrations for $\tau\in \Sigma_K\setminus\{\sigma_1,\sigma_2\}$ arbitrarily.
This gives a $(\varphi,\Gamma_K)$-module $D$.
By Remark \ref{rem:induced_fil}, $\Fil_{\id}^\bullet(D)$ induces the ordering $(h_{1,\sigma_1}>h_{2,\sigma_1})$ on the $\sigma_1$-Sen weights and the ordering $(h_{2,\sigma_2}<h_{1,\sigma_2})$ on the $\sigma_2$-Sen weights, but for $w_0 = (1\mbox{ }\mbox{ }2)\in S_2$, $\Fil_{w_0}^\bullet(D)$ induces the ordering $(h_{1,\sigma_i}>h_{2,\sigma_i})$ on the $\sigma_i$-Sen weights for both $i=1,2$.
This implies that $\Fil_{\id}^\bullet(D)$ is non-split, for  in the split case, $\Fil_{w_0}^\bullet(D)$ would induce the orderings $(h_{2,\sigma_1}<h_{1,\sigma_1})$ on the $\sigma_1$-Sen weights and $(h_{1,\sigma_2}>h_{2,\sigma_2})$ on the $\sigma_2$-Sen weights.
\qedhere
\end{enumerate}
\end{proof}

\begin{rem}
For a critical generic crystalline $(\varphi,\Gamma_{K})$-module $D$, it may be strongly non-split for some non-critical refinement $\Fil^\bullet(D)$, and it may have no strongly non-split non-critical refinement at all.
Hence, ``strongly non-split'' depends on the choice of triangulation.
\end{rem}

\subsection{More general cases}

\begin{rem}
We end this section with a heuristic procedure for locating possible triangulations of a fixed $(\varphi,\Gamma_K)$-module $D$ over $\cR_{K,E}$, where $E$ is a finite extension of $\bbQ_p$ containing the Galois closure $K^\norm$ of $K$.
This discussion will not be used later.
As the first step, consider the locus $X_1$ in $\cT_E(E)$ where the character $\delta_1 = \delta$ makes the space 
$$
\set{f_1\in H^0(D^\vee(\delta_1))}{\forall \sigma\in \Sigma_K, \overline{f_1} \neq 0 \in H^0(D^\vee(\delta_1)/t_\sigma)}
$$
nonempty, where $\overline{f_1}$ denotes the image of $f_1$ under the mod-$t_\sigma$ reduction map
$$
\Hom_{\varphi,\Gamma_K}(D,\cR_{K,E}(\delta)) \rightarrow{} \Hom_{\varphi,\Gamma_K}(D/t_\sigma,\cR_{K,E}(\delta)/t_\sigma).
$$
By the classification of $(\varphi,\Gamma_K)$-submodules of $\cR_{K,E}(\delta_1)$ in Lemma \ref{lem:rank_one_coh}(iii), any $f_1$ in the space above corresponds to a map $f_1:D\to \cR_{K,E}(\delta_1)$ of $(\varphi,\Gamma_K)$-modules that is surjective, and vice versa.
To the data of $(\delta_1,f_1)$, we associate a $(\varphi,\Gamma_K)$-submodule of $D$ of rank $n-1$ as
$$
0 \to D_{n-1} := \ker(f_1) \to D \xrightarrow{f_1} \cR_{K,E}(\delta_1) \to 0.
$$
The next step is to consider, for all $(f_1,\delta_1)$ found in Step 1, the locus $X_2$ in $\cT_E(E)$ where the character $\delta_2 = \delta$ makes the space
$$
\set{f_2\in H^0(D_{n-1}^\vee(\delta_2))}{\forall \sigma\in \Sigma_K, \overline{f_2} \neq 0 \in H^0(D_{n-1}^\vee(\delta_2)/t_\sigma)}
$$
nonempty.
Any $f_2$ from this space corresponds to a surjective map from $D_{n-1}$ to $\cR_{K,E}(\delta_2)$ of $(\varphi,\Gamma_K)$-modules, whose kernel (we denote by $D_{n-2}$) satisfies
$$
0 \to D_{n-2} := \ker(f_2) \to D_{n-1} \xrightarrow{f_2} \cR_{K,E}(\delta_2) \to 0.
$$
Continuing this way, we get data $(\delta_1,f_1,\delta_2,f_2,\dots, \delta_{n-1},f_{n-1},\delta_n)$,
where for $D_{n-i}:=\ker(f_i)$ with $D_n = D$, $f_i$ are chosen from
$$
\set{f_i\in H^0(D_{n-i+1}^\vee(\delta_i))}{\forall \sigma\in \Sigma_K, \overline{f_i} \neq 0 \in H^0(D_{n-i+1}^\vee(\delta_i)/t_\sigma)},
$$
from which we know there is a triangulation $\Fil^\bullet(D) = D_{\bullet}$ on $D$ with parameters $(\delta_n,\dots,\delta_1)$.
\end{rem}

\section{Pullback operations}
\label{section:pullback_stack}

Consider the stack $\fX_n = \fX_{\GL_n}$ of $G_K$-equivariant vector bundles of rank $n$ on the Fargues-Fontaine curve $X_{\overline{K}}$, the stack $\fX_P$ of $G_K$-equivariant $P$-bundles on $X_\ovK$ for a standard parabolic subgroup $P\subset \GL_n$ containing the upper triangular Borel subgroup $B\subset \GL_n$ with Levi quotient $M \cong \GL_{n_1} \times \dots \times \GL_{n_r}$, and the stack $\fX_M$ of $G_K$-equivariant $M$-bundles on $X_\ovK$, all over the category of rigid $E$-analytic spaces $\Rig_E$ equipped with the Tate-fpqc topology, cf. \cite[\S5.1, \S5.3]{EGH_CatLLC}.
Note that $\fX_M \cong \fX_{n_1}\times \dots \times \fX_{n_r}$.
When $P=B$, the Levi quotient $M$ is the diagonal torus, which we denote by $T \cong (\GL_1)^n$ so that $\fX_T = \fX_1 \times \dots \times \fX_1$.

By \cite[5.1.1 and 5.1.5]{EGH_CatLLC}, a point of $\fX_n(\Sp(A))$ is a $(\varphi,\Gamma_K)$-module $D$   of rank $n$ over $\cR_{K,A}$, and a point of $\fX_P(\Sp(A))$ is a pair $(D,\Fil^\bullet(D))$ consisting of a $(\varphi,\Gamma_K)$-module $D$ of rank $n$ over $\cR_{K,A}$ together with a filtration 
$$
0 = \Fil^0(D) \subset \Fil^1(D) \subset \dots \subset \Fil^r(D) = D
$$
by $(\varphi,\Gamma_K)$-submodules with the subquotients $\Fil^i(D)/\Fil^{i-1}(D)$ being locally free over $\cR_{K,A}$ of rank $n_i$ for all $1\leq i\leq r$, in which case we say that $\Fil^\bullet(D)$ is a {\bf $P$-parabolization} on $D$.
In general, a $(\varphi,\Gamma_K)$-module $D$ of rank $n$ is called {\bf paraboline} if it admits a $P$-parabolization for some proper parabolic subgroup $P$ of $\GL_n$ containing $B$.

For $P_i$ the parabolic subgroup containing $B$ of Levi quotient $M_i \cong \GL_{n-i}\times \GL_{i}$, we may consider the pullback map $p_{i,\sigma}: \fX_{P_i} \to \fX_{P_i}$, sending a $P_i$-bundle
$$
(0\to D_{n-i} \to D \to D^i:= D/D_{n-i} \to 0),
$$
where $D_{n-i}$ is a saturated $(\varphi,\Gamma_K)$-submodule of $D$ of rank $n-i$, to the $P_i$-subbundle
$$
(0 \to D_{n-i} \to p_{i,\sigma}(D) \to t_\sigma D^i \to 0).
$$
This makes sense on $\fX_Q$ for any parabolic $Q$ such that $B \subset Q \subset P_i$, where $p_{i,\sigma}$ acts similarly.

\begin{rem}
On $\fX_B$, consider $p_{i,\sigma}: \fX_B \to \fX_B$. Over the sublocus of $\cT^n$ where $\wt_{\sigma}(\delta_{j}/\delta_{k}) \neq 1$ for all $1\leq j\leq n-i<k \leq n$ with fixed $i$, and for a non-split pair $(D,\Fil^\bullet(D))\in \fX_B(A)$, it follows from Theorem \ref{thm:p_k-invertible} that from $p_{i,\sigma}(D,\Fil^\bullet(D))$ together with its induced triangulation, we can recover $D$ together with its triangulation.
\end{rem}

\begin{lem}
\label{lem:p_i_commutes_with_p_j}
For any $1\leq i < j\leq n$ and $\sigma,\tau\in \Sigma_K$ (with possibly $\sigma = \tau$), we have that 
$$ 
p_{i,\sigma} \circ p_{j,\tau} = p_{j,\tau} \circ p_{i,\sigma}
$$ 
as maps from $\fX_Q$ to $\fX_Q$, for any parabolic $Q$ such that $B\subset Q\subset P_i\cap P_j$.
\end{lem}
\begin{proof}
We may assume $Q=P_i\cap P_j$.
Then, a $Q$-bundle is $D\in \fX_n(A)$ together with a filtration
$$
0 \subset D_{n-j}:=\Fil^{n-j}(D) \subset D_{n-i} := \Fil^{n-i}(D) \subset D_n = D
$$
by $(\varphi,\Gamma_K)$-submodules over $\cR_{K,A}$ such that $D^\bullet:=D/D_{n-\bullet}$ is locally free of rank $\bullet\in \{i,j\}$, and $p_{i,\sigma}(D,\Fil^\bullet(D)) = (D_{n-i} - t_\sigma D^i)$, which has the induced $Q$-filtration 
$$
0 \subset D_{n-j} \subset D_{n-i} \subset p_{i,\sigma}(D,\Fil^\bullet(D)).
$$
From this we conclude that, because the twist by $t_{\tau}$ commutes with pullback via $t_\sigma$,
\begin{align*}
p_{j,\tau}(p_{i,\sigma}(D,\Fil^\bullet(D))) 
&= (D_{n-j} - t_{\tau}(D_{n-i}/D_{n-j} - t_\sigma D^i))
\\&=
(D_{n-j}-t_{\tau}(D_{n-i}/D_{n-j})-t_\tau t_{\sigma} D^i)
\\&=
(D_{n-j}-t_{\tau}(D_{n-i}/D_{n-j})-t_\sigma (t_{\tau} D^i))
=
p_{i,\sigma}(p_{j,\tau}(D,\Fil^\bullet(D))).
\qedhere
\end{align*}
\end{proof}

\subsection{Substacks}

We define several substacks of $\fX_n$ in this subsection. 
We will see in Theorem \ref{thm:comparing_Wu}(i) that $p_{i,\sigma}$ descend to them.

\begin{defn}
\begin{enumerate}[(i)]
\item 
Recall we have the morphisms of stacks
$$
\begin{tikzcd}
& \arrow[dl, "\beta_P"'] \fX_P \arrow[dr, "\alpha_P"] & \\
\fX_n &  & \fX_M ,
\end{tikzcd}
$$
where for a paraboline $(\varphi,\Gamma_K)$-module $(D,\Fil^\bullet(D)) \in \fX_P(A)$, we have the forgetful map 
$$
\beta_P(D,\Fil^\bullet(D)) = D,
$$
and the map 
$$
\alpha_P(D,\Fil^\bullet(D)) = (\Fil^i(D)/\Fil^{i-1}(D))_i
$$
taking the successive quotients of the $P$-structure $\Fil^\bullet(D)$ on $D$.
\item 
For $P=B$, we define the {\bf weight maps} 
$$
\omega_T: \fX_T \to \Res_{K/\bbQ_p}(\bbA^{n,\rig})_E \cong (\bbA_E^{n})^{\Sigma_K,\rig}, \quad (\delta_1,\dots,\delta_n) \mapsto (\wt(\delta_1),\dots,\wt(\delta_n)),
$$
and
$$
\omega_B: \fX_B \xrightarrow{\alpha_B} \fX_T \xrightarrow{\omega_T} (\bbA_E^{n})^{\Sigma_K,\rig}, \quad (D,\Fil^\bullet(D)) \mapsto (\wt(\delta_1),\dots,\wt(\delta_n)),
$$
where $\delta_i$ are the unique characters determined by the rank 1 subquotients
$$
\Fil^i(D)/\Fil^{i-1}(D) \cong \cR_{K,A}(\delta_i){\otimes}_{\cO_{\Sp(A)}} \cL_i.
$$
With the projection
$$\pr_\sigma:(\bbA_E^{n})^{\Sigma_K,\rig} \to \bbA^{n,\rig}_E$$
onto the $\sigma$-component for $\sigma\in \Sigma_K$, we define the {\bf $\sigma$-weight maps} by 
$$
\omega_{T,\sigma}:= \pr_\sigma\circ \omega_T : \fX_T \to \bbA^{n,\rig}_E , \quad (\delta_1,\dots,\delta_n) \mapsto (\wt_\sigma(\delta_1),\dots,\wt_\sigma(\delta_n)),
$$
and
$$
\omega_{B,\sigma}:= \omega_{T,\sigma}\circ \alpha_B: \fX_B \to \bbA_E^{n,\rig}, \quad (D,\Fil^\bullet(D)) \mapsto (\wt_\sigma(\delta_1),\dots,\wt_\sigma(\delta_n)).
$$
\end{enumerate}
\end{defn}

\begin{defn}
\label{defn:uniform_tri}
Let $S\subset \Sigma_K$ be a subset of embeddings of $K$ into $E$.
\begin{enumerate}[(i)]
\item 
For a trianguline $(\varphi,\Gamma_K)$-module $D$ over $\cR_{K,A}$, we say $D$ is {\bf $S$-weight-uniform trianguline} if for every rigid point $z\in \Sp(A)$ and every finite extension $L/k(z)$, the $\sigma$-weight map $\omega_{B,\sigma}$ is constant on $\beta_B^{-1}(D_z\otimes_{k(z)}L)$ for every $\sigma\in S$. That is, all triangulations on the fiber $D_z\otimes_{k(z)}L$ induce the same ordering on its $\sigma$-Sen weights for all $\sigma\in S$.
\item 
As the $S$-weight-uniform trianguline condition is defined via geometric fibers, it is stable under base change. Denote by $\fX_B^{S\text{-}\unif}$ this substack of $\fX_B$, and by $\fX_n^{S\text{-}\unif}$ the sheafification of its image under $\beta_B$ in $\fX_n$ with respect to the Tate-fpqc topology on $\Rig_E$.
We refer to $\fX_B^{S\text{-}\unif}$ and $\fX_n^{S\text{-}\unif}$ as the {\bf $S$-weight-uniform substacks}.
\item 
If $S=\Sigma_K$, we simply write $\fX_B^{\unif}$ and $\fX_n^{\unif}$, and call them {\bf weight-uniform substacks}.
\end{enumerate}
\end{defn}

\begin{rem}
Sheafification is necessary for descent-theoretic reasons. \cite[Example 5.3.2]{EGH_CatLLC} gives a filtered $\varphi$-module on $\bbP^1$ defining a $(\varphi,\Gamma_{\bbQ_p})$-module of rank 2 on $\bbP^1$, which has no global triangulation on $\bbP^1$, but is globally trianguline over the two copies of $\bbA^1$.
\end{rem}

\begin{defn}
\label{defn:Xn_SwuIk}
\begin{enumerate}[(i)]
\item 
Fix $\sigma\in \Sigma_K$ and $i\in \{1,\dots,n\}$.
\begin{enumerate}
\item
Let $U_{i}$ be the Zariski-open subspace of $\bbA_E^{n,\rig}$ that is complement to the hyperplanes $\{T_{j}-T_{k}=0\}$ for all $1\leq j \leq n-i$ and $n-i+1\leq k\leq n$, where we denote by $\set{T_{l}}{1\leq l\leq n}$ the standard coordinates on $\bbA_E^{n,\rig}$.
\item
Let $\fX_n^{\sigma\text{-}\unif,i}$ denote the sheafification of the image $\beta_B(\fX_B^{\sigma\text{-}\unif} \cap \omega_{B,\sigma}^{-1}(U_{i}))$, which is the substack of $\fX_n^{\sigma\text{-}\unif}$ characterized by the property that $D\in \fX_n^{\sigma\text{-}\unif,i}(A)$ if and only if for any $z\in \Sp(A)$, any finite field extension $L/k(z)$, and any\footnote{By $\sigma$-weight-uniformity, all triangulations of $D_z \otimes_{k(z)} L$ induce the same ordering on $\sigma$-Sen weights. Since membership in $\omega_{B,\sigma}^{-1}(U_i)$ is precisely the condition $h_{j, \sigma} \neq h_{k, \sigma}$ for $j \leq n-i<k$, the weight-separation condition holds for one triangulation if and only if it holds for every triangulation.} triangulation on the fiber $D_z\otimes_{k(z)} L$, its first $n-i$ ordered $\sigma$-Sen weights are disjoint from its last $i$ ordered $\sigma$-Sen weights in $\overline{k(z)}$.
\item 
For $a\leq 0 \leq  b \in \bbZ$, let $\fX_n^{\sigma\text{-}\unif,i,[a,b]}$ be the substack of $\fX_n^{\sigma\text{-}\unif}$ given by the condition that $D\in \fX_n^{\sigma\text{-}\unif,i,[a,b]}(A)$ if and only if for any $z\in \Sp(A)$, any finite field extension $L/k(z)$, and any triangulation on $D_z\otimes_{k(z)}L$ with parameters $(\delta_{1,z},\dots,\delta_{n,z})$, we have
$$
\set{\wt_\sigma(\delta_{j,z})}{1\leq j\leq n-i} 
\cap 
\set{\wt_\sigma(\delta_{k,z})+h }{n-i+1\leq k \leq n, h\in [a,b]\cap \bbZ}
= 
\emptyset
$$
in $\overline{k(z)}$.
In particular, $\fX_n^{\sigma\text{-}\unif,i}$ defined in (b) is equal to $\fX_n^{\sigma\text{-}\unif,i,[0,0]}$.
\end{enumerate}
\item 
Fix $S\subset \Sigma_K$, and for each $\sigma\in S$, choose a subset 
$$
I_\sigma := \{1\leq i_1<\dots<i_{d_\sigma}\leq n\}
$$
of $\{1,\dots,n\}$ of size $d_\sigma$.
Set $i_0 = 0$ and $i_{d_{\sigma}+1} = n$.
For ${\bf k} = (k_{\sigma,i})_{\sigma\in S,i\in I_\sigma} \in  \bbN^{\sum_{\sigma\in S}d_\sigma}$, let
$$
\bigcap_{\sigma\in S} \bigcap_{m=1}^{d_\sigma} \fX_{n}^{\sigma\text{-}\unif, i_m, \left[0,\sum_{r=1}^m k_{\sigma,i_r}\right]} \subset \fX_{n}^{S\text{-}\unif, I, {\bf k}} \subset \bigcap_{\sigma\in S} \bigcap_{m=1}^{d_\sigma} \fX_{n}^{\sigma\text{-}\unif, i_m, \left[0,k_{\sigma,i_m}\right]}
$$
be the substack of $\fX_{n}^{S\text{-}\unif}$ given by the condition that $D\in \fX_n^{S\text{-}\unif,I, {\bf k}}(A)$ if and only if for any $z\in \Sp(A)$, any finite field extension $L/k(z)$, and any triangulation on $D_z\otimes_{k(z)}L$ with parameters $(\delta_{1,z},\dots,\delta_{n,z})$, setting $(h_{i,\sigma})_z:=\wt_\sigma(\delta_{i,z})$, for all $\sigma\in S$ we have that
$$
\forall 0\leq m \leq d_\sigma, \forall n+1-i_{m+1} \leq j< n+1-i_m, \forall 0\leq m' < m, \forall n+1-i_{m'+1} \leq k < n+1-i_{m'}
$$
$$
(h_{j,\sigma})_z \notin \set{(h_{k,\sigma})_z+a}{a\in\bbN, 0\leq a \leq \sum_{r=m'+1}^{m} k_{\sigma, i_r}}
$$
in $\overline{k(z)}$.
On the other hand, for $-{\bf k} = (-k_{\sigma,i})\in \bbZ_{\leq 0}^{\sum_{\sigma\in S} d_\sigma}$, let 
$$
\bigcap_{\sigma\in S} \bigcap_{m=1}^{d_\sigma} \fX_{n}^{\sigma\text{-}\unif, i_m, \left[-\sum_{r=1}^m k_{\sigma,i_r},0\right]} \subset \fX_{n}^{S\text{-}\unif, I, {\bf -k}} \subset \bigcap_{\sigma\in S} \bigcap_{m=1}^{d_\sigma} \fX_{n}^{\sigma\text{-}\unif, i_m, \left[-k_{\sigma,i_m},0\right]}
$$
be the substack of $\fX_{n}^{S\text{-}\unif}$ given by the condition that $D\in \fX_n^{S\text{-}\unif,I, {\bf -k}}(A)$ if and only if for any $z\in \Sp(A)$, any finite field extension $L/k(z)$, and any triangulation on $D_z\otimes_{k(z)}L$ with parameters $(\delta_{1,z},\dots,\delta_{n,z})$, setting $(h_{i,\sigma})_z:=\wt_\sigma(\delta_{i,z})$, for all $\sigma\in S$ we have that
$$
\forall 0\leq m \leq d_\sigma, \forall n+1-i_{m+1} \leq j< n+1-i_m, \forall 0\leq m' < m, \forall n+1-i_{m'+1} \leq k < n+1-i_{m'}
$$
$$
(h_{j,\sigma})_z \notin \set{(h_{k,\sigma})_z -a}{a\in\bbN, 0\leq a \leq \sum_{r=m'+1}^{m} k_{\sigma, i_r}}.
$$
\item 
When $S=\Sigma_K$ and $I_\sigma = \{1,\dots,n\}$ for all $\sigma$, we simply write $\fX_{n}^{\unif, {\bf k}}$ to ease the notation.
\end{enumerate}
\end{defn}

\begin{rem}
\label{rem:SwuIk}
One should think of the conditions defining $\fX_{n}^{S\text{-}\unif, I, {\bf -k}}$ as the requirements on weights such that changing the weights through the pullback $\prod_{\sigma\in S, i \in I_\sigma} p_{i,\sigma}^{k_{\sigma,i}}$ does not meet any of the relevant walls in the weight space.
\end{rem}

\begin{cor}
The following two classes of $(\varphi,\Gamma_K)$-modules of rank $n$ over $\cR_{K,A}$ give rise to $\Sp(A)$-valued points of $\fX_n^\unif$:
\begin{enumerate}[(i)]
\item 
if $D\in \fX_n(A)$ is trianguline and $D_z$ is non-critical crystabelline at all $z\in \Sp(A)$.
\item 
if $D\in \fX_n(A)$ has a strongly non-split triangulation $\Fil^\bullet(D)$ with parameters in $\cT^n_\circ(A)$.
\end{enumerate}
\end{cor}
\begin{proof}
(i) follows from Definition \ref{defn:ncr}, and (ii) follows from Corollary \ref{cor:generic-unique}.
\end{proof}

\subsection{Operation $p_{i,\sigma}$ on $\fX_n^{\sigma\mh\unif,i}$}

We begin by recalling Wu's result obtained in \cite{Wu_translation} that we rely on, and then explain how it extends $p_{i,\sigma}$ and allows us to define morphisms on our stacks.

\begin{lem}
\label{lem:coprime-polynomials}
\begin{enumerate}[(i)]
\item 
Let $R$ be any commutative unital ring, and let $Q,S,Q',S'\in R[T]$ be  monic polynomials such that $\deg(Q)= \deg(Q'), \deg(S) = \deg(S')$, $Q(T)S(T)=Q'(T)S'(T)$ and $(Q',S) = (1)$.
Then, $Q=Q'$ and $S=S'$.
\item 
Let $k$ be a field of characteristic $0$, and let $Q,S,Q',S'\in k[T]$ be monic polynomials such that $QS=Q'S'$ and $Q(T-1)S(T)=Q'(T-1)S'(T)$.
Then, $Q=Q'$ and $S=S'$.
\item 
Let $A$ be an affinoid $\bbQ_p$-algebra, and let $Q,S\in A[T]$ be polynomials.
Then, $(Q(T), S(T))=(1)$ if and only if for any $z \in \Sp(A)$, the sets of roots of $Q(T) \otimes_A k(z)$ and of $S(T) \otimes_A k(z)$ in $\overline{k(z)}$ have empty intersection.
\end{enumerate}
\end{lem}

\begin{proof}
\begin{enumerate}[(i)]
\item 
Subtracting $Q'(T)S(T)$ from both sides of $Q(T)S(T)=Q'(T)S'(T)$, we get
$$
(Q-Q')S=Q'(S'-S).
$$
Since $(Q',S)=1$, we find polynomials $A,B\in R[T]$ such that $AQ'+BS=1$.
Multiplying both sides of $AQ'+BS=1$ by $(S'-S)$, we get
$$
AQ'(S'-S)+BS(S'-S) = S'-S.
$$
Each of the two terms on the left is divisible by $S$, so is the right-hand side $S'-S$.
Thus there exists $C\in R[T]$ such that $S'-S =CS$, or equivalently,
$$
S' = (1+C)S.
$$
But by assumption, $S$ and $S'$ are monic polynomials of the same degree.
Hence, $C=0$ and $S=S'$.
Then $(Q-Q')S=0$ in $R[T]$.
Since $S$ is monic, we have $Q=Q'$.
\item 
The relation $QS=Q'S'$ implies $\frac{Q}{Q'}=\frac{S'}{S}\in k(T)$.
The relation $Q(T-1)S(T) = Q'(T-1)S'(T)$ implies $\frac{Q(T-1)}{Q'(T-1)} = \frac{S'(T)}{S(T)}$.
Thus, the rational function $C(T):=Q(T)/Q'(T)\in k(T)$ satisfies $C(T)=C(T-1)$, and hence
$
C(T)=C(T-1)=C(T-2) = C(T-3) = \dots.
$
Since $k$ has characteristic $0$, the rational function $C(T)=Q(T)/Q'(T)$ equals a constant $c\in k$. Then, $c=1$ since both $Q$ and $Q'$ are monic.
Thus, $1=\frac{Q}{Q'}=\frac{S'}{S}$ as desired.
\item 
This is \cite[Lemma 3.15]{Wu_translation}.
\qedhere
\end{enumerate}
\end{proof}

Let $D$ be a $(\varphi,\Gamma_K)$-module of rank $n$ over $\cR_{K,A}$. 
Fix an embedding $\sigma \in \Sigma_K$ and denote by $P_{\Sen}(T) \in \left(K \otimes_{\mathbb{Q}_p} A\right)[T]$ the Sen polynomial of $D$, and by $P_{\Sen, \sigma}(T)$ the $\sigma$-Sen polynomial, i.e., the $\sigma$-component of $P_{\text {Sen}}(T)$ via $(K \otimes_{\mathbb{Q}_p} A)[T] \simeq \prod_{\sigma\in \Sigma_K} A[T]$.

\begin{thm}[{\cite[Proposition 3.16]{Wu_translation}}]
\label{thm:Wu_3.16}
Suppose that the $\sigma$-Sen polynomial $P_{\text {Sen}, \sigma}(T)$ of $D$ admits a decomposition $P_{\text {Sen}, \sigma}(T)=Q(T) S(T)$ in $A[T]$ by comaximal monic polynomials. 
\begin{enumerate}[(i)]
\item
There exists a unique ($\varphi, \Gamma_K$)-module $D^{\prime}$ over $\mathcal{R}_{K,A}$ contained in $D$ and containing $t_\sigma D$ such that the Sen polynomial of $D^{\prime}$ is equal to 
$$
Q(T-1) S(T) \prod_{\tau \neq \sigma} P_{\mathrm{Sen}, \tau}(T) \in \prod_{\sigma \in \Sigma_K} A[T].
$$
\item 
There exists a unique ($\varphi, \Gamma_K$)-module $D^{''}$ over $\mathcal{R}_{K,A}$ contained in $t_\sigma^{-1}D$ and containing $D$ such that the Sen polynomial of $D^{''}$ is equal to 
$$
Q(T+1) S(T) \prod_{\tau \neq \sigma} P_{\mathrm{Sen}, \tau}(T) \in \prod_{\sigma \in \Sigma_K} A[T].
$$
\item 
If the image of each difference between the roots of $Q(T)$ and $S(T)$ in $\overline{k(z)}$ never belongs to $\{-1,0,1\}$ for all $z\in \Sp(A)$, then the operations in $(i)$ and $(ii)$ are mutual inverses.
\end{enumerate}
\end{thm}
\begin{proof}
\begin{enumerate}[(i)]
\item 
This part is {\cite[Proposition 3.16]{Wu_translation}} except that it is stated there that the submodule $D'$ contains $tD$, instead of $t_\sigma D$.
Let us sketch Wu's construction and the proof of its uniqueness, and then indicate why his proof shows that $t_\sigma D \subset D' \subset D$, which is {\it a priori} stronger than $tD \subset D' \subset D$.
\vskip8pt

By Beauville-Laszlo gluing \cite[Proposition A.3]{Wu_translation}, it suffices to look at the $\sigma$-component $D_{\dif,\sigma}^{+}(D)$ of the ``localization'' $D_{\dif}^{+}(D):= D\otimes_{\cR_{K,A}} (K_\infty\otimes_{\bbQ_p} A)[[t]]$, which is finite projective of rank $n$ over $(K_\infty\otimes_{K,\sigma}A)[[t]]$ with a semilinear $\Gamma_K$-action, and show that $D_{\dif,\sigma}^{+}(D)$ has a unique $\Gamma_K$-stable projective submodule $M$ containing $t D_{\dif,\sigma}^{+}(D)$ such that the Sen operator $\Theta_\sigma$ on $D_{\Sen,\sigma}(M) = M/tM$ has $\sigma$-Sen polynomial $Q(T-1)S(T)$. 
Then, the $(\varphi,\Gamma_K)$-submodule $D'$  corresponds to the $\Gamma_K$-stable lattice in $D_{\dif}^+(D)[1/t]$ given by
$$
D_{\dif,\tau}^+(D') = 
\begin{cases}
D_{\dif,\tau}^+(D) &\text{ if } \tau \in \Sigma_K\setminus \{\sigma\},\\
M &\text{ if } \tau = \sigma.
\end{cases}
$$
This submodule $M$ is explicitly given by
$$
M:= \ker(D_{\dif,\sigma}^{+}(D) \to D_{\dif,\sigma}^{+}(D)/t \to \ker(Q(\Theta_\sigma)|D_{\Sen,\sigma}(D))),
$$
where we have the canonical decomposition as $A[T]$-module on which $T$ acts by $\Theta_\sigma$,
$$
D_{\Sen,\sigma}(D) = \ker(Q(\Theta_\sigma)|D_{\Sen,\sigma}(D)) \oplus \ker(S(\Theta_\sigma)|D_{\Sen,\sigma}(D)),
$$
thanks to the assumption that $P_{\Sen,\sigma}(T)=Q(T)S(T)$ and $(Q(T),S(T))=1$.

The uniqueness also follows from this canonical decomposition of $D_{\Sen,\sigma}(D)$.
Indeed, let $R:=K_\infty\otimes_{K,\sigma}A$ and let $M$ be any $\Gamma_K$-stable projective $R[[t]]$-submodule of $D_{\dif,\sigma}^{+}(D)$ containing $t D_{\dif,\sigma}^{+}(D)$ so that the Sen operator $\Theta_\sigma$ on $D_{\Sen,\sigma}(M) = M/tM$ has characteristic polynomial $Q(T-1)S(T)$.
Then from 
$$
t^2 D_{\dif,\sigma}^{+}(D) \subset t M \subset t D_{\dif,\sigma}^{+}(D) \subset M \subset D_{\dif,\sigma}^{+}(D),
$$
we obtain two exact sequences of $\Gamma_K$-stable $R$-modules
\begin{equation}
\label{eqn:M/tM}
0 \to tD_{\dif,\sigma}^{+}(D)/tM \to M/tM \to M/tD_{\dif,\sigma}^{+}(D) \to 0,
\end{equation}
\begin{equation}
\label{eqn:Ddif+/tDdif+}
0 \to M/tD_{\dif,\sigma}^{+}(D) \to D_{\dif,\sigma}^{+}(D)/tD_{\dif,\sigma}^{+}(D) \to D_{\dif,\sigma}^{+}(D)/M \to 0,
\end{equation}
which, as we now show, following \cite[Lemma 1.1.5]{Zhu_affineGrassmann}, consist of projective $R$-modules.
Since $M$ is projective over $R[[t]]$, $M/tM$ is projective over $R = R[[t]]/tR[[t]]$.
Similarly, $D_{\dif,\sigma}^{+}(D)/tD_{\dif,\sigma}^{+}(D)$ is projective over $R$.
Since $M$ is an $R[[t]]$-lattice in the projective $R((t))$-module $D_{\dif,\sigma}(D) := D_{\dif,\sigma}^+(D)[1/t]$ of rank $n$, it follows that
$$
D_{\dif,\sigma}(D)/ M = \bigoplus_{k\geq 0} t^{-(k+1)}M/t^{-k}M
$$
is projective over $R$ since $M/tM$ and so $t^{-(k+1)}M/t^{-k}M$ are projective.
Similarly, 
$$
D_{\dif,\sigma}(D)/ D_{\dif,\sigma}^+(D) = \bigoplus_{k\geq 0} t^{-(k+1)}D_{\dif,\sigma}^+(D)/t^{-k}D_{\dif,\sigma}^+(D)
$$
is projective over $R$ because $D_{\dif,\sigma}^{+}(D)/tD_{\dif,\sigma}^{+}(D)$ is projective over $R$.
Hence, from
$$
0 \to D_{\dif,\sigma}^+(D)/M \to D_{\dif,\sigma}(D)/D_{\dif,\sigma}^+(D) \to D_{\dif,\sigma}(D)/M \to 0,
$$
we conclude that $D_{\dif,\sigma}^+(D)/M$ is projective.
By \eqref{eqn:Ddif+/tDdif+}, $M/tD_{\dif,\sigma}^{+}(D)$ is projective.
For the sequences \eqref{eqn:M/tM} and \eqref{eqn:Ddif+/tDdif+} of finite projective $\Gamma_K$-stable $R$-modules, we look at the derivative of $\Gamma_K$-action, i.e., the Sen operator $\Theta_\sigma$, to get the factorizations
\begin{align*}
P_{\Sen,M/tM }(T) &= P_{\Sen, tD_{\dif,\sigma}^{+}(D)/tM}(T) \cdot P_{\Sen, M/tD_{\dif,\sigma}^{+}(D)}(T),\\
P_{\Sen,D_{\dif,\sigma}^{+}(D)/tD_{\dif,\sigma}^{+}(D)}(T) &= P_{\Sen, M/tD_{\dif,\sigma}^{+}(D)}(T) \cdot P_{\Sen, D_{\dif,\sigma}^{+}(D)/M}(T).
\end{align*}
Since the derivative of $\Gamma_K$ on $D_{\dif,\sigma}^+(D)$ is a derivation $\nabla: D_{\dif,\sigma}^+(D)\to D_{\dif,\sigma}^+(D)$ over the derivation $t\frac{\partial}{\partial t}$ on the coefficient $R[[t]]$, cf. \cite[\S3.4]{Fontaine_arithmetique}, it follows that 
$$
\nabla(ty) = ty + t\nabla y
$$
for any $y \in D_{\dif,\sigma}^+(D)$, from which we see
$$
P_{\Sen, tD_{\dif,\sigma}^{+}(D)/tM}(T) = P_{\Sen, D_{\dif,\sigma}^{+}(D)/M}(T-1).
$$
Thus, $Q'(T):=P_{\Sen, D_{\dif,\sigma}^{+}(D)/M}(T)$ and $S'(T):=P_{\Sen, M/tD_{\dif,\sigma}^{+}(D)}(T)$ are monic polynomials in $A[T]$ satisfying the equalities
\begin{align*}
Q'(T)S'(T) &= P_{\Sen,\sigma}(T) = Q(T)S(T),
\\
Q'(T-1)S'(T) &= P_{\Sen, M/tM}(T) = Q(T-1)S(T).    
\end{align*}
For each $z\in \Sp(A)$, we can specialize to the polynomial ring $k(z)[T]$ over the residue field $k(z)$ and apply Lemma \ref{lem:coprime-polynomials}(ii) to deduce that $Q_z(T)=Q'_z(T)$ and $S_z(T)=S'_z(T)$.
Hence, $\deg(Q')=\deg(Q)$ and $\deg(S)=\deg(S')$, and by Lemma \ref{lem:coprime-polynomials}(iii), we deduce comaximality $(Q,S)=(Q',S)=(Q,S')=(Q',S')=A[T]$, and that $Q'(T)=Q(T)$ and $S'(T)=S(T)$ in $A[T]$ by Lemma \ref{lem:coprime-polynomials}(i).
By Chinese remainder theorem, we have a natural isomorphism
$$
A[T] \xrightarrow[\sim]{f} A[T]/Q(T) \times A[T]/S(T),
$$
and we let $e_Q =f^{-1}(1,0)$ and $e_S = f^{-1}(0,1)$ be the idempotents in $A[T]$.
Then, for any $R$-linear section $s': D_{\dif,\sigma}^+(D)/M \to D_{\dif,\sigma}^+(D)/t$ that splits \eqref{eqn:Ddif+/tDdif+} as $R$-modules, we get an $R[\Theta_\sigma]$-linear section $s:=e_{Q}(\Theta_\sigma)\circ s'$ that splits \eqref{eqn:Ddif+/tDdif+} as Sen-modules by the choice of $e_Q$ and the fact that $Q'=Q$ and $S=S'$.
This must be the canonical decomposition 
$$
D_{\Sen,\sigma}(D) = \ker(Q(\Theta_\sigma)|D_{\Sen,\sigma}(D)) \oplus \ker(S(\Theta_\sigma)|D_{\Sen,\sigma}(D)).
$$
Therefore, $M/tD_{\dif,\sigma}^+(D) = \ker(S(\Theta_\sigma)|D_{\Sen,\sigma}(D))$ and
$$
M = \ker(D_{\dif,\sigma}^{+}(D) \to D_{\dif,\sigma}^{+}(D)/t \to \ker(Q(\Theta_\sigma)|D_{\Sen,\sigma}(D))),
$$
which proves the uniqueness of such $\Gamma_K$-stable projective submodule.

\vskip8pt

By construction, $t D_\dif^+(D) \subset D_\dif^+(D') \subset D_\dif^+(D)$.
Hence, we know that $tD \subset D' \subset D$ by Beauville-Laszlo.
To show that $t_\sigma D \subset D' \subset D$, it suffices to note\footnote{One way to see them is to apply the uniqueness proven in the preceding paragraph to $\cR_{K,E}$, which yields  $t_\sigma\cR_{K,E}$ as it has the correct Sen polynomial and satisfies $t\cR_{K,E} \subset t_\sigma\cR_{K,E} \subset \cR_{K,E}$.
Comparing $D_{\dif,\bullet}^+$ yields
$$
t_\sigma (K_\infty \otimes_{K,\tau} E [[t]]) = 
\begin{cases}
K_\infty \otimes_{K,\tau} E [[t]] &\text{ if } \tau \neq \sigma, \\
t (K_\infty \otimes_{K,\tau} E [[t]]) &\text{ if } \tau = \sigma,
\end{cases}
$$
which implies the desired result.} 
that 
\begin{itemize}
\item
$
t_\sigma 
$
acts as a unit on the components $D_{\dif,\tau}^+(D)$ for $\tau \neq \sigma$, and
\item 
$t_\sigma$ acts as $t$ on the $\sigma$-component $D_{\dif,\sigma}^+(D)$.
\end{itemize}
This means that in the construction of $M \subset D_{\dif,\sigma}^+(D)$ above and the proof of uniqueness of such $M$, we may replace $t$ by $t_\sigma$ and obtain the same result.

\item 
This is a corollary to (i), after replacing $D$ by $t_\sigma^{-1}D$ and switching the roles played by $Q(T)$ and $S(T)$ in the statement of (i):
the Sen polynomial of $t_\sigma^{-1}D$ is exactly
$$
P_{\Sen, \sigma}(T+1) \prod_{\tau \neq \sigma} P_{\mathrm{Sen}, \tau}(T) = Q(T+1) S(T+1) \prod_{\tau \neq \sigma} P_{\mathrm{Sen}, \tau}(T) 
\in \prod_{\sigma \in \Sigma_K} A[T].
$$
The $\sigma$-component of $D_{\dif}^+(t_\sigma^{-1}D) = t_\sigma^{-1}D_\dif^+(D)$ contains the $\Gamma_K$-subrepresentation
\begin{align*}
M &= \ker(D_{\dif,\sigma}^+(t_\sigma^{-1}D) \to D_{\dif,\sigma}^+(t_\sigma^{-1}D)/t \to \ker(S(\Theta_\sigma+1)|D_{\Sen,\sigma}(t_{\sigma}^{-1}D))),
\end{align*}
which is the unique $\Gamma_K$-stable projective submodule containing 
$$
t D_{\dif,\sigma}^{+}(t_\sigma^{-1}D) = t_\sigma D_{\dif,\sigma}^{+}(t_\sigma^{-1}D) = D_\dif^+(D)
$$
such that $\Theta_\sigma$ acts on $D_{\Sen,\sigma}(M)$ with Sen polynomial $$
Q(T+1)S(T+1-1) = Q(T+1)S(T)
$$
by (i).
The result then follows.
\item 
The assumption ensures $(Q(T),S(T)) = 1 = (Q(T-1),S(T)) = (Q(T+1),S(T))$. 
The desired conclusion follows from the uniqueness part of (i) and (ii).
\qedhere
\end{enumerate}
\end{proof}

\begin{prop}
\label{prop:Wu_trianguline}
Let $(D,\Fil^\bullet(D)) \in \fX_B(A)$ be a triangulated $(\varphi,\Gamma_K)$-module over $\cR_{K,A}$ with parameters $(\delta_1,\dots,\delta_n)$.
Then,
$$
P_{\Sen,\sigma}(T) = \prod_{i=1}^n (T-\wt_\sigma(\delta_i)) \in A[T]
$$
factors naturally.
Suppose there exists a nonempty subset $I \subset \{1,\dots,n\}$ such that 
$$
Q_I(T):= \prod_{m\in I} (T-\wt_\sigma(\delta_m)), \quad S_I(T):= \prod_{m \in \{1,\dots,n\}\setminus I} (T-\wt_\sigma(\delta_m)),
$$
are comaximal in $A[T]$.
Then, applying Theorem \ref{thm:Wu_3.16}(i) (resp. Theorem \ref{thm:Wu_3.16}(ii)) to
$$
P_{\Sen,\sigma}(T) = Q_I(T) S_I(T)
$$
produces a $(\varphi,\Gamma_K)$-module $D' \subset D$ (resp. $D''\supset D$) that is again trianguline. 
\end{prop}
\begin{proof}
We reduce to the case where all the subquotients $\Fil^i(D)/\Fil^{i-1}(D)$ are free of rank 1 over $\cR_{K,A}$ as follows.
Choose a finite admissible cover $\{\Sp(A_l)\}_{l=1}^r$ of $\Sp(A)$ such that the line bundles
$
\cL_i = H^0(\gr^i(\Fil^\bullet(D))(\delta_i^{-1})) \cong \Hom_{\varphi,\Gamma_K}(\cR_{K,A}(\delta_i),\gr^i(\Fil^\bullet(D)))
$
are trivial on each $\Sp(A_l)$.
Below we construct  $(\varphi,\Gamma_K)$-submodules $(D\otimes_A A_l)'$ with triangulations $\Fil^\bullet((D\otimes_A A_l)')$ satisfying the condition\footnote{That is, the $\sigma$-Sen polynomial of $\Fil^i((D\otimes_A A_l)')$ is the product of $Q_{i,I}(T-1)=\prod_{j\leq i, j\in I}(T-(\wt_\sigma(\delta_j)+1))$ and $S_{i,I}(T)=\prod_{j\leq i, j\notin I}(T-\wt_\sigma(\delta_j))$, and $t_\sigma \Fil^i(D\otimes_A A_l) \subset \Fil^i((D\otimes_A A_l)')\subset \Fil^i(D\otimes_A A_l)$, for each $i$.} on Sen polynomials, for $1\leq l\leq r$.
By the uniqueness part of Theorem \ref{thm:Wu_3.16} applied to the $(\varphi,\Gamma_K)$-module $\Fil^i(D)$ of rank $i$, the $(\varphi,\Gamma_K)$-modules $\{\Fil^i((D\otimes_{A}A_l)')\}_{l=1}^r$ over $\coprod_{l=1}^r \Sp(A_l)$ glue to the $(\varphi,\Gamma_K)$-submodule $(\Fil^i(D))' \subset \Fil^i(D)$ over $\Sp(A)$ given by the existence part of Theorem \ref{thm:Wu_3.16}, for $1\leq i\leq n$. 
Hence, the submodule $D'\subset D$ is trianguline with triangulation $\Fil^i(D') = (\Fil^i(D))'$ for $1\leq i\leq n$.
The same argument holds for $D''\supset D$.

Assume from now on that all the subquotients $\Fil^i(D)/\Fil^{i-1}(D) \cong \cR_{K,A}(\delta_i)$ are free.
For a $(\varphi,\Gamma_K)$-module $R_{K,A}(\delta)$ of rank $1$, $D_{\dif}^+(\delta)$ and $D_{\Sen}(\delta)= D_\dif^+(\delta)/t$ are free of rank $1$.
By the functoriality of $D_{\dif}^+$, $D_{\dif}^+(D)$ is free of rank $n$ over $(K_\infty \otimes_{\bbQ_p} A)[[t]]$.
The triangulation $\Fil^\bullet(D)$ induces a basis $\{e_1,\dots, e_n\}$ of $D_{\dif}^+(D)$ such that $\mathrm{Span}(e_1,\dots, e_m) = D_\dif^+(\Fil^m(D))$ for $1\leq m\leq n$.
Using $D_\dif^+(D) = \prod_{\tau\in \Sigma_K} D^+_{\dif,\tau}(D)$ and $D_{\Sen}(D) = D_\dif^+(D)/t = \prod_{\tau\in \Sigma_K} (D^+_{\dif,\tau}(D)/t) = \prod_{\tau\in \Sigma_K} (D^+_{\dif,\tau}(D)/t_\tau) = \prod_{\tau\in \Sigma_K} D_{\Sen,\tau}(D)$, we get an induced basis $\{e_{1,\tau},\dots,e_{n,\tau}\}$ of $D^+_{\dif,\tau}(D)$, which modulo $t_\tau$ reduces to a basis $\{\overline{e}_{1,\tau},\dots, \overline{e}_{n,\tau}\}$ of $D^+_{\Sen,\tau}(D)$, for each $\tau\in \Sigma_K$.

By the comaximal factorization $P_{\Sen,\sigma}(T)=Q_I(T)S_I(T)$, we have the decomposition
$$
D_{\Sen,\sigma}(D) = \ker(Q_I(\Theta_\sigma)|D_{\Sen,\sigma}(D)) \oplus \ker(S_I(\Theta_\sigma)|D_{\Sen,\sigma}(D)).
$$
Thanks to the triangulation, we have the decomposition of $\Theta_\sigma|_{\Fil^m(D)}$ on the filtration steps
\begin{align*}
D_{\Sen,\sigma}(\Fil^m D) 
=
\ker(Q_I(\Theta_\sigma)|D_{\Sen,\sigma}(\Fil^m D)) \oplus \ker(S_I(\Theta_\sigma)|D_{\Sen,\sigma}(\Fil^m D))
\end{align*}
compatible with the filtration, which implies that
$$
\Fil^m D_{\Sen,\sigma}(D) = \Fil^m  \ker(Q_I(\Theta_\sigma)|D_{\Sen,\sigma}(D)) \oplus \Fil^m \ker(S_I(\Theta_\sigma)|D_{\Sen,\sigma}(D))
$$
for each $1\leq m \leq n$.
We can choose a new basis $\{e_1',\dots,e_n'\}$ of $D_{\dif}^+(D)$ such that
\begin{enumerate}
\item 
$\mathrm{Span}(e_1',\dots, e_m') = D_\dif^+(\Fil^m(D))$ for $1\leq m\leq n$, and  
\item
we have
$$
\Fil^m (\ker(Q_I(\Theta_\sigma)|D_{\Sen,\sigma}(D))) = \Span\set{\overline{e}_{i,\sigma}'}{i\in \{1,\dots,m\} \cap I},
$$
$$
\Fil^m (\ker(S_I(\Theta_\sigma)|D_{\Sen,\sigma}(D))) = \Span\set{\overline{e}_{i,\sigma}'}{i\in \{1,\dots,m\}\setminus I},
$$
for each $1\leq m\leq n$.
\end{enumerate}
Indeed, we obtain $\{e_1',\dots,e_n'\}$ from $\{e_{1},\dots,e_{n}\}$ by induction on $n$.
When $n=1$, $e_1$ reduces to an eigenvector $\overline{e}_{1,\sigma}$ of $\Theta_\sigma$ with eigenvalue being the unique root of $P_{\Sen,\Fil^1(D)_\sigma}(T)$.
Thus, we can set $e_1':=e_1$, which satisfies (1) and (2) for the $(\varphi,\Gamma_K)$-module $\Fil^1(D)$ of rank $1$.
Suppose that, for $1<m\leq n$, we have chosen $\{e_1',\dots,e_{m-1}'\}$ that satisfies (1) and (2) for the $(\varphi,\Gamma_K)$-module $\Fil^l(D)$ of rank $l$, for all $l\leq m-1$.
For $\Fil^m(D)$, the basis $\{e_1',\dots,e_{m-1}',e_m\}$ satisfies (1) but not necessarily (2).
Choosing $B,C\in A[T]$ such that $BQ_I+CS_I=1$, we have
$$
\overline{e}_{m,\sigma} = (BQ_I+CS_I)(\Theta_\sigma)(\overline{e}_{m,\sigma}) = \underbrace{BQ_I(\Theta_\sigma)(\overline{e}_{m,\sigma})}_{\in \ker(S_I(\Theta_\sigma))} + \underbrace{CS_I(\Theta_\sigma)(\overline{e}_{m,\sigma})}_{\in \ker(Q_I(\Theta_\sigma))} ,
$$
and there are two cases for us to discuss:
\begin{itemize}
\item 
If $m\in I$, then $\Fil^m (\ker(S_I(\Theta_\sigma)|D_{\Sen,\sigma}(D)))= \Fil^{m-1} (\ker(S_I(\Theta_\sigma)|D_{\Sen,\sigma}(D)))$, which has a basis given by $\set{\overline{e}_{i,\sigma}'}{i\in\{1,\dots,m-1\}\setminus I}$ by the induction hypothesis.
Thus,
$$
BQ_I(\Theta_\sigma)(\overline{e}_{m,\sigma}) = \sum_{i\in\{1,\dots,m-1\}\setminus I}a_{i,\sigma} \overline{e}_{i,\sigma}'
$$
belongs to
$\Span_{K_\infty \otimes_{K,\sigma} A}\set{\overline{e}_{i,\sigma}'}{i\in\{1,\dots,m-1\}\setminus I}$, for some $a_{i,\sigma} \in K_\infty \otimes_{K,\sigma} A$.
Let $a_i \in K_\infty \otimes_{\bbQ_p}A \cong \prod_{\tau\in \Sigma_K} K_\infty \otimes_{K,\tau}A$ be the element $(0,\dots,0,a_{i,\sigma},0,\dots, 0)$, and set 
$$
e_m':= e_m - \sum_{i\in\{1,\dots,m-1\}\setminus I} a_i e_i'.
$$
Then, $\{e_1',\dots,e_m'\}$ is a basis of $D^{+}_{\dif}(\Fil^m(D))$, and 
$
\overline{e}_{m,\sigma}' = \overline{e}_{m,\sigma} - BQ_I(\Theta_\sigma)(\overline{e}_{m,\sigma}) = CS_I(\Theta_\sigma)(\overline{e}_{m,\sigma})\in \Fil^m (\ker(Q_I(\Theta_\sigma)|D_{\Sen,\sigma}(D)))$.
Hence, (1) and (2) are satisfied.
\item 
If $m\notin I$, then $\Fil^m (\ker(Q_I(\Theta_\sigma)|D_{\Sen,\sigma}(D)))= \Fil^{m-1} (\ker(Q_I(\Theta_\sigma)|D_{\Sen,\sigma}(D)))$, which has a basis given by $\set{\overline{e}_{i,\sigma}'}{i\in\{1,\dots,m-1\}\cap I}$ by the induction hypothesis.
Similar to the above case when $m\in I$, we can subtract from $e_m$ a suitable $(K_\infty\otimes_{\bbQ_p} A)$-linear combination of $e_1',\dots,e_{m-1}'$ to get $e_m'$ such that $\{e_1',\dots,e_m'\}$ satisfies (1) and (2).
\end{itemize}
By induction, we conclude the existence of a basis of $D_\dif^+(D)$ satisfying (1) and (2). 
Without loss of generality, we may assume that our chosen basis $\{e_1,\dots,e_n\}$ satisfies (1) and (2).

By the construction of $D''\supset D$, it suffices to consider the submodule $D'\subset D$.
By Wu's construction given in the proof of Theorem \ref{thm:Wu_3.16}, the submodule $D'$ has $D_{\dif,\tau}^+(D') = D_{\dif,\tau}^+(D)$ for all $\tau \in \Sigma_K\setminus \{\sigma\}$, and using our modified basis, it is easy to see that
$$
D_{\dif,\sigma}^+(D') = \ker(D_{\dif,\sigma}^{+}(D) \to D_{\dif,\sigma}^{+}(D)/t \to \ker(Q_I(\Theta_\sigma)|D_{\Sen,\sigma}(D))) = \Span\{e'_1,\dots,e_n'\},
$$
where $e_m' := t_\sigma e_m$ if $m\in I$ and $e_m' := e_m$ if $m\notin I$.
Since $\mathrm{Span}(e_1,\dots, e_m) = D_\dif^+(\Fil^m(D))$ are $\Gamma_K$-stable for $1\leq m \leq n$, the only way for $D_{\dif,\sigma}^+(D')$ to be $\Gamma_K$-stable is that $\Span(e_1',\dots,e_m')$ are $\Gamma_K$-stable for $1\leq m \leq n$.
Using the triangulation $\Fil^\bullet(D)$ on $D$, Beauville-Laszlo gluing \cite[Proposition A.3]{Wu_translation} implies that the $\Gamma_K$-stable flag on $D_\dif^+(D')$
$$\{\Span(e_1',\dots,e_m')| 1\leq m \leq n\} $$
determines a triangulation $\Fil^\bullet(D')$ on $D'$ such that $D_{\dif}^+(\Fil^m (D'))=\Span(e_1',\dots,e_m')$ and
$$
\Fil^{m}(D')/\Fil^{m-1}(D') \cong 
\begin{cases}
\cR_{K,A}(x_\sigma \delta_m) &\text{ if } m\in I,\\
\cR_{K,A}(\delta_m) &\text{ if } m\notin I,
\end{cases}
$$
for all $1\leq m \leq n$.
\end{proof}

\begin{thm}
\label{thm:comparing_Wu}
\begin{enumerate}[(i)]
\item 
Let $D\in \fX_n^{\sigma\text{-}\unif,i}(A)$. 
Tate-fpqc locally, choose any triangulation $\Fil^\bullet(D)$ on $D$ over $\cR_{K,A}$ with parameters $\delta = (\delta_1,\dots,\delta_n)$ so that 
$$
\Fil^i(D)/\Fil^{i-1}(D) \cong \cR_{K,A}(\delta_i) \otimes_{\cO_{\Sp(A)}} \cL_i
$$
for some line bundles $\cL_i$ on $\Sp(A)$.
Over the Tate-fpqc cover, $p_{i,\sigma}(D,\Fil^\bullet(D))$ equals the unique submodule of $D$ defined in Theorem \ref{thm:Wu_3.16}(i) and is thus independent of the choice of local triangulation $\Fil^\bullet(D)$.
The submodule $p_{i,\sigma}(D,\Fil^\bullet(D))$ descends to the unique submodule $D'$ of $D$ over $\Sp(A)$ defined in Theorem \ref{thm:Wu_3.16}(i).
\item
Denote the operation $D\mapsto D''$ in Theorem \ref{thm:Wu_3.16}(ii) by $q_{i,\sigma}$.
Then for $D\in \fX_n^{\sigma\mh\unif,i,1}(A) = \fX_n^{\sigma\mh\unif,i,[0,1]}(A)$, we have $p_{i,\sigma}(D) \in \fX^{\sigma\mh\unif,i,-1}_n(A) = \fX_n^{\sigma\mh\unif,i,[-1,0]}(A)$.
\item
For any subset $S\subset\Sigma_K$, and subsets $I_\sigma \subset \{1,\dots, n\}$ for $\sigma\in S$, and
${\bf k } = (k_{i,\sigma})_{\sigma\in S, i\in I_\sigma} \in \bbN^{\sum_{\sigma\in S}|I_\sigma|}$ any tuple of nonnegative integers, there is an isomorphism
$$
p_{\bf k} : \fX_n^{S\text{-}\unif, I, {\bf k}} \longrightarrow \fX_n^{S\text{-}\unif, I, {\bf -k}}
$$
given by the composite of elements from $\set{p_{i,\sigma}}{1\leq i\leq n, \sigma\in \Sigma_K}$ in any order in which $p_{i,\sigma}$ appears with multiplicity $k_{i,\sigma}$, whose inverse is given by
$$
q_{\bf k}: \fX_n^{S\mh\unif,I,-{\bf k}} \longrightarrow \fX_n^{S\mh\unif,I,{\bf k}}
$$
that is defined as the composite of elements from $\set{q_{i,\sigma}}{1\leq i\leq n, \sigma\in \Sigma_K}$ in any order in which $q_{i,\sigma}$ appears with multiplicity $k_{i,\sigma}$.
\end{enumerate}
\end{thm}
\begin{proof}
For any $(\varphi,\Gamma_K)$-module $D$ over $\Sp(A)$ that belongs to the stackification of the image of $\beta_B: \fX_B\to \fX_n$, it becomes trianguline after pullback to a Tate-fpqc cover $\{X_j \to \Sp(A)\}_{j\in J}$.
The coverings in the Tate-fpqc topology are generated by the usual (admissible) Tate coverings and the
morphisms $\Sp(B) \to \Sp(A)$ of rigid spaces for faithfully flat maps $A\to B$ of affinoid
$E$-algebras, cf. \cite[\S5.1.7]{EGH_CatLLC}.
Since $\Sp(A)$ is quasi-compact, $\{X_j \to \Sp(A)\}_{j\in J}$ can be refined to an fpqc cover $\Sp(B)\twoheadrightarrow \Sp(A)$.
In the proof of (i)-(iii) below, we implicitly take an fpqc morphism $f: \Sp(B) \twoheadrightarrow \Sp(A)$ such that $f^*D = D \widehat{\otimes}_{A} B$ is trianguline and give our constructions of $p_{i,\sigma}$ and $q_{i,\sigma}$ over $\Sp(B)$.
For any point $y\in \Sp(B)$ and $z:= f(y)\in \Sp(A)$, the fibers satisfy $(f^*D)_{k(y)} = D_{z} \otimes_{k(z)} k(y)$.
Since $\fX_n$ satisfies Tate-fpqc descent, our constructions on $f^*D$ over $\Sp(B)$ descend to  constructions on $D$ over $\Sp(A)$, because the uniqueness part of Theorem \ref{thm:Wu_3.16} implies that $p_{i,\sigma}$ and $q_{i,\sigma}$ respect the isomorphisms in the descent data. 
Thus, $p_{i,\sigma}(f^*D) \subset f^*D$ descends to $p_{i,\sigma}(D)\subset D$ in (i), and $q_{i,\sigma}(f^*D) \supset f^*D$ descends to $q_{i,\sigma}(D)\supset D$ in (ii).

\vskip8pt

\noindent
(i) 
Passing to a finite admissible cover $\{\Sp(A_l)\}_{l=1}^r$ of $\Sp(A)$, we assume the $\cL_i$ are trivial.
By the definition of $U_{i}$, we have that
$$
\{ \wt_\sigma(\delta_{1,z}),\dots, \wt_\sigma(\delta_{n-i,z})\} \cap \{\wt_\sigma(\delta_{n-i+1,z}),\dots, \wt_\sigma(\delta_{n,z})\} = \emptyset
$$
inside the residue field $k(z)$ for all $z\in \Sp(A)$.
Thus, for $Q_{i,\sigma}(T):=
\prod_{j=n-i+1}^{n} (T-\wt_\sigma(\delta_j))$ and $S_{i,\sigma}(T):=\prod_{j={1}}^{n-i} (T-\wt_\sigma(\delta_j))$, we have $(Q_{i,\sigma}(T),S_{i,\sigma}(T)) = (1)$ and a factorization 
$$
P_{\Sen,\sigma}(T) = Q_{i,\sigma}(T)S_{i,\sigma}(T) \in A[T],
$$
which is independent of the choice of $\Fil^\bullet(D)$:
indeed, if $\{\delta'_m\}_{m=1}^n$ are the parameters attached to another triangulation of $D$ over $\cR_{K,A}$, then we get another comaximal factorization
$$
P_{\Sen,\sigma}(T) = Q'_{i,\sigma}(T)S'_{i,\sigma}(T) \in A[T],
$$
where we similarly put $Q'_{i,\sigma}(T):=\prod_{j=n-i+1}^{n} (T-\wt_\sigma(\delta'_j))$ and $S'_{i,\sigma}(T):=\prod_{j=1}^{n-i} (T-\wt_\sigma(\delta'_j))$.
Moreover, we have $(Q'_{i,\sigma}(T),S_{i,\sigma}(T)) = (Q_{i,\sigma}(T), S'_{i,\sigma}(T)) = 1$ by Lemma \ref{lem:coprime-polynomials}(iii). 
We conclude that $Q_{i,\sigma}(T)=Q'_{i,\sigma}(T)$ and $S_{i,\sigma}(T)=S'_{i,\sigma}(T)$ by Lemma \ref{lem:coprime-polynomials}(i).

Let $\{e_1,\dots, e_{n-i},\dots, e_n\}$ be a basis of $D_{\dif}^+(D)$ as in the proof of Proposition \ref{prop:Wu_trianguline} with 
$$
I := \set{k}{n-i+1\leq k\leq n} \subset \{1,\dots, n\}
$$
such that for the splitting of $D_{\Sen,\sigma}(D)$ as
$$
\ker(Q(\Theta_\sigma)|D_{\Sen,\sigma}(D)) \oplus \ker(S(\Theta_\sigma)|D_{\Sen,\sigma}(D)) \cong D_{\Sen,\sigma}(D/\Fil^{n-i}(D)) \oplus D_{\Sen,\sigma}(\Fil^{n-i}(D))
$$
according to the derivative $\Theta_\sigma$ of the $\Gamma_K$-action, the mod-$t_\sigma$ reduction $\{\overline{e_1},\dots, \overline{e_{n-i}}\}$ is a basis of $\ker(S(\Theta_\sigma)|D_{\Sen,\sigma}(D))$ and $\{\overline{e_{n-i+1}},\dots, \overline{e_n}\}$ is a basis of $\ker(Q(\Theta_\sigma)|D_{\Sen,\sigma}(D))$.
From the proof of Proposition \ref{prop:Wu_trianguline}, we see that the submodule $N$ with basis 
$$
\{e_1,\dots, e_{n-i}, t_\sigma e_{n-i+1},\dots, t_\sigma e_n\}
$$ 
in $D_\dif^+(D)$ corresponds to the module $M$ given by Wu in Theorem \ref{thm:Wu_3.16}(i).
Thus, $p_{i,\sigma}(D,\Fil^\bullet(D))$ corresponds to $M$ since the localization of $p_{i,\sigma}(D,\Fil^\bullet(D))$ clearly has the same basis 
$$
\{e_1,\dots, e_{n-i}, t_\sigma e_{n-i+1},\dots, t_\sigma e_n\}
$$
by the definition of pullback, and (i) is proven.

\vskip8pt

\noindent
(ii)
Keep the notations in the proof of (i), in particular the basis $e_1,\dots,e_n$.
By (i), $q_{i,\sigma}(D)$ can be computed using any triangulation $\Fil^\bullet(D)$ on $D$ with respect to which our basis is chosen, because all triangulations induce the same factorization of $\sigma$-Sen polynomial.
We know that $q_{i,\sigma}(D)$ is trianguline by Proposition \ref{prop:Wu_trianguline}.
Let $h_1,\dots,h_n$ be the ordered $\sigma$-Sen weights of $D$ with respect to any triangulation $\Fil^\bullet(D)$.
Then,
$$
h'_m = \begin{cases}
h_m &\text{ if } 1\leq m \leq n-i,\\
h_m+1 &\text{ if } n-i+1\leq m \leq n,
\end{cases}
$$
are the ordered $\sigma$-Sen weights of an induced triangulation on the pullback $D':=p_{i,\sigma}(D)$.
Since we assume $D\in \fX_n^{\sigma\mh\unif,i,1}(A) = \fX_n^{\sigma\mh\unif,i,[0,1]}(A)$, we have
$$
h_j' \not\equiv h_k' \pmod{\frm_z}, \quad  h_j' \not\equiv h_k'-1 \pmod{\frm_z}
$$
for all $z\in \Sp(A)$ with the corresponding maximal ideal $\frm_z$ and for all $1\leq j\leq n-i< k \leq n$.
In particular, we can apply $q_{i,\sigma}$ to $D'$.
It remains to show that $D'$ is $\sigma$-weight-uniform.
Assume, for the sake of contradiction, that there exists another triangulation on $D'$ such that at some point $z\in \Sp(A)$, the induced ordering on $\sigma$-Sen weights is 
$$
(h'_{w(1),z},\dots, h'_{w(n),z})
$$
for some $w\in S_n$ such that for some $1\leq m_0\leq n$, $h'_{w(m_0),z}\neq h'_{m_0,z}\in k(z)$.
By Proposition \ref{prop:Wu_trianguline}, this other triangulation on $D'$ induces on $q_{i,\sigma}(D') = q_{i,\sigma}(p_{i,\sigma}(D))= D$ a triangulation whose ordered $\sigma$-Sen weights are $(h''_{1},\dots, h''_{n})$ with
$$
h_m'' = 
\begin{cases}
h'_{w(m)}-1 = h_{w(m)} &\text{ if } w(m)\in I := \set{k}{n-i+1\leq k \leq n}, \\
h'_{w(m)} = h_{w(m)}&\text{ if } w(m)\notin I.
\end{cases}
$$
\underline{Case (a)}.
If $\{m_0,w(m_0)\}$ is contained in either $I$ or $I^c$, then $h'_{w(m_0),z}\neq h'_{m_0,z}\in k(z)$ implies 
$$
h_{m_0,z}'' = h_{w(m_0),z} \neq h_{m_0,z} \in k(z),
$$
contrary to the assumption that $D$ is $\sigma$-weight-uniform.
\\
\underline{Case (b)}.
Otherwise, either $m_0\in I$ but $w(m_0) \notin I$, or $m_0\notin I$ but $w(m_0) \in I$. 
Then we have
$$
h_{m_0,z}'' = h_{w(m_0),z} \neq h_{m_0,z} \in k(z)
$$ 
because $D\in \fX_n^{\sigma\mh\unif,i}(A)$, contrary to the assumption that $D$ is $\sigma$-weight-uniform.

Hence, we conclude that $p_{i,\sigma}: \fX_n^{\sigma\mh\unif,i,1} \to \fX_n^{\sigma\mh\unif, i,-1}$ preserves $\sigma$-weight-uniformity.
\vskip8pt

\noindent
(iii)
By Lemma \ref{lem:p_i_commutes_with_p_j}, the independence (i) of triangulations for the pullback interpretation of $p_{i,\sigma}$ under the given assumption on weights, and (ii) that $p_{i,\sigma}$ preserves weight-uniformity under the given assumption on weights, it follows that any way of composing $p_{i,\sigma}$ with multiplicity $k_{i,\sigma}$ as $i$ varies in $\{1,\dots, n\}$ and $\sigma$ varies in $\Sigma_K$ produces the same result, which is the morphism $p_{\bf k}$.
Since $q_{i,\sigma}$ is the inverse of $p_{i,\sigma}$ by Theorem \ref{thm:Wu_3.16}(iii), it follows that we can define $q_{\bf k}$ by composing $q_{i,\sigma}$'s with multiplicity $k_{i,\sigma}$ in any order.
It is clear that $p_{\bf k}$ and $q_{\bf k}$ are mutual inverses.
\end{proof}

\subsubsection{}
To summarize, we have $p_{i,\sigma}: \fX_B \to \fX_B$, which is often invertible by Theorem \ref{thm:p_k-invertible}.
By Theorem \ref{thm:comparing_Wu}, this induces a ``change of weights" morphism of stacks 
$$
p_{i,\sigma} : \fX_{n}^{\sigma\mh\unif,i} \to \fX_n,
$$
which, if the ``regularity" of the Sen weights is unchanged after the weight change, is invertible:
$$
p_{i,\sigma} : \fX_{n}^{\sigma\mh\unif,i,1} \xrightarrow{\sim} \fX_{n}^{\sigma\mh\unif,i,-1}
$$
by Theorem \ref{thm:comparing_Wu}(iii).
Moreover, we know by Lemma \ref{lem:p_i_commutes_with_p_j} that for any $1\leq i,j\leq n$ and $\tau,\sigma \in \Sigma_K$, 
$$
(p_{j,\tau} \circ p_{i,\sigma})(D) = (p_{i,\sigma} \circ p_{j,\tau})(D)
$$
if the pullbacks are independent of the choice of triangulation.
More generally, one obtains
$$
p_{\bf k}: \fX_n^{\unif, {\bf k}} \xrightarrow{\sim} \fX_n^{\unif, {\bf -k}}
$$
as long as the ``regularity" of the Sen weights is unchanged after the change.
See \cite[Example 3.18]{Wu_translation} for the general (non-invertible) case of changing the Sen weights on $\fX_n$ in a similar way.

One family of ``nice'' $D$ are those non-critical crystabelline $(\varphi,\Gamma_K)$-modules $D \in \fX_n^\unif(A)$ of regular Sen weights.
When $A=E$ is a finite extension of $\bbQ_p$, we verify in the last section that $p_{i,\sigma}$ on $D$ corresponds to the translation functor on $\pi_{\fs}(D)$, for Ding's construction $\pi_{\fs}(D)$.

\subsection{Étaleness}
According to Fontaine and others, there is an equivalence of categories $D_{\rig}$ between continuous $E$-linear $G(\overline{K}/K)$-representations and étale $(\varphi,\Gamma_K)$-modules over $\cR_{K,E}$. 
We study in this subsection how the pullback operations interplay with étaleness on certain very generic or non-critical crystabelline $(\varphi,\Gamma_K)$-modules.

In this subsection, we write $\cR:=\cR_{K,E}$.
Let $\cC_E$ be the category of local Artinian $E$-algebras $(A,\frm_A)$ with residue field $A/\frm_A \cong E$.
By \cite[Lemma 2.2.5]{BeCh_families}, our discussion also applies to $A$-linear Galois representations and $(\varphi,\Gamma_K)$-modules over $\cR_{K,A}$ for $A\in \cC_E$.

Let us recall Kedlaya's theory of slopes.
For a $(\varphi,\Gamma_K)$-module $D$ of rank $n$ over $\cR$, we define its {\bf degree} by $\deg(D):=\deg(\wedge^n D)$, which is the $p$-adic valuation of a ``$\varphi$-eigenvalue'' on $\wedge^n D$.
Then, the {\bf slope} of $D$ is defined by $\mu(D):= \deg(D)/\rank(D)$.
We say that $D$ is {\bf semistable} if for any finite free $\varphi$-submodule $M$ of $D$ satisfying $\varphi^*M\cong M$, one has $\mu(M)\geq \mu(D)$.
Then, $D$ is called {\bf étale} if it is semistable of slope zero, cf. \cite[p.8]{Liu_cohomology}.
\begin{rem}
\label{rem:slope_basic}
The following properties clearly hold for $(\varphi,\Gamma_K)$-modules.
\begin{enumerate}[(i)]
\item 
If $0 \rightarrow D' \rightarrow D \rightarrow D'' \rightarrow 0$ is exact, then  $\operatorname{deg}(D)=\operatorname{deg}\left(D'\right)+\operatorname{deg}\left(D''\right)$.
\item
One has $\mu\left(D_1 \otimes D_2\right)=\mu\left(D_1\right)+\mu\left(D_2\right)$.
\item 
We have $\operatorname{deg}\left(D^{\vee}\right)=-\operatorname{deg}(D)$ and $\mu\left(D^{\vee}\right)=-\mu(D)$.
\end{enumerate}
\end{rem}

\begin{ex}
\label{ex:slopes}
By \cite[Construction 6.2.4]{KPX_cohomology}, for any continuous character $\delta: K^\times \to E^\times$, if we factorize ${\delta}={\delta}^{\unr}{\delta}^{\wt}$ such that ${\delta}^{\unr}(\pi_K)=\delta(\pi_K)$ and ${\delta}^\wt|_{\cO_K^\times} = {\delta}|_{\cO_K^\times}$, then 
$
\cR(\delta) = \cR(\delta^{\unr}) \otimes_\cR \cR(\delta^{\wt}).
$
Since $\cR(\delta^{\wt})$ is étale and $\varphi^f=\delta(\pi_K)$ on $\cR(\delta^{\unr}) = D_{f,\delta(\pi_K)}\otimes_{K_0\otimes E} \cR$, where $D_{f,\delta(\pi_K)}=K_0\otimes_{\bbQ_p} E$ carries a Frobenius semilinear endomorphism $\varphi$ such that $\varphi^f = \delta(\pi_K)$, cf. \cite[Lemma 6.2.3]{KPX_cohomology}, it follows from Remark \ref{rem:slope_basic}(ii) that
$$
\mu(\cR(\delta)) = \deg(\cR(\delta)) = \frac{1}{f}v_p(\delta(\pi_K)),
$$
where $v_p$ is the $p$-adic valuation on $\overline{\bbQ_p}$ normalized by the condition $v_p(p)=1$.

\end{ex}
The following is the slope filtration theorem by Kedlaya.
\begin{thm}[{\cite[Theorem 2.4]{Liu_cohomology}}]
\label{thm:slope_filtration}
Every $(\varphi,\Gamma_K)$-module $D$ over $\cR_{K,E}$ admits a unique filtration $0=D_0 \subset$ $D_1 \subset \cdots \subset D_l=D$ by saturated $(\varphi,\Gamma_K)$-submodules whose successive quotients are semistable with increasing slopes $\mu(D_1 / D_0)<\cdots<\mu(D_{l} / D_{l-1})$.
\end{thm}
In view of Theorem \ref{thm:slope_filtration} and Remark \ref{rem:slope_basic}(i), we deduce
\begin{cor}
\label{cor:etale_criterion}
A $(\varphi,\Gamma_K)$-module $D$ over $\cR_{K,E}$ is étale if and only if $\mu(D)=0$ and it does not contain any saturated $(\varphi,\Gamma_K)$-submodule of strictly negative slope.
\end{cor}

\subsubsection{}
Suppose $(D,\Fil^\bullet(D))$ is strongly non-split of parameters in $\cT^n_\circ(E)$. 
By Corollaries \ref{cor:generic-unique} and \ref{cor:etale_criterion}, $D$ is étale if and only if $\Fil^i(D)$ has non-negative slopes and $D$ has slope zero.
Thus, we obtain the following proposition.
\begin{prop}
\label{prop:vgen_etale}
\begin{enumerate}[(i)]
\item 
If $(D,\Fil^\bullet(D))$ is strongly non-split of parameter $(\delta_1,\dots,\delta_n) \in \cT^n_\circ(E)$, then $D$ is étale if and only if 
$\sum_{i=1}^m v_p(\delta_i(\pi_K)) \geq 0$ for $1\leq m < n$ and $\sum_{i=1}^n v_p(\delta_i(\pi_K)) = 0$.
\item 
If $(D,\Fil^\bullet(D))$ is strongly non-split of parameter $(\delta_1,\dots,\delta_n)\in\cT^n_\circ(E)$ and étale, then $p_{j,\sigma}(D)$ is étale up to twist if and only if  
$$
\begin{cases}
\displaystyle
\sum_{i=1}^{m} v_p(\delta_i(\pi_K)) \geq \frac{mj}{ne} & \quad 1\leq m \leq n-j, \\
\displaystyle
\sum_{i=1}^{m} v_p(\delta_i(\pi_K)) + \frac{m-(n-j)}{e} \geq \frac{mj}{ne} &\quad   n-j<m<n.
\end{cases}
$$
\end{enumerate}
\end{prop}

\begin{proof}
\begin{enumerate}[(i)]
\item 
By Remark \ref{rem:slope_basic}(i) and the definition of slope, 
$$
\mu(\Fil^i(D)) =  \deg(\Fil^i(D))/\rank(\Fil^i(D)) = \frac{1}{i}\sum_{j=1}^i \deg(\cR(\delta_j)).
$$
By Example \ref{ex:slopes}, we see that $\displaystyle \mu(\Fil^i(D)) = \frac{1}{fi}\sum_{j=1}^i v_p(\delta_j(\pi_K))$, for all $i$.
\item 
The new parameter of $p_{j,\sigma}(D)$ is $(\delta_1,\dots,\delta_{n-j}, x_\sigma \delta_{n-j+1}, x_\sigma \delta_{n-j+2},\dots, x_\sigma \delta_{n})$, and 
$$
v_p(x_\sigma(\pi_K)) = v_p(\sigma(\pi_K)) = v_p(\pi_K) = 1/e.
$$
Let $\chi:K^\times \to E^\times$ be any continuous character.
The twist $p_{j,\sigma}(D)(\chi)$ is again non-split and very generic of parameter $(\delta_1\chi,\dots, \delta_{n-j}\chi, x_\sigma \delta_{n-j+1}\chi, x_\sigma \delta_{n-j+2}\chi,\dots, x_\sigma \delta_{n}\chi)$, which is étale if and only if we have
$$
\begin{cases}
\sum_{i=1}^m v_p(\delta_i(\pi_K)) + m v_p(\chi(\pi_K)) \geq 0 &\text{ if } 1\leq m \leq n-j, \\
\sum_{i=1}^m v_p(\delta_i(\pi_K)) + (m-(n-j))v_p(\pi_K) + m v_p(\chi(\pi_K))\geq 0 &\text{ if } n-j < m \leq n-1, \\
\sum_{i=1}^{n} v_p(\delta_i(\pi_K)) + j v_p(\pi_K) + n v_p(\chi(\pi_K)) = 0
&\text{ if } m = n.
\end{cases}
$$
by (i) above.
Since $D$ is étale, ${\sum_{i=1}^{n} v_p(\delta_i(\pi_K))} + j v_p(\pi_K) + n v_p(\chi(\pi_K)) = 0$ implies that $v_p(\chi(\pi_K)) = -j/(ne)$,
which allows us to complete the proof by taking unramified $\chi$.
\qedhere
\end{enumerate}
\end{proof}

\subsubsection{}
In this paragraph, we consider the étaleness for $D\in \PGnc$.
By Corollary \ref{cor:etale_criterion} and Theorem \ref{thm:D_cris_WD_equivalence}, all saturated $(\varphi,\Gamma_K)$-submodules of $D$ are $\Fil_w^i(D)=D_{w,i}$ as in Proposition \ref{prop:crys_generic_triang} for $0\leq i \leq n$ and $w\in S_n$.
Let $\{\alpha_i :=\phi_i(\pi_K)\}_{1\leq i\leq n}$ be the $\varphi^f$-eigenvalues on $\DcrisKm(D)$.
Then, by Proposition \ref{prop:crys_generic_triang}(ii) and Example \ref{ex:slopes}, we have 
\begin{equation}
\label{eqn:crystabelline_nc_slope}
v_p(\delta_{w,i}(\pi_K)) = v_p\left(\prod_{\sigma} \sigma(\pi_K)^{h_{i,\sigma}} \phi_{w(i)}(\pi_K)\right) = v_p(\alpha_{w(i)})+\sum_\sigma h_{i,\sigma} \frac{1}{e}.
\end{equation}
\begin{prop}
\label{prop:crys_etale}
Suppose $D\in \PGnc$.
Let $\tau\in S_n$ be a refinement such that
$$
v_p(\alpha_{\tau(1)}) \leq v_p(\alpha_{\tau(2)}) \leq \dots \leq v_p(\alpha_{\tau(n)}).
$$
Then, $D$ is étale if and only if 
$$
\begin{cases}
\displaystyle
\sum_{i=1}^j v_p(\alpha_{\tau(i)}) \geq -\frac{1}{e}\sum_{\sigma} \sum_{i=1}^j h_{i,\sigma}  &\text{ if } 1\leq j < n, \\
\displaystyle
\sum_{i=1}^n v_p(\alpha_{i}) = -\frac{1}{e}\sum_{\sigma, i}  h_{i,\sigma} &\text{ if } j = n.
\end{cases}
$$
\end{prop}
\begin{proof}
It follows immediately from Corollary \ref{cor:etale_criterion}, our choice of $\tau$, and Equation \eqref{eqn:crystabelline_nc_slope}.
\end{proof}

\begin{rem}
Alternatively, to prove Proposition \ref{prop:crys_etale} we can observe that, by our choice of $\tau$, the $n$ inequalities amount to having that the Hodge filtration on $D_\cris(D)$ is weakly admissible in the sense of \cite[Definition 4.4.3]{Fontaine_p-adiques_semi-stables}, or concretely \cite[Equations (4) and (5)]{BS_functoriality}.
Indeed, by \cite[Proposition 3.1.1.5]{BM-multiplicity} and \cite[Proposition 4.4.9]{Fontaine_p-adiques_semi-stables}, it suffices to check $t_H(D')\geq t_N(D')$ for all {\it $E$-filtered} $(\varphi,G(K_m/K))$-submodules over $K$ of $D_\cris(D)$ with the induced filtration.
But these $D'$ are in bijection with the $n!$ refinements, so their Newton numbers are clear, and their Hodge numbers are clear by the non-criticalness assumption.
A computation similar to that for \cite[Proposition 3.2, $(i)\Rightarrow(ii)$]{BS_functoriality} gives the $n$ inequalities.
\end{rem}

\begin{cor}
Suppose $D\in \fX_n(E)$  satisfies either the hypothesis of Proposition \ref{prop:vgen_etale} or that of Proposition \ref{prop:crys_etale} and is étale of regular Sen weights.
Suppose that $p_{i,\sigma}(D)$ again satisfies the hypothesis of Proposition \ref{prop:vgen_etale} or that of Proposition \ref{prop:crys_etale} and is étale after twist by a character $\chi: K^\times \to E^\times$ (whether this is possible can be checked by Proposition \ref{prop:vgen_etale} or \ref{prop:crys_etale}).

Then, for any local Artinian E-algebra $A\in \cC_E$ of residue field $E$ and for any deformation $D_A$ of $D$ to an element in $\fX_n(A)$, $p_{i,\sigma}(D_A)$ is again étale after twisting by the same character $\chi$.
\end{cor}

\begin{proof}
First of all, by $p_{i,\sigma}(D_A)$ we mean the result of applying Wu's construction (Theorem \ref{thm:Wu_3.16}) to $D_A$, which is applicable and deforms $p_{i,\sigma}(D)$ from $\cR_{K,E}$ to $\cR_{K,A}$ since at the unique closed point $z\in \Sp(A)$, $D_{A,z}=D$ is assumed to have regular Sen weights.

The étaleness of $p_{i,\sigma}(D_A)(\chi)$ follows from the assumption that $p_{i,\sigma}(D)(\chi)$ is étale and the fact that extensions of pure slope-$s$ $\varphi$-modules are pure of slope $s$ by \cite[Lemma 2.2.5]{BeCh_families}.
\end{proof}

\section{Relation with translation functors}
\label{section:translation_functor}

\subsection{Trianguline deformations of $(\varphi,\Gamma_K)$-modules}
We discuss certain deformations of $(\varphi,\Gamma_K)$-modules from $\cR_{K,E}$ to $\cR_{K,A}$ when $A$ is a local Artinian $E$-algebra of residue field $E$.

\begin{defn}
Let $D$ be a $(\varphi,\Gamma_K)$-module of rank $n$ over $\cR_{K,E}$, and let $A\in \cC_E$ be a local Artinian $E$-algebra with maximal ideal $\frm_A$ such that $A/\frm_A\cong E$.
\begin{enumerate}[(i)]
\item 
By a {\bf trianguline deformation} of $D$ to $\cR_{K,A}$ we mean a trianguline $(\varphi,\Gamma_K)$-module $D_A$ over $\cR_{K,A}$ such that $D_A\otimes_A E = D_A/\frm_A D_A \cong D$.
Note that if $\Fil^\bullet(D_A)$ is a triangulation on $D_A$, then we have an induced triangulation $\Fil^\bullet(D_A) \otimes_{A} E$ on $D_A\otimes_A E \cong D$.
\item 
For any given triangulation $\cF^\bullet$ on $D$, an {\bf $\cF^\bullet$-trianguline deformation} of $D$ to $R_{K,A}$ is a pair $(D_A,\Fil^\bullet(D_A))$ consisting of a deformation $D_A$ of $D$ to $\cR_{K,A}$ and a triangulation $\Fil^\bullet(D_A)$ such that $\Fil^\bullet(D_A)\otimes_{A} E = \cF^\bullet$ via the map $D_A \twoheadrightarrow D_A/\frm_A D_A \cong D$.
\end{enumerate}

\end{defn}

The next lemma generalizes \cite[Proposition 2.3.6]{BeCh_families} from the case $K=\bbQ_p$ to arbitrary finite extension $K/\bbQ_p$, with essentially the same argument.
\begin{lem}
\label{lem:unique_deform}
Let $D$ be a $(\varphi,\Gamma_K)$-module of rank $n$ over $\cR_{K,E}$. 
Let $A\in \cC_E$ be a local Artinian $E$-algebra with maximal ideal $\frm_A$ such that $A/\frm_A\cong E$.
\begin{enumerate}[(i)]
\item 
Suppose that for a continuous character $\delta:K^\times \to E^\times$, $D$ has a saturated $(\varphi,\Gamma_K)$-submodule $D_0$ of rank 1 such that $D_0\cong \cR_{K,E}(\delta)$ and $\Hom_{\varphi,\Gamma_K}(\cR_{K,E}(\delta),D/D_0)=0$.
Then for any deformation $D_A$ of $D$ to $\cR_{K,A}$, there can exist at most one continuous character $\delta_A:K^\times \to A^\times$ such that $D_A$ contains an $\cR_{K,E}$-saturated $(\varphi,\Gamma_K)$-submodule $D_{A,0}$ over $\cR_{K,A}$
and $(\delta_A \mod \frm_A):K^\times \to E^\times$ equals $\delta$.
Moreover, if such $\delta_A$ exists, $D_{A,0}$ is also the unique $\cR_{K,E}$-saturated $(\varphi,\Gamma_K)$-submodule of $D_A$ isomorphic to $\cR_{K,A}(\delta_A)$.
\item 
Suppose $D$ is trianguline with a triangulation $\cF^\bullet= (0=\cF^0 \subsetneq \cF^1\subsetneq \dots\subsetneq \cF^n = D)$ with parameter $(\delta_1,\dots,\delta_n)$ such that for all $1\leq i<j\leq n$, $\delta_i/\delta_j \notin \set{x^{\bf k}}{{\bf k} \in \bbN^{\Sigma_K}}$.
Then, for any trianguline deformation $D_A$ of $D$, $D_A$ can admit at most one triangulation $\Fil^\bullet(D_A)$ such that $(D_A, \Fil^\bullet(D_A))$ is an $\cF^\bullet$-trianguline deformation of $D$.
\end{enumerate}
\end{lem}
\begin{proof}
\begin{enumerate}[(i)]
\item 
Given $\delta:K^\times \to E^\times$, let $\delta_A^{(1)},\delta_A^{(2)}: K^\times \to A^\times$ be two lifts of $\delta$, and suppose $D_{A,0}^{(i)}$ is an $\cR_{K,E}$-saturated $(\varphi,\Gamma_K)$-submodule over $\cR_{K,A}$ of $D_A$ such that $D_{A,0}^{(i)}\cong \cR_{K,A}(\delta_A^{(i)})$, for $i=1,2$.
We need to show that $D_{A,0}^{(1)}=D_{A,0}^{(2)}$ as submodules of $D_A$ and thus $\delta_A^{(1)}=\delta_A^{(2)}$.
It suffices to show that, for any $i,j\in \{1,2\}$, we have 
\begin{equation}
\label{eqn:vanishing}
\Hom_{\varphi,\Gamma_K}(\cR_{K,A}(\delta_A^{(j)}), D_A/D_{A,0}^{(i)}) = 0,
\end{equation}
from which we can conclude that $D_{A,0}^{(j)}\subset D_{A,0}^{(i)}$, and hence by symmetry, $D_{A,0}^{(1)}=D_{A,0}^{(2)}$.

Since $D_{A,0}^{(i)}$ is a free $A$-submodule of the free $A$-module $D_A$, the quotient $D_A/D_{A,0}^{(i)}$ is free over $A$ by \cite[Lemma 2.2.3(i)]{BeCh_families}.
Then, $D_{A,0}^{(i)}$ is a direct summand of $D_A$ as $A$-module, and hence $D_{A,0}^{(i)}\cap \frm_A D_A = \frm_A D_{A,0}^{(i)}$, by which 
$
D_{A,0}^{(i)}/\frm_A D_{A,0}^{(i)} \cong \cR_{K,E}(\delta)
$
is a saturated $(\varphi,\Gamma_K)$-submodule of $D_A/\frm_A D_A \cong D$ over $\cR_{K,E}$.
Since $\Hom_{\varphi,\Gamma_K}(R_{K,E}(\delta), D/D_0)$ vanishes, $D_{A,0}^{(i)}/\frm_A D_{A,0}^{(i)}$ is a saturated submodule of $D_0$ of rank 1.
Therefore, $D_{A,0}^{(i)}/\frm_A D_{A,0}^{(i)} = D_0$.
For $M_i:= D_A/D_{A,0}^{(i)}$, we prove \eqref{eqn:vanishing} by dévissage through the $\frm_A$-adic filtration:
$$
0 = \frm_A^N M_i \subset \dots \subset \frm_A^2 M_i \subset \frm_AM_i \subset M_i
$$
for some large $N>0$ such that $\frm_A^N=0$, whose graded pieces are
$$
\frm_A^l M_i / \frm_A^{l+1}M_i \cong M_i\otimes_A (\frm_A^l/\frm_A^{l+1}) \cong (D/D_0)^{\dim_E(\frm_A^l/\frm_A^{l+1})}
$$
by the projectivity of $M_i$ over $A$.
By left exactness of $\Hom_{\varphi,\Gamma_K}(\cR_{K,A}(\delta_A^{(j)}), - )$, it is enough to show  $\Hom_{\varphi,\Gamma_K}(\cR_{K,A}(\delta_A^{(j)}), M_i\otimes_A (\frm_A^l/\frm_A^{l+1})) = 0$ for all $0\leq l< N$.
Any nonzero map in $\Hom_{\varphi,\Gamma_K}(\cR_{K,A}(\delta_A^{(j)}), M_i\otimes_A (\frm_A^l/\frm_A^{l+1}))$ factors through $\cR_{K,A}(\delta_A^{(j)})\otimes_A E = \cR_{K,E}(\delta)$, which contradicts the fact that $\Hom_{\varphi,\Gamma_K}(\cR_{K,E}(\delta),D/D_0)=0$.
Thus, \eqref{eqn:vanishing} is proved.
\item 
We claim that $\Fil^i(D_A)/\Fil^{i-1}(D_A)$ is the unique $(\varphi,\Gamma_K)$-submodule of $D_A/\Fil^{i-1}(D_A)$ that lifts the submodule $\cF^{i}/\cF^{i-1}$ of rank $1$ in $D/\cF^{i-1}$, for all $1\leq i\leq n$, which is enough for us to conclude the uniqueness of the $\cF^\bullet$-triangulation on $D_A$.
Under the assumption that for all $1\leq i<j\leq n$, $\delta_j/\delta_i \notin \set{x^{-\bf k}}{{\bf k} \in \bbN^{\Sigma_K}}$, we deduce from Lemma \ref{lem:rank_one_coh}(ii) that
\begin{align*}
\Hom_{\varphi,\Gamma_K}(\cR_{K,E}(\delta_i), (D/\cF^{i-1})/(\cF^{i}/\cF^{i-1})) 
&=\Hom_{\varphi,\Gamma_K}(\cR_{K,E}(\delta_i), D/\cF^{i})
\\&=
H^0((D/\cF^{i})(\delta_i^{-1})) = 0.
\end{align*}
Therefore, we can apply (i) to the deformation $D_{A}/\Fil^{i-1}(D_A)$ and the character $\delta_i$ to conclude the uniqueness of $\Fil^i(D_A)/\Fil^{i-1}(D_A)$ in $D_A/\Fil^{i-1}(D_A)$.
\qedhere
\end{enumerate} 
\end{proof}

\begin{rem}
\label{rem:triang_deform}
For a $(\varphi,\Gamma_K)$-module $D$ of rank $n$ over $\cR_{K,E}$,
if the local Artinian $E$-algebra $A$ is $E[\vep]/\vep^2$, we use $\widetilde{D}$ instead of $D_{E[\vep]/\vep^2}$ to denote the deformations of $D$ to $\cR_{K,E[\vep]/\vep^2}$.
We can identify {\it elements} in $\Ext^1(D,D)$ with deformations of $D$ over $R_{K,E[\vep]/\vep^2}$ as follows. 
Let 
$$
0 \to D \xrightarrow{\iota} \widetilde{D} \xrightarrow{\pi} D \to 0 
$$
be a self-extension of $D$ as a $(\varphi,\Gamma_K)$-module over $\cR_{K,E}$.
We define an action of $\vep$ on $\widetilde{D}$ by 
$$
\vep_{\widetilde{D}}: \widetilde{D}\xrightarrow{\pi} D \xrightarrow{\id} D \xrightarrow{\iota} \widetilde{D}.
$$
Since $\iota$ and $\pi$ are $(\varphi,\Gamma_K)$-equivariant, $\vep_{\widetilde{D}}$ commutes with $\varphi$ and $\Gamma_K$.
Also, $\vep_{\widetilde{D}}^2=\iota \circ \pi \circ \iota\circ \pi = \iota\circ 0 \circ \pi = 0$, which shows that $\widetilde{D}$ is a $(\varphi,\Gamma_K)$-module of rank $n$ over $\cR_{K,E[\vep]/\vep^2}$ with $\widetilde{D}/\vep \widetilde{D}\cong D$.
Thus, $\widetilde{D}$ is a deformation of $D$ to $R_{K,E[\vep]/\vep^2}$.
Conversely, suppose $\widetilde{D}$ is a $(\varphi,\Gamma_K)$-module of rank $n$ over $\cR_{K,E[\vep]/\vep^2}$ such that $\widetilde{D}/\vep\widetilde{D} \cong D$.
Then,
$
\vep: \widetilde{D}/\vep\widetilde{D} \xrightarrow{} \vep \widetilde{D}
$
is an isomorphism of $(\varphi,\Gamma_K)$-modules over $\cR_{K,E}$ because, as $\cR_{K,E[\vep]/\vep^2}$-module, $\widetilde{D}\cong (\cR_{K,E[\vep]/\vep^2})^n = (\cR_{K,E})^n\otimes_E E[\vep]/\vep^2$.
We therefore get a self-extension 
$$
0 \to \vep \widetilde{D}\cong D \hookrightarrow \widetilde{D} \twoheadrightarrow \widetilde{D}/\vep \widetilde{D}\cong D \to 0
$$
of $D$, which is an element in $\Ext^1(D,D)$.
\end{rem}
\begin{defn}
\label{defn:subspace_trianguline}
For any given triangulation $\cF^\bullet$ on a $(\varphi,\Gamma_K)$-module $D$ of rank $n$
over $\cR_{K,E}$, we define the {\bf subspace of $\cF^\bullet$-trianguline deformations} $\Ext^1_{\cF^\bullet}(D,D)$ to be the subspace of
$ 
\Ext^1(D,D)
$
consisting of deformations $\widetilde{D}$ of $D$ to $\cR_{K,E[\vep]/\vep^2}$ that admit an $\cF^\bullet$-triangulation $\Fil^\bullet(\widetilde{D})$ under the identification explained in Remark \ref{rem:triang_deform}.
\end{defn}

\begin{rem}
In Definition \ref{defn:subspace_trianguline}, for a deformation $\widetilde{D}$ of $D$ that is trianguline, there may be several  triangulations $\Fil^\bullet(\widetilde{D})$ on $\widetilde{D}$ that lift $\cF^\bullet$.
If the triangulation parameters $(\delta_1,\dots,\delta_n)$ of $\cF^\bullet$ satisfy the condition that $\delta_i/\delta_j \notin \set{x^{\bf k}}{{\bf k} \in \bbN^{\Sigma_K}}$ for $1\leq i<j\leq n$, then Lemma \ref{lem:unique_deform}(ii) applies to guarantee the uniqueness of the $\cF^\bullet$-triangulation $\Fil^\bullet(\widetilde{D})$ on $\widetilde{D}$.
\end{rem}

To verify that $\Ext^1_{\cF^\bullet}(D,D)$ is an {\it $E$-subspace} of $\Ext^1(D,D)$, we recall the definition of Baer sum.
Let $B,C$ be two $(\varphi,\Gamma_K)$-modules over $\cR_{K,E}$.
For any two extensions in $\Ext^1(B,C) $, 
$$
0 \to C \to D_1 \to B \to 0, \quad 0 \to C \to D_2 \to B \to 0,
$$
their Baer sum
``$D_1+D_2$'' is calculated by taking the direct sum
$$
0 \to C\oplus C \to D_1\oplus D_2 \to B\oplus B \to 0,
$$
pushing out via the sum $C\oplus C \xrightarrow{\Sigma} C$, and then pulling back via the diagonal $B \xrightarrow{\Delta} B\oplus B$:
\begin{equation}
\label{eqn:Baer}
\begin{tikzcd}
0 \arrow[r] & C\oplus C \arrow[r] \arrow[d, "\Sigma"] & D_1\oplus D_2 \arrow[r] \arrow[d] & B \oplus B \arrow[r] \arrow[d, Rightarrow, no head]& 0 \\
0 \arrow[r] & C \arrow[r] & C \sqcup_{C\oplus C} (D_1\oplus D_2) \arrow[r] & B \oplus B \arrow[r] & 0 \\
0 \arrow[r] & C \arrow[r] \arrow[u, Rightarrow, no head] & D_1+D_2 \arrow[r] \arrow[u] & B \arrow[r] \arrow[u, "\Delta"]& 0.
\end{tikzcd}
\end{equation}

Moreover, if there are further extensions of $(\varphi,\Gamma_K)$-modules
$$
0 \to C' \to C \to C'' \to 0, \quad 0 \to D_i' \to D_i \to D_i'' \to 0, \quad  0 \to B' \to B \to B'' \to 0,
$$
for $i=1,2$ that induce $9$-term commutative diagrams with exact rows and columns
\begin{equation}
\label{eqn:9-term}
\begin{tikzcd}
            & 0 \arrow[d]             & 0 \arrow[d]               & 0 \arrow[d]             &   \\
0 \arrow[r] & C' \arrow[r] \arrow[d]  & D_i' \arrow[r] \arrow[d]  & B' \arrow[r] \arrow[d]  & 0 \\
0 \arrow[r] & C \arrow[r] \arrow[d]   & D_i \arrow[r] \arrow[d]   & B \arrow[r] \arrow[d]   & 0 \\
0 \arrow[r] & C'' \arrow[r] \arrow[d] & D_i'' \arrow[r] \arrow[d] & B'' \arrow[r] \arrow[d] & 0 \\
            & 0                       & 0                         & 0                       &  
\end{tikzcd}
\end{equation}
for $i=1,2,$ then we naturally have a complex
\begin{equation}
\label{eqn:Baer&subquot}
0 \to (D_1'+D_2') \to (D_1+D_2) \to (D_1''+D_2'') \to 0,
\end{equation}
which is exact, since after taking pushout and pullback, we have a commutative diagram
$$
\begin{tikzcd}
            & 0 \arrow[d]             & 0 \arrow[d]                     & 0 \arrow[d]             &   \\
0 \arrow[r] & C' \arrow[r] \arrow[d]  & D_1'+D_2' \arrow[r] \arrow[d]   & B' \arrow[r] \arrow[d]  & 0 \\
0 \arrow[r] & C \arrow[r] \arrow[d]   & D_1+D_2 \arrow[r] \arrow[d]     & B \arrow[r] \arrow[d]   & 0 \\
0 \arrow[r] & C'' \arrow[r] \arrow[d] & D_1''+D_2'' \arrow[r] \arrow[d] & B'' \arrow[r] \arrow[d] & 0 \\
            & 0                       & 0                               & 0                       &  
\end{tikzcd}
$$
with exact rows and two exact columns. 
By the 9-lemma, the middle column \eqref{eqn:Baer&subquot} is exact.

\begin{prop}
\label{prop:triang_subspace}
For any triangulation $\cF^\bullet$ on a $(\varphi,\Gamma_K)$-module $D$ of rank $n$
over $\cR_{K,E}$, the {subspace of $\cF^\bullet$-trianguline deformations}
$
\Ext^1_{\cF^\bullet}(D,D)
$ 
is an $E$-linear subspace of $\Ext^1(D,D)$.
\end{prop}
\begin{proof}
Let $\widetilde{D}_1$ and $\widetilde{D}_2$ be two deformations of $D$ to $\cR_{K,E[\vep]/\vep^2}$ lying in the subset $\Ext^1_{\cF^\bullet}(D,D)$, respectively equipped with $\cF^\bullet$-triangulations $\Fil^\bullet(\widetilde{D}_1)$ and $\Fil^\bullet(\widetilde{D}_2)$.
For each $1\leq j\leq n$ and $i=1,2$, we naturally have a commutative diagram similar to  \eqref{eqn:9-term} by $$
\begin{tikzcd}
            & 0 \arrow[d]                 & 0 \arrow[d]                                                 & 0 \arrow[d]                 &   \\
0 \arrow[r] & \cF^j \arrow[r] \arrow[d]   & \Fil^j(\widetilde{D}_i) \arrow[r] \arrow[d]                 & \cF^j \arrow[r] \arrow[d]   & 0 \\
0 \arrow[r] & D \arrow[r] \arrow[d]       & \widetilde{D}_i \arrow[r] \arrow[d]                         & D \arrow[r] \arrow[d]       & 0 \\
0 \arrow[r] & D/\cF^j \arrow[r] \arrow[d] & \widetilde{D}_i/\Fil^j(\widetilde{D}_i) \arrow[r] \arrow[d] & D/\cF^j \arrow[r] \arrow[d] & 0. \\
            & 0                           & 0                                                           & 0                           &  
\end{tikzcd}
$$
Hence, by the exactness of \eqref{eqn:Baer&subquot}, we have exact sequences
$$
0
\to 
\Fil^j(\widetilde{D}_1)+\Fil^j(\widetilde{D}_2)
\to
\widetilde{D}_1+\widetilde{D}_2
\to 
(\widetilde{D}_1/\Fil^j(\widetilde{D}_1))+(\widetilde{D}_2/\Fil^j(\widetilde{D}_2))
\to 
0
$$
for $1\leq j\leq n$, which shows that $\Fil^\bullet(\widetilde{D}_1)+\Fil^\bullet(\widetilde{D}_2)$ is an $\cF^\bullet$-triangulation on $\widetilde{D}_1+\widetilde{D}_2$.
Note that by applying Baer sum to the following commutative diagrams for $i=1,2,$
$$
\begin{tikzcd}
            & 0 \arrow[d]                            & 0 \arrow[d]                                                & 0 \arrow[d]                            &   \\
0 \arrow[r] & \cF^{j-1} \arrow[r] \arrow[d]          & \Fil^{j-1}(\widetilde{D}_i) \arrow[r] \arrow[d]            & \cF^{j-1} \arrow[r] \arrow[d]          & 0 \\
0 \arrow[r] & \cF^j \arrow[r] \arrow[d]              & \Fil^{j}(\widetilde{D}_i) \arrow[r] \arrow[d]              & \cF^j \arrow[r] \arrow[d]              & 0 \\
0 \arrow[r] & \gr^j(\cF^\bullet) \arrow[r] \arrow[d] & \gr^j(\Fil^{\bullet}(\widetilde{D}_i)) \arrow[r] \arrow[d] & \gr^j(\cF^\bullet) \arrow[r] \arrow[d] & 0 \\
            & 0                                      & 0                                                          & 0                                      &  
\end{tikzcd}
$$
we obtain that
$
\gr^j(\Fil^\bullet(\widetilde{D}_1)+\Fil^\bullet(\widetilde{D}_2)) = \gr^j(\Fil^\bullet(\widetilde{D}_1))+\gr^j(\Fil^\bullet(\widetilde{D}_2))
$ for $1\leq j\leq n$.
To show that $\Ext^1_{\cF}(D,D)$ is stable under $E$, let $0\to D\to \widetilde{D}\to D\to 0$ be a deformation of $D$ equipped with an $\cF^\bullet$-triangulation $\Fil^\bullet(\widetilde{D})$.
For $\alpha\in E$, the scalar multiple $\alpha[\widetilde{D}]$ is the extension obtained by pushing $0\to D\to \widetilde{D}\to D\to 0$ out along the map $D\xrightarrow{\alpha}D$.
Therefore, $\alpha[\Fil^\bullet(\widetilde{D})]$ is an $\cF^\bullet$-triangulation on $\alpha[\widetilde{D}]$, and $\Ext^1_{\cF^\bullet}(D,D)$ is an $E$-linear subspace of $\Ext^1(D,D)$.
\end{proof}

For a $(\varphi,\Gamma_K)$-module $D$ of rank $1$, its deformation to $\cR_{K,E[\vep]/\vep^2}$ has the following description.
\begin{prop}
\label{prop:rank_one_deform}
For any continuous character $\delta: K^\times \to E^\times$, we have a canonical $E$-linear isomorphism
$$
\Hom(K^\times,E)\xrightarrow{\sim }
\Ext^1(\cR_{K,E}(\delta),\cR_{K,E}(\delta)), \quad \psi \mapsto \cR(\delta(1+\psi\vep)),
$$
where $\Hom(K^\times,E)$ denotes the vector space of continuous homomorphisms from $K^\times$ to $E$.
\end{prop}
\begin{proof}
Twisting by $\delta^{-1}$ gives an isomorphism $\Ext^1(\cR_{K,E}(\delta),\cR_{K,E}(\delta))\xrightarrow{\sim}\Ext^1(\cR_{K,E},\cR_{K,E})$.
Thus, we may assume that $\delta=1$.
Given a continuous character $\psi: K^\times \to E$, we have
$$
(1+\psi(a)\vep)(1+\psi(b)\vep) = 1 + (\psi(a)+\psi(b))\vep = 1+\psi(ab)\vep
$$
for any $a,b\in K^\times$.
Hence, $1+\psi\vep: K^\times \to (E[\vep]/\vep^2)^\times$ is a continuous character that lifts the trivial character $1$, and $\cR_{K,E[\vep]/\vep^2}(1+\psi\vep)$ deforms $\cR_{K,E}$.
Conversely, given any deformation $\widetilde{D}$ of $\cR_{K,E}$ to $\cR_{K,E[\vep]/\vep^2}$, by Theorem \ref{thm:rank-one_classification} there is a unique character $\widetilde{\delta}: K^\times \to (E[\vep]/\vep^2)^\times$ such that $(\widetilde{\delta} \mod \vep) = 1$ and $\widetilde{D}\cong \cR_{K,E[\vep]/\vep^2}(\widetilde{\delta})$.
Hence, $\widetilde{\delta}=1+\psi\vep$ for a continuous map $\psi:K^\times\to E$, which is a homomorphism because $\widetilde{\delta}$ is.
The map $\Hom(K^\times,E)\to \Ext^1(\cR_{K,E},\cR_{K,E})$ is an $E$-linear map because it is the comparison isomorphism \cite[Proposition 2.3.7]{KPX_cohomology}
$$
H^1_{\rm cont}(G_K, V) \xrightarrow[\sim]{D_{\rig}} H^1(D_{\rig}(V)) 
$$
between continuous Galois cohomology and $(\varphi,\Gamma_K)$-cohomology for the trivial 1-dimensional representation $V:=E$ of $G_K$ and the associated $(\varphi,\Gamma_K)$-module $D_{\rig}(V)=\cR_{K,E}$. 
Indeed, by local class field theory, we have isomorphisms
$$
H^1_{\rm cont}(G_K, E)=\Hom_{\rm cont}(\widehat{K^\times},E) = \Hom(K^\times,E),
$$
where a continuous character $\psi\in \Hom(K^\times,E)$ corresponds to the deformation $E(1+\psi\vep)\in \Ext^1_{G_K}(E,E) \cong H^1_{\rm cont}(G_K,E)$ of the trivial representation $E$.
Under $D_{\rig}$, this gives our map 
\begin{equation*}
\Hom_{\rm cont}(K^\times,E)\to \Ext^1(\cR_{K,E},\cR_{K,E}), \quad \psi \mapsto \cR_{K,E}(1+\psi\vep) = D_{\rig}(E(1+\psi\vep)).
\qedhere
\end{equation*}
\end{proof}

\subsection{Deformations of locally analytic principal series}
\label{subsection:translation_foundation}

We slightly generalize results on translation of locally analytic principal series in \cite{JLS_translation} to Artinian coefficients. 
Let $A\in \cC_E$ be a local Artinian $E$-algebra of maximal ideal $\frm_A$ and residue field $A/\frm_A\cong E$.

\subsubsection{}
Let $\bfG \supset \bfP \supset \bfB \supset \bfT$ be a split connected reductive group over $K$ together with a parabolic subgroup, a Borel subgroup and a $K$-split maximal split torus.
Write $G:=\bfG(K), P:=\bfP(K), B:=\bfB(K), T:=\bfT(K)$ for their groups of $K$-points.
Let $\cG \supset \cP \supset \cB \supset \cT$ be the Weil restrictions of scalars from $K$ to $\bbQ_p$, i.e., $\cG := \Res^K_{\bbQ_p}(\bfG), \cB:=\Res^K_{\bbQ_p}(\bfB)$, etc.
Then, $\cG(\bbQ_p)=G, \cP(\bbQ_p)=P, \cB(\bbQ_p)=B, \cT(\bbQ_p)=T$ are locally $K$-analytic groups, whose $\bbQ_p$-Lie algebras are naturally $K$-vector spaces that we denote by $\frg, \frp,\frb,\frt$ respectively.
Denote their base changes from $\bbQ_p$ to $E$ by $\frg_E, \frp_E, \frb_E, \frt_E$, and let $\frt_E^*$ be the {\it weight space}, which is the $E$-linear dual of $\frt_E$.
Via the isomorphism $K\otimes_{\bbQ_p}E \cong \prod_{\sigma\in \Sigma_K} E$, we have $\frg_E=\prod_{\sigma} \frg_E\otimes_{K,\sigma}E$, $\frb_E=\prod_{\sigma} \frb_E\otimes_{K,\sigma}E$, etc.
Then, $\frg_E$ is a split reductive Lie algebra over $E$ with Cartan subalgebra $\frt_E$ because $\cG\times_{\bbQ_p} E \cong \prod_{\sigma\in\Sigma_K} \bfG \times_{K,\sigma} E$ is a split connected reductive group over $E$ of Lie algebra $\frg_E$ with split maximal torus $\cT \times_{\bbQ_p} E$.
Via the inclusion 
$$
\bfG(K) =\cG(\bbQ_p) \subset \cG(E) = (\cG\times_{\bbQ_p}E)(E),
$$
we can use the highest weight theory for the split $E$-group $\cG\times_{\bbQ_p}E$: if $L(\lambda)$ is the irreducible algebraic $\cG\times_{\bbQ_p}E$-representation of highest weight $\lambda = (\lambda_\sigma)_{\sigma\in \Sigma_K} \in \frt_E^* \cong \prod_{\sigma} \frt^* \otimes_{K,\sigma} E$, then its restriction to $G=\bfG(K)$ is still irreducible by Zariski density of $\cG(\bbQ_p)$ in $\cG(E)$, and its differential is the irreducible $U(\frg_E)$-module of highest weight $\lambda$, where $U(\frg_E)$ is the universal enveloping algebra of the Lie algebra $\frg_E$.

Let $D(G)$ be the locally $\bbQ_p$-analytic distribution algebra of $G$ over $E$, which is defined to be the strong dual of the locally convex topological $E$-vector space $C^{\la}(G,E)$ of locally $\bbQ_p$-analytic $E$-valued functions on $G$.
Similarly, we define $D(B),D(P)$ with $D(B)\subset D(P)\subset D(G)$.
Let $D(\frg,P)$ be the smallest subring of $D(G)$ containing both $U(\frg_E)$ and $D(P)$, cf. \cite[\S2.1.4]{JLS_translation}.
We have a canonical embedding $U(\frg_E) \hookrightarrow D(G)$ such that the center $\frz_E:=Z(U(\frg_E))$ embeds into the center $Z(D(G))$.
We have the {\it translation functor} on $\frz$-finite $D(G)$-modules
$$
T_\lambda^\mu : D(G)\text{-mod}_{\frz\text{-fin}}\to D(G)\text{-mod}_{\frz\text{-fin}}, \quad M \mapsto \pr_{|\mu|}(L(\overline{\nu})\otimes_E \pr_{|\lambda|}(M)),
$$
whenever $\lambda,\mu\in\frt_E^*$ are such that $\mu-\lambda$ lifts to an $E^\times$-valued algebraic character of $\bfT$.
Here, $\pr_{|\mu|}$ is the projection to the generalized eigenspace of $\frz_E$ for its eigencharacter $\chi_\mu$ determined by the action of $\frz_E$ on the Verma module $U(\frg_E)\otimes_{U(\frb_E)} E_\mu$, and $\overline{\nu}$ is the unique dominant weight in the Weyl orbit of the integral weight $\nu:=\mu-\lambda\in \frt_E^*$.
For any $\lambda \in \frt^*_E$, let $\frm_{\chi_\lambda}:=\ker(\chi_\lambda) \subset Z(U(\frg_E))$ be the maximal ideal given by the kernel of $\chi_\lambda$.
For a locally $\bbQ_p$-analytic representation $V$ of $G$ on topological $E$-vector space of compact type, let $V'_b$ denote its strong dual, which is naturally a separately continuous $D(G)$-module.
For a locally $\bbQ_p$-analytic representation $V$ of $P$ over $E$, let $\Ind_P^G(V)^\la$ be the locally $\bbQ_p$-analytic parabolic induction
$$
\operatorname{Ind}_P^G(V)^\la=\set{f: G \rightarrow V \text { locally } \bbQ_p \text {-analytic }}{ f(hg)=h.f(g), \forall g \in G, \forall h \in P}.
$$

\begin{rem}
This is slightly more general than the setup of \cite{JLS_translation}, where the authors take a split connected reductive group $\bfH$ over $F$ and consider locally $F$-analytic representations of $\bfH(F)$.
We are taking $F=\bbQ_p$ and $\bfH = \Res^{K}_{\bbQ_p}(\bfG)$ so that $\bfH(F) = G$, which may not be split over $\bbQ_p$.
But since we always consider locally $\bbQ_p$-analytic representations on $E$-vector spaces, and $\bfH \times_{\bbQ_p}E \cong \prod_{\sigma\in \Sigma_K} \bfG\times_{K,\sigma}E$ is a split connected reductive group over $E$, the results in \cite{JLS_translation} can be proved in this case using the same arguments.
Cf. \cite[\S2]{Breuil_socle1} and \cite[\S1]{Ding_wallcross}.
\end{rem}

\begin{ex}
\label{ex:4.1.4/7}
Here we generalize \cite[Example 4.1.4, Example 4.1.7]{JLS_translation}.
Let $\tau_A: P \to A^\times$ be a locally $\bbQ_p$-analytic\footnote{Equivalently, a continuous character $\tau_A:P\to A^\times$ by \cite[Proposition 8.3]{Buzzard_eigenvar}.} character.
Let $d\tau_A: \frp \to A$ be the $\bbQ_p$-linear derivative of $\tau_A$, which induces an $E$-linear map $d\tau_A:\frp_E\to A$ that defines a $U(\frp_E)$-module structure on $A$.
More generally, the locally analytic function $\tau_A \in C^\la(P,A)\cong C^\la(P,E)\otimes_E A$ induces an $E$-linear map $\chi_{\tau_A}:D(P)=C^\la(P,E)'_b \to A$.
One can verify that $\chi_{\tau_A}$ is an $E$-algebra map as $\tau_A$ is a group homomorphism.
Let $M$ be a finitely generated $A$-module.
\begin{enumerate}[(i)]
\item 
We write $M_{\tau_A}$ for the $P$-representation given by the composite $P \xrightarrow{\tau_A} A^\times \to \Aut_E(M)$, which is also a $D(P)$-module by the composite $D(P) \xrightarrow{\chi_{\tau_A}} A \to \End_E(M)$.
Let $(M_{\tau_A})^\vee := \Hom_E(M_{\tau_A},E)$ be the $E$-linear dual on which $P$ acts by the contragredient action. 
\item 
We write $M_{d\tau_A}$ for the $U(\frp_E)$-module structure on $M$ given by the composite $U(\frp_E)\hookrightarrow D(P) \xrightarrow{\chi_{\tau_A}} A \to \End_E(M)$.
\item 
The induction $U(\frg_E)\otimes_{U(\frp_E)} M_{d\tau_A}$ is a $U(\frg_E)$-module equipped with a locally $\bbQ_p$-analytic $P$-action given by the formula $h(u\otimes m) = {\rm Ad}(h)(u) \otimes \tau_A(h)m$ for $h\in P, u\in U(\frg_E)$ and $m\in M_{d\tau_A}$, where ${\rm Ad}$ denotes the adjoint action of $G$ on $U(\frg_E)$.
The actions of $U(\frg_E)$ and $P$ are compatible in the sense that $U(\frg_E)\otimes_{U(\frp_E)} M_{d\tau_A}$ is an object of the category $\cO^{P,\infty}$ \cite[\S4.1.2]{JLS_translation}, which is naturally a $D(\frg,P)$-module by \cite[Lemma 7.4.2]{Agrawal-Strauch}.
\item 
The same argument given in \cite[Example 4.1.7]{JLS_translation} shows that the inclusion $U(\frg_E)\hookrightarrow D(\frg,P)$ induces an $A$-linear isomorphism of $D(\frg,P)$-modules
$$
U(\frg_E)\otimes_{U(\frp_E)} M_{d\tau_A} \xrightarrow{\sim} D(\frg,P)\otimes_{D(P)} M_{\tau_A},
$$
and hence an $A$-linear isomorphism of $D(G)$-modules
$$
\check{\cF_P^G}(U(\frg_E)\otimes_{U(\frp_E)} M_{d\tau_A}):= D(G)\otimes_{D(\frg,P)}(U(\frg_E)\otimes_{U(\frp_E)} M_{d\tau_A}) \xrightarrow{\sim} D(G)\otimes_{D(P)} M_{\tau_A}.
$$
\end{enumerate}
\end{ex}

\begin{lem}
\label{lem:ST_dual}
For a finite $A$-module $M$ and locally $\bbQ_p$-analytic character $\tau_A: P\to A^\times$, we have the following isomorphism of $D(G)$-modules
\begin{equation}
\label{eqn:ST_dual}
D(G)\otimes_{D(P)} (M_{\tau_A})^\vee \xrightarrow{\sim} (\Ind_P^G(M_{\tau_A})^{\la})'_b, \quad (\mu \otimes \ell) \mapsto (f \mapsto \mu(\ell\circ f)),
\end{equation}
for any distribution $\mu \in D(G)$, linear functional $\ell \in (M_{\tau_A})^\vee$, and function $f\in \Ind_P^G(M_{\tau_A})^\la$.
\end{lem}
\begin{proof}
Since $A$ is local Artinian and $M$ is finitely generated over $A$, $M$ is of finite length over $A$ with successive quotients given by $A/\frm_A\cong E$.
Thus, the submodule $M[\frm_A]$ is nonzero, and any nonzero $m\in M[\frm_A]$ generates a $1$-dimensional $E$-vector space $Am=Em$ on which $P$ acts by the mod-$\frm_A$ reduction $\tau:=(\tau_A \mod \frm_A)$ of $\tau_A$.
Consider
$$
0 \to E_{\tau} \to {M}_{\tau_A} \to (M/Am)_{\tau_A} \to 0 ,
$$
which is an exact sequence of $D(P)$-modules.
\begin{enumerate}[(i)]
\item 
Applying the locally analytic parabolic induction and continuous dual, which are exact, we get an exact sequence
$$
0 \to (\Ind_P^G((M/Am)_{\tau_A})^\la)'_b \to (\Ind_{P}^G(M_{\tau_A})^\la)'_b \to (\Ind_{P}^G(E_{\tau})^\la)'_b \to 0
$$
of $D(G)$-modules.
\item 
Applying the exact $E$-linear dual and $D(G)\otimes_{D(P)}(-)$, we get a sequence
$$
0 \to D(G)\otimes_{D(P)} ((M/Am)_{\tau_A})^\vee \to D(G)\otimes_{D(P)} (M_{\tau_A})^\vee \to D(G)\otimes_{D(P)} (E_{\tau})^\vee \to 0
$$
which is exact, since ${\rm Tor}_1^{D(P)}(D(G), ((M/Am)_{\tau_A})^\vee) =0$ by \cite[Lemma 6.3(i)]{Schneider-Teitelbaum_duality}, where the hypothesis (FIN) in {\it loc. cit.} is satisfied by $((M/Am)_{\tau_A})^\vee$, as $(E_{\tau})^\vee\cong E_{\tau^{-1}}$ satisfies (FIN) by \cite[Lemma 6.5]{Schneider-Teitelbaum_duality}, $((M/Am)_{\tau_A})^\vee$ is a successive extension of $(E_{\tau})^\vee$ as $D(P)$-module, and (FIN) is stable under extensions by \cite[Horseshoe Lemma 2.2.8]{Weibel_HA}.
\end{enumerate}
The natural map \eqref{eqn:ST_dual} defines a commutative diagram with exact rows
$$
\begin{tikzcd}[scale cd=0.9]
0 \arrow[r] & D(G)\otimes_{D(P)} ((M/Am)_{\tau_A})^\vee \arrow[d] \arrow[r] & D(G)\otimes_{D(P)} (M_{\tau_A})^\vee \arrow[r] \arrow[d, "\eqref{eqn:ST_dual}"] & D(G)\otimes_{D(P)} (E_{\tau})^\vee \arrow[r] \arrow[d] & 0 \\
0 \arrow[r] & (\Ind_P^G((M/Am)_{\tau_A})^\la)'_b \arrow[r]                        &  (\Ind_{P}^G(M_{\tau_A})^\la)'_b \arrow[r]                                            & (\Ind_{P}^G(E_{\tau})^\la)'_b  \arrow[r]                     & 0.
\end{tikzcd}
$$

The quotient $M/Am$ has shorter length than $M$ as an $A$-module, hence we have by induction on $\dim_E(M)$ that \eqref{eqn:ST_dual} holds for $(M/Am)_{\tau_A}$.
The validity of \eqref{eqn:ST_dual} when $\dim_E(M)=1$ follows from \cite[Lemma 6.1(iv)]{Schneider-Teitelbaum_duality}.
By the 5-lemma, \eqref{eqn:ST_dual} holds for $M_{\tau_A}$.
\end{proof}

\begin{rem}
\label{rem:dual_E[vep]}
Suppose the local Artinian $E$-algebra $A$ is {\it Frobenius} \cite[Definition 4.2.5]{Weibel_HA}, i.e., there is an $A$-module isomorphism $A\xrightarrow{\sim} A^\vee = \Hom_E(A,E)$.
For a finite free $A$-module $M$, we have $(M_{\tau_A})^\vee \cong M_{\tau_{A}^{-1}}$ as $P$-modules.
For example, $E[\vep]/\vep^2$ is a Frobenius $E$-algebra since 
$$
E[\vep]/\vep^2 \times E[\vep]/\vep^2 \to E,\quad (a+b\vep, c+d\vep) \mapsto ad+bc,
$$
is an $E[\vep]/\vep^2$-equivariant non-degenerate bilinear pairing, showing that $E[\vep]/\vep^2\cong (E[\vep]/\vep^2)^\vee$ as $E[\vep]/\vep^2$-modules.
We therefore have $(A_{\tau_A})^\vee\cong A_{\tau_A^{-1}}$ for $A=E[\vep]/\vep^2$.
\end{rem}

\begin{lem}
\label{lem:T_FGP}
Let $\lambda,\mu \in \frt_E^*$ be such that $\lambda-\mu$ is the weight of an algebraic character of $T$. 
Let $L_P$ be the Levi subgroup of $P$, and $\Rep_E^{\rm sm}(L_P)^{\text{s-adm}}$ be the category of strongly admissible smooth representations of $L_P$.
For any pair $(M,V)\in \cO^{P,\infty}\times \Rep_E^{\rm sm}(L_P)^{\text{s-adm}}$ such that $M$ is an $A$-module where $A$ acts by $D(\frg,P)$-linear endomorphisms, and $\check{\cF}_P^G(M,V):=D(G)\otimes_{D(\frg,P)}(M\otimes_E V'_b)$, there is a canonical $A$-linear isomorphism 
$$
T_\lambda^\mu (\check{\cF}_P^G(M,V)) \xrightarrow{\sim } \check{\cF}_P^G(T_\lambda^\mu(M),V)
$$
of $D(G)$-modules that is natural in $M$ and $V$.
\end{lem}
\begin{proof}
The isomorphism comes from {\cite[Theorem 4.1.12]{JLS_translation}}, and it is natural in $M\in \cO^{P,\infty}$ and $V\in \Rep_E^{\rm sm}(L_P)^{\text{s-adm}}$.
Since the $A$-actions on both sides are induced by that on $M$, which commutes with the action of $D(\frg,P)$, this isomorphism is $A$-linear by its naturality in $M$.
\end{proof}

\begin{lem}
\label{lem:translate_D(G)_induction}
Let $\widetilde{\lambda}_A, \widetilde{\mu}_A: T\to A^\times$ be locally $\bbQ_p$-analytic characters of $T$, let $\widetilde{\lambda}, \widetilde{\mu}: T\to E^\times$ be their reductions modulo $\frm_A$, with weights $\lambda:=d\widetilde{\lambda},\mu:=d\widetilde{\mu} \in \frt_E^*$.
Suppose that 
\begin{enumerate}[(i)]
\item 
$\lambda$ and $\mu$ are both anti-dominant, 
\item 
$\widetilde{\nu}:=\widetilde{\mu}_A/\widetilde{\lambda}_A$ is an $E^\times$-valued algebraic character of $T$, and
\item 
the condition (4.2.8)\footnote{It states that ``$\mu^{\natural}$ lies in the closure of the facet $\mathscr{F} \subset \mathscr{E}(\lambda)$ which contains $\lambda^{\natural}$."} of \cite{JLS_translation} is satisfied by $\lambda$ and $\mu$.
\end{enumerate}
Then, for any $w\in W_{[\lambda]}=W_{[\mu]}$ and finitely generated $A$-module $M$, 
there exists an $A$-linear isomorphism of $D(G)$-modules
$$
T_\lambda^\mu(D(G)\otimes_{D(B)} M_{w\cdot \widetilde{\lambda}_A}) \cong D(G)\otimes_{D(B)} M_{w\cdot \widetilde{\mu}_A}.
$$
Here, the dot action of $w\in W$ on the locally $\bbQ_p$-analytic character $\widetilde{\lambda}_A:T\to A^\times$ is defined via the same formula given in \cite[\S4.2.3]{JLS_translation}, so that the mod-$\frm_A$ reduction of $d(w\cdot \widetilde{\lambda}_A)$ equals $d(w\cdot \widetilde{\lambda}) = w\cdot \lambda$, which is the usual dot action of $w\in W$ on the weight $\lambda \in \frt_E^*$.
\end{lem}
\begin{proof}
By Example \ref{ex:4.1.4/7}(iv) and Lemma \ref{lem:T_FGP} (here we take $V=E$, the trivial representation of $L_B=T$), it suffices to prove that we have an $A$-linear isomorphism of $D(\frg,B)$-modules
$$
T_\lambda^\mu(U(\frg_E)\otimes_{U(\frb_E)}M_{d(w\cdot \widetilde{\lambda}_A)}) \cong U(\frg_E)\otimes_{U(\frb_E)} M_{d(w\cdot \widetilde{\mu}_A)}.
$$
Let $\overline{\nu}$ be the dominant weight in the Weyl orbit of $\nu:=\mu-\lambda$, and let $L(\overline{\nu})$ be the irreducible $G$-representation of highest weight $\overline{\nu}$.
Recall that $T_\lambda^\mu(N)=\pr_{|\mu|}(L(\overline{\nu})\otimes_E \pr_{|\lambda|}N)$.
For $N_A:=U(\frg_E)\otimes_{U(\frb_E)}M_{d(w\cdot \widetilde{\lambda}_A)}$, it is a successive extension of $N_A\otimes_A E = U(\frg_E)\otimes_{U(\frb_E)}M_{d(w\cdot \widetilde{\lambda})}$, which lies in $U(\frg_E)\text{-mod}_{|\lambda|}$.
Since $U(\frg_E)\text{-mod}_{|\lambda|}$ is stable under extensions, we get $\pr_{|\lambda|}(N_A)=N_A$.
Then, we filter $L(\overline{\nu})|_B$ by $B$-stable subspaces with 1-dimensional graded pieces on which $B$ acts through algebraic characters $\widetilde{\xi}_i:T\to E^\times$ as in \cite[Proposition 4.2.10]{JLS_translation}, where $\xi_i\in \frt_E^*$ runs through all the weights of $L(\overline{\nu})$ counting multiplicities.
Tensoring this filtration with $N_A$, we get a filtration on $L(\overline{\nu})\otimes_E N_A$ by $D(\frg,B)$-submodules with graded pieces $D(\frg,B)\otimes_{D(B)}(E_{\widetilde{\xi}_i} \otimes_E M_{w\cdot \widetilde{\lambda}_A}) = D(\frg,B)\otimes_{D(B)}(M_{\widetilde{\xi}_i(w\cdot \widetilde{\lambda}_A)}) \cong U(\frg_E)\otimes_{U(\frb_E)} M_{\xi_i + d(w\cdot \widetilde{\lambda}_A)}$, by Example \ref{ex:4.1.4/7}(iv).
Since $\pr_{|\mu|}$ is exact, $T_\lambda^\mu(N_A)$ has a filtration with graded pieces $\pr_{|\mu|}(U(\frg_E)\otimes_{U(\frb_E)} M_{\xi_i + d(w\cdot \widetilde{\lambda}_A)})$ appearing with multiplicity $\dim_E(L(\overline{\nu})_{\xi_i})$ for each $\xi_i$.
Again because the mod-$\frm_A$ reduction of $U(\frg_E)\otimes_{U(\frb_E)} M_{\xi_i + d(w\cdot \widetilde{\lambda}_A)}$, and hence $U(\frg_E)\otimes_{U(\frb_E)} M_{\xi_i + d(w\cdot \widetilde{\lambda}_A)}$ itself, lies in $U(\frg_E)\text{-mod}_{|\xi_i+w\cdot \lambda|}$, its projection to $U(\frg_E)\text{-mod}_{|\mu|}$ is nonzero if and only if $\xi_i+w\cdot \lambda$ lies in the orbit of $\mu$ under the dot-action of $W$.
Under the assumed condition (4.2.8) of \cite{JLS_translation}, by \cite[Corollary 4.2.9]{JLS_translation}, $\xi_i+w\cdot \lambda$ lies in the orbit of $\mu$ under the dot-action of $W$ if and only if $\xi_i = w(\nu)$, in which case we have multiplicity $\dim_E(L(\overline{\nu})_{w(\nu)})=\dim_E(L(\overline{\nu})_{\overline{\nu}}) =1$.
Therefore, 
\begin{align*}
T_\lambda^\mu(U(\frg_E)\otimes_{U(\frb_E)}M_{d(w\cdot \widetilde{\lambda}_A)}) 
&=
\pr_{|\mu|}(U(\frg_E)\otimes_{U(\frb_E)} M_{w({\nu}) + d(w\cdot \widetilde{\lambda}_A)})
\\&= 
U(\frg_E)\otimes_{U(\frb_E)} M_{w({\nu}) + d(w\cdot \widetilde{\lambda}_A)}
\\&\cong
D(\frg,B)\otimes_{D(B)} M_{w(\widetilde{\nu})(w\cdot\widetilde{\lambda}_A)} 
\\&= 
D(\frg,B)\otimes_{D(B)} M_{w\cdot \widetilde{\mu}_A}
\\&\cong
U(\frg_E)\otimes_{U(\frb_E)} M_{d(w\cdot \widetilde{\mu}_A)},
\end{align*}
is the desired isomorphism of $D(\frg,B)$-modules, which is $A$-linear as each step is $A$-linear.
\end{proof}

\subsection{Ding's correspondence and the proof of Theorem B}
\label{subsection:intertwining}

Let $K/\bbQ_p$ be a finite extension. 
Following \cite{Ding_crystabelline}, we take $\bfG=\GL_{n,K}$ in the setup of \S\ref{subsection:translation_foundation}, with $\bfB$ the upper triangular Borel subgroup containing $\bfT$ the diagonal torus in $\bfG$. 
Let $\bfB^-$ be the lower triangular Borel subgroup.
Then $G=\bfG(K) = \GL_n(K) \supset B = \bfB(K) \supset T = \bfT(K)$ and $B^{-}=\bfB^-(K)$.

Let $D\in \PGnc$ be non-critical crystabelline of regular Sen weights ${\bf h} = (h_{1,\sigma}>\dots > h_{n,\sigma})_{\sigma\in \Sigma_K}$.
Let ${\bf k}= (k_{i,\sigma})\in (\bbN^n)^{\Sigma_K}$ and let $p_{\bf k}(D)$ be obtained from $D$ by $p_{\bf k}=\prod_{i,\sigma}p_{i,\sigma}^{k_{i,\sigma}}$ such that the resulting weights ${\bf h'} = (h'_{1,\sigma}>\dots > h'_{n,\sigma})_{\sigma\in \Sigma_K}$ are still regular.\footnote{Note that $p_{\bf k}(D)$ is independent of the order in which the pullback operators $\prod_{i,\sigma}p_{i,\sigma}^{k_{i,\sigma}}$ are applied, because here $D$ is crystabelline of regular Sen weights ${\bf h}$, and ${\bf k}$ is chosen such that $p_{\bf k}(D)$ is of regular Sen weights ${\bf h'}$.
Namely, as $E$-filtered $(\varphi,G(L/K))$ modules, $D_{\cris}^L(p_{\bf k}(D))$ is the {\it unique} subobject of  $D_{\cris}^L(D)$ whose Hodge filtration has the same underlying flag as $\Fil^\bullet(D_{\dR}(D))$, but with jumps occurring at ${-\bf h'}$.
By Theorem \ref{thm:D_cris_WD_equivalence}, $p_{\bf k}(D)$ is the {\it unique} $(\varphi,\Gamma_K)$-submodule of $D$ with Sen weights ${\bf h'}$.} 
In this subsection, we recall the constructions of $\pi_{\min}(D)\subset \pi_{\fs}(D)$ from \cite{Ding_crystabelline} and prove Theorem B from \S\ref{thm:B}.

\subsubsection{}
We recall Ding's parameter map $\kappa_w$ on trianguline deformation spaces of $D$ to $E[\vep]/\vep^2$.
For any refinement $w\in S_n$ of $\phi$, let $\Fil_w^\bullet(D)$ be the corresponding non-critical triangulation of $D$ whose triangulation parameters are $\delta_{w,i} = \phi_{w(i)}\prod_{\sigma\in \Sigma_K} x_\sigma^{h_{i,\sigma}}$ for $1\leq i\leq n$.
Let $\Ext^1_w(D,D) \subset \Ext^1(D,D)$ be the subspace of $\Fil^\bullet_w(D)$-trianguline deformations to $\cR_{K,E[\vep]/\vep^2}$, cf. Remark \ref{rem:triang_deform} and Definition \ref{defn:subspace_trianguline}(ii).
Since $\phi$ is generic, we can apply Lemma \ref{lem:unique_deform}(ii) to conclude that on any deformation $\widetilde{D}\in \Ext^1_w(D,D)$, there is a unique triangulation $\Fil^\bullet(\widetilde{D})$ that lifts $\Fil^\bullet_w(D)$.
The graded pieces of $\Fil^\bullet(\widetilde{D})$ deform $\cR_{K,E}(\delta_{w,i})$, and hence are of the form $\cR_{K,E[\vep]/\vep^2}(\delta_{w,i}(1+\psi_i \vep))$ for continuous characters $\psi_i \in \Hom(K^\times,E)$ by Proposition \ref{prop:rank_one_deform}.
The tuple $\psi:=(\psi_1,\dots,\psi_n)$ defines a continuous character $\psi: T \to E$, and following \cite[(2.12)]{Ding_crystabelline} we define 
$$
\kappa_w : \Ext^1_w(D,D) \to \Hom(T,E) , \quad \widetilde{D}\mapsto \psi,
$$
which is an $E$-linear map by Proposition \ref{prop:rank_one_deform}.
By \cite[Proposition 2.10, Lemma 2.11]{Ding_crystabelline}, $\kappa_w$ is surjective, and the subspace $\ker(\kappa_w)$ of $\Ext^1(D,D)$ is independent of $w\in S_n$.
We denote this common subspace by $\Ext_0^1(D,D)$.
For any subspace $X\subset \Ext^1(D,D)$ containing $\Ext^1_0(D,D)$, we write  $\overline{X} := X/\Ext^1_0(D,D)$ for the quotient. 
Thus, we get an isomorphism 
$$
\kappa_w : \overline{\Ext}_w^1(D,D) \xrightarrow{\sim} \Hom(T,E) .
$$
The sum over all refinements generates the full deformation space: the diagram
$$
\begin{tikzcd}
\bigoplus_{w\in S_n} \Ext^1_w(D,D) \arrow[r, two heads] \arrow[d, two heads] & \Ext^1(D,D) \arrow[d, two heads]\\
\bigoplus_{w\in S_n} \overline{\Ext}^1_w(D,D) \arrow[r, two heads] & \overline{\Ext}^1(D,D)
\end{tikzcd}
$$
commutes with surjective rows and columns by \cite[Proposition 2.12]{Ding_crystabelline}.

\begin{lem}
\label{lem:fil_gr}
Let $A$ be a noetherian ring, $V$ a finite projective $A$-module, and $\Fil^\bullet V$ a separated and exhaustive descending filtration on $V$ by $A$-submodules with projective graded pieces $\gr_\Fil^k V:= \Fil^kV / \Fil^{k+1}V$.
Let $\fS$ be an $A$-algebra containing a non-zerodivisor $t$. 
For each jump $k$ in the filtration, choose an $A$-linear splitting $\Fil^kV = \Fil^{k+1}V\oplus \Gr_{k}$ such that $\Gr_k \cong \gr^k_{\Fil}V$.
Then,  $V=\bigoplus_{k}\Gr_k$ as $A$-modules, and 
$$
\sum_{k\in \bbZ} \Fil^kV\otimes_{A} t^{-k}\fS = \sum_{k\in \bbZ} \Gr_k \otimes_{A} t^{-k} \fS
$$
as $\fS$-submodules inside $V\otimes_{A}\fS[1/t]$.
\end{lem}
\begin{proof}
Since $A$ is noetherian, $V$ is finite projective over $A$, and $\gr^k_\Fil V$ are projective over $A$, it follows that only finitely many $\gr^k_\Fil V$ are nonzero. 
After choosing complementary subspaces $\Gr_k\subset \Fil^{k}V$ to $\Fil^{k+1}V$, we have $V=\bigoplus_{k\in \bbZ} \Gr_k$ and $\Fil^rV = \bigoplus_{k\geq r} \Gr_k$. 
Thus,
\begin{align*}
\sum_{k\in \bbZ} \Fil^kV\otimes_{A} t^{-k}\fS 
&=
\sum_{k\in \bbZ} (\bigoplus_{k'\geq k}\Gr_k')\otimes_{A} t^{-k}\fS 
=
\sum_{k\in \bbZ} \Gr_k\otimes_{A} t^{-k}\fS 
\end{align*}
as $\fS$-submodules of $V\otimes_A \fS[1/t]$,
because $t^{-k}\fS \subset t^{-k'}\fS$ for $k'> k$.
\end{proof}

\begin{lem}
\label{lem:pullback_parameter}
\begin{enumerate}[(i)]
\item
The operation $p_{\bf k}= \prod_{i,\sigma}p_{i,\sigma}^{k_{i,\sigma}}$ induces a well-defined map
$$
p_{\bf k}: \Ext^1(D,D) \xrightarrow{} \Ext^1(p_{\bf k}(D), p_{\bf k}(D)),\quad \widetilde{D}\mapsto p_{\bf k}(\widetilde{D}),
$$
independent of the order in which $p_{i,\sigma}$ are applied, where each $p_{i,\sigma}$ means the operation in Theorem \ref{thm:Wu_3.16}(i) that increases the smallest (modulo $\vep$) $i$ $\sigma$-Sen weights of $\widetilde{D}$ by $1$.
In fact, $p_{\bf k}$ equals the change of weights map $f_{{\bf h},{\bf h'}}: (\fX_n)^{\wedge}_{{\bf h}}\to (\fX_n)^{\wedge}_{{\bf h'}}$ of \cite[\S1.3]{Wu_translation} on $\Ext^1(D,D)$, in the sense that the following diagram is commutative
$$
\begin{tikzcd}
{\Ext^1(D,D)} \arrow[d, "\text{Remark \ref{rem:triang_deform}}", hook] \arrow[r, "p_{\bf k}"] & {\Ext^1(p_{\bf k}(D), p_{\bf k}(D))} \arrow[d, "\text{Remark \ref{rem:triang_deform}}", hook] \\
{(\fX_n)^{\wedge}_{{\bf h}}(E[\vep]/\vep^2)} \arrow[r, "{f_{{\bf h},{\bf h'}}}"]              & {(\fX_n)^{\wedge}_{{\bf h'}}(E[\vep]/\vep^2)}   .  
\end{tikzcd}
$$
\item 
The map $
p_{\bf k}: \Ext^1(D,D) \xrightarrow{\sim} \Ext^1(p_{\bf k}(D), p_{\bf k}(D))
$
is an $E$-linear isomorphism such that for every refinement $w$, it restricts to the pullback map on trianguline deformations
$$
p_{\bf k}: \Ext^1_w(D,D) \xrightarrow{\sim} \Ext^1_w(p_{\bf k}(D), p_{\bf k}(D)) .
$$
Moreover, the following diagram is commutative
$$
\begin{tikzcd}
{\Ext}^1_w(D,D) \arrow[r, "\kappa_w"] \arrow[d, "p_{\bf k}"] & \Hom(T,E) \arrow[d, "\id", Rightarrow, no head]\\
\Ext^1_w(p_{\bf k}(D),p_{\bf k}(D)) \arrow[r, "\kappa_w"] & \Hom(T,E) .
\end{tikzcd}
$$
In other words, the pullback $p_{\bf k}$ does not change the parameter $\psi$.
\end{enumerate}
\end{lem}

\begin{proof}
(i) 
Under the identification of $\Ext^1(D,D)$ with deformations of $D$ to $E[\vep]/\vep^2$ in Remark \ref{rem:triang_deform}, since $D$ is of regular integral weights ${\bf h}$ and $\Sp E$ is the only closed point of $\Sp(E[\vep]/\vep^2)$, for any $\widetilde{D}\in \Ext^1(D,D)$, the $(\varphi,\Gamma_K)$-module $\widetilde{D}$ belongs to $(\fX_n)^{\wedge}_{{\bf h}}(E[\vep]/\vep^2)$.
Recall from \cite[\S1.3, Definition 3.4, Lemma 3.10, Proposition 3.12, Appendix A]{Wu_translation} that there is an isomorphism 
$$
\Psi_{\bf h} = (f_{\bf h}, D_{\pdR}) : (\fX_n)^{\wedge}_{{\bf h}} \xrightarrow{\simeq} (\fX_n)^{\wedge}_{{\bf 0}}\times_{\frg/G}\widetilde{\frg}_{\bf h}/G
$$ 
of stacks over $\Rig_E$, such that for a $(\varphi,\Gamma_K)$-module $D_A$ over an affinoid rigid space $\Sp(A)$ whose Sen weights are ${\bf h}$ modulo the nilradical of $A$, $f_{{\bf h}}(D_A)$ is the unique $(\varphi,\Gamma_K)$-submodule of Sen weight ${\bf 0}$ in $D_A[1/t]$ such that $f_{{\bf h}}(D_A)[1/t]=D_A[1/t]$, and $D_{\pdR}(D_A)$ is the triple $(D_{\pdR}(D_{\dif}^+(D_A)), \nu_A, \Fil^\bullet D_{\pdR}(D_{\dif}^+(D_A))) \in (\widetilde{g}_{\bf h}/G)(\Sp(A))$.
Here, $D_{\pdR}(D_{\dif}^+(D_A))$ is the finite projective $(K\otimes_{\bbQ_p} A)$-module of rank $n$ defined by
$$
D_{\pdR}(D_\dif^+(D_A)) := (D_\dif^+(D_A)[1/t] \otimes_{(K_\infty \otimes_{\bbQ_p} A)[[t]]} (K_\infty \otimes_{\bbQ_p} A)[[t]][\log(t)])^{\Gamma_K} , 
$$
where $\Gamma_K$ acts on $\log(t)$ by $\gamma.\log (t)=\log (\vep(\gamma))+\log (t)$ for $\gamma \in \Gamma_K$, cf. \cite[Definition A.5]{Wu_translation}.
Then, $\nu_A: D_{\pdR}(D_\dif^+(D_A))\to D_{\pdR}(D_\dif^+(D_A))$ is the $A$-linear endomorphism induced by the $(K_\infty \otimes_{\bbQ_p} A)[[t]]$-linear derivation on $(K_\infty \otimes_{\bbQ_p} A)[[t]][\log(t)]$ satisfying $\nu_A(\log(t))=-1$, and $\Fil^\bullet D_{\pdR}(D_{\dif}^+(D_A))$ is the $\nu_A$-stable filtration
$$
\Fil^iD_{\pdR}(D_\dif^+(D_A)) := (t^i D_\dif^+(D_A)[\log(t)])^{\Gamma_K}, \quad \forall i\in \bbZ,
$$
on $D_{\pdR}(D_\dif^+(D_A))$ by projective $A$-submodules {\it of type {\bf h}}, cf. \cite[Definition 3.4(1)]{Wu_translation}, i.e., the $\sigma$-component has $A$-projective graded pieces with jumps occurring at $\{-h_{1,\sigma}<\dots<-h_{n,\sigma}\}$ for all $\sigma\in \Sigma_K$.
The change of weights map in \cite[\S1.3]{Wu_translation} between two regular integral weights ${\bf h}$ and ${\bf h'}$ is then defined to be the composite
$$
\begin{tikzcd}
(\fX_n)^{\wedge}_{{\bf h}} \arrow[r, "\Psi_{{\bf h}}"] \arrow[d, "{f_{{\bf h},{\bf h'}}}", dashed] & (\fX_n)^{\wedge}_{{\bf 0}}\times_{\frg/G}\widetilde{\frg}_{\bf h}/G \arrow[d, "{(\id, q_{{\bf h},{\bf h'}})}"] \\
(\fX_n)^{\wedge}_{{\bf h'}}                                                                        & (\fX_n)^{\wedge}_{{\bf 0}}\times_{\frg/G}\widetilde{\frg}_{\bf h'}/G \arrow[l, "\Psi_{{\bf h'}}^{-1}"]  \,,       
\end{tikzcd}
$$
where $q_{{\bf h},{\bf h'}} : \widetilde{\frg}_{\bf h}/G \to \widetilde{\frg}_{\bf h'}/G$ is the isomorphism of stacks mapping $(D_{\pdR,X},\nu_X,\Fil_{{\bf h},X}^\bullet)$ to the triple $(D_{\pdR,X},\nu_X,\Fil_{{\bf h'},X}^\bullet)$ for any $X\in \Rig_E$, cf. \cite[Definition 3.4(2)]{Wu_translation}, where $\Fil_{{\bf h'},X}^\bullet$ has the same underlying $\nu_X$-stable flag in $D_{\pdR,X}$ as $\Fil_{{\bf h},X}^\bullet$, except that the filtration is shifted such that its jumps occur at $\set{-h'_{i,\sigma}}{1\leq i\leq n, \sigma\in \Sigma_K}$.
From this description, it is clear that for any regular integral weights ${\bf h}, {\bf h'}$ and ${\bf h''}$, the change of weights maps satisfy $f_{{\bf h}, {\bf h''}} = f_{{\bf h'}, {\bf h''}}\circ f_{{\bf h},{\bf h'}}$.
Consequently, to show that $p_{\bf k}= \prod_{i,\sigma}p_{i,\sigma}^{k_{i,\sigma}}$ is well-defined on $\Ext^1(D,D)$, it suffices to show that $p_{i,\sigma}(\widetilde{D}) = f_{{\bf h}, {\bf h'}}(\widetilde{D})$ when ${\bf h'}$ is the mod-$\vep$ reduction of the Sen weights of $p_{i,\sigma}(\widetilde{D})$. 
By Theorem \ref{thm:Wu_3.16}(i), we need to show that $t_\sigma D_\dif^+(\widetilde{D}) \subset D_\dif^+(f_{{\bf h}, {\bf h'}}(\widetilde{D})) \subset D_\dif^+(\widetilde{D})$, and that $D_\dif^+(f_{{\bf h}, {\bf h'}}(\widetilde{D}))$ has the same Sen polynomial as that of $p_{i,\sigma}(\widetilde{D})$.

We write $A:=E[\vep]/\vep^2$ and $D_A:=\widetilde{D} \in \Ext^1(D,D)$. 
Let $V:= D_{\pdR}(D_\dif^+(D_A))$, which comes equipped with the $A$-linear endomorphism $\nu:V\to V$ and the $\nu$-stable filtration $\Fil^\bullet_{\bf h}(V)$.
For any $\tau\in \Sigma_K$, let $V_\tau = V \otimes_{K\otimes_{\bbQ_p}A,\tau\otimes 1} A$ be the $\tau$-component of $V$, which is a finite projective $A$-module of rank $n$ equipped with the endomorphism $\nu_\tau:V_\tau\to V_\tau$ and the $\nu_\tau$-stable filtration $\Fil^\bullet_{\bf h_\tau}(V)$.
By \cite[Proposition A.6(2)]{Wu_translation}, the $\Gamma_K$-representation $D_{\dif,\tau}^+(D_A)$ can be recovered by
$$
D_{\dif,\tau}^+(D_A) = \left( \sum_{k\in \bbZ} \Fil_{\bf h_\tau}^k(V_\tau) \otimes_{A} t^{-k}(K_\infty\otimes_{K,\tau} A)[[t]][\log(t)]\right)^{\nu_\tau = 0}.
$$
If we forget the $\Gamma_K$-action, then we have a $(K_\infty \otimes_{K,\tau} A)$-linear isomorphism
$$
\iota: \sum_{k\in \bbZ} \Fil^k_{{\bf h}_\tau}(V_\tau)\otimes_{A}  t^{-k}(K_\infty\otimes_{K,\tau} A)[[t]] \longrightarrow \left( \sum_{k\in \bbZ} \Fil_{\bf h_\tau}^k(V_\tau) \otimes_{A} t^{-k}(K_\infty\otimes_{K,\tau} A)[[t]][\log(t)]\right)^{\nu_\tau = 0}
$$
defined by $\iota(x) := \sum_{i\geq 0}\frac{1}{i!} \nu^i(x) \log(t)^i$ such that under the isomorphism $\iota$, the derivative action $\nabla_\tau$ of $\rm{Lie}(\Gamma_{K})$ on $D_{\dif,\tau}^+(D_A)$ (whose mod-$t$ reduction gives the $\tau$-Sen operator on $D_{\Sen,\tau}(D_A)$) becomes $(1\otimes_A t\frac{\partial}{\partial t} + \nu_\tau \otimes_A 1)$ on $\sum_{k\in \bbZ} \Fil^k_{{\bf h}_\tau}(V_\tau)\otimes_{A}  t^{-k}(K_\infty\otimes_{K,\tau} A)[[t]]$, cf. Step 1 of the proof of \cite[Proposition A.6]{Wu_translation}.
By Lemma \ref{lem:fil_gr}, we have via $\iota^{-1}$ the following identification
$$
D_{\dif,\tau}^+(D_A) =
\sum_{k\in \bbZ} \Fil^k_{{\bf h}_\tau}(V_\tau)\otimes_{A}  t^{-k}(K_\infty\otimes_{K,\tau} A)[[t]] = \bigoplus_{j=1}^n \Gr_{j,\tau} \otimes_A t^{h_{j,\tau}}(K_\infty\otimes_{K,\tau} A)[[t]]
$$
as $(K_\infty\otimes_{K,\tau} A)[[t]]$-modules, where $\Fil_{{\bf h}_\tau}^{-h_{j,\tau}}(V_\tau) = \Fil_{{\bf h}_\tau}^{-h_{j,\tau}+1}(V_\tau) \oplus \Gr_{j,\tau}$ is any chosen $A$-linear splitting such that $\Gr_{j,\tau} \cong \gr_{\Fil_{{\bf h}_\tau}}^{-h_{j,\tau}}(V_\tau) = \Fil_{{\bf h}_\tau}^{-h_{j,\tau}}(V_\tau)/\Fil_{{\bf h}_\tau}^{-h_{j,\tau}+1}(V_\tau)$ is projective of rank $1$ over $A$.
The $\tau$-Sen polynomial of $D_A$ is the characteristic polynomial of $[(1\otimes_A t\frac{\partial}{\partial t} + \nu_\tau \otimes_A 1) \pmod{t}]$.
Since $\nu_\tau$ stabilizes the filtration $\Fil_{\bf h_\tau}^\bullet(V_\tau)$, the matrix of $[(1\otimes_A t\frac{\partial}{\partial t} + \nu_\tau \otimes_A 1) \pmod{t}]$ is lower triangular.
Let $\overline{\nu}_{j,\tau} \in \End_A(\gr_{\Fil_{{\bf h}_\tau}}^{-h_{j,\tau}}(V_\tau))$ be the induced endomorphism on the nonzero graded piece $\gr_{\Fil_{{\bf h}_\tau}}^{-h_{j,\tau}}(V_\tau)$ of rank $1$.
Then, $\overline{\nu}_{j,\tau}$ is just multiplication by an element of $A$, which we again denote by $\overline{\nu}_{j,\tau}$.
Since $\nu_\tau$ is nilpotent by construction, $\overline{\nu}_{j,\tau} \in \Nil(A)$.
Since $t\frac{\partial}{\partial t}(t^{h})=ht^h$, it follows that
$
\prod_{j=1}^n (T-(h_{j,\tau}+\overline{\nu}_{j,\tau}))
$
is the $\tau$-Sen polynomial of $D_A$.
From the description of $f_{{\bf h},{\bf h'}}$ given above, we deduce again via $\iota^{-1}$ the following identification
$$
D_{\dif,\tau}^+(f_{{\bf h},{\bf h'}}(D_A))  = \bigoplus_{j=1}^n \Gr_{j,\tau} \otimes_A t^{h'_{j,\tau}}(K_\infty\otimes_{K,\tau} A)[[t]].
$$
Thus, we see that the $\tau$-Sen polynomial of $f_{{\bf h},{\bf h'}}(D_A)$ is 
$
\prod_{j=1}^n (T-(h'_{j,\tau}+\overline{\nu}_{j,\tau}))
$, $t D_{\dif,\sigma}^+(D_A) \subset D_{\dif,\sigma}^+(f_{{\bf h},{\bf h'}}(D_A)) \subset D_{\dif,\sigma}^+(D_A)$, and $D_{\dif,\tau}^+(f_{{\bf h},{\bf h'}}(D_A)) = D_{\dif,\tau}^+(D_A)$ for all $\tau \in \Sigma_K\setminus \{\sigma\}$.
Then, $p_{i,\sigma}(\widetilde{D}) = f_{{\bf h}, {\bf h'}}(\widetilde{D})$ by Theorem \ref{thm:Wu_3.16}(i).
Consequently, $p_{\bf k}=\prod_{i,\sigma}p_{i,\sigma}^{k_{i,\sigma}}$ is well-defined as a map from $\Ext^1(D,D)$ to $\Ext^1(p_{\bf k}(D), p_{\bf k}(D))$.

\vskip8pt

(ii) To verify these properties, by induction we may assume $p_{\bf k} = p_{i,\sigma}$ for some $i$ and $\sigma$, as in Theorem \ref{thm:Wu_3.16}(i).
By Theorem \ref{thm:Wu_3.16}(iii), $p_{i,\sigma}$ is bijective.
As for $E$-linearity, for $\widetilde{D}\in \Ext^1(D,D)$, we have a commutative diagram with exact rows
$$
\begin{tikzcd}
0 \arrow[r] & D \arrow[r] & \widetilde{D} \arrow[r] & D \arrow[r] & 0 \\ 
0 \arrow[r] & p_{i,\sigma}(D) \arrow[r] \arrow[u,hook] & p_{i,\sigma}(\widetilde{D}) \arrow[r] \arrow[u,hook] & p_{i,\sigma}(D) \arrow[r] \arrow[u,hook] & 0,
\end{tikzcd}
$$
where the vertical arrows are inclusions.
For $a\in E$, applying pushout by multiplication by $a$ on the subobjects, we see that $p_{i,\sigma}(a[\widetilde{D}]) = a[p_{i,\sigma}(\widetilde{D})]$.
As for additivity, we need to show that for any $\widetilde{D}_1,\widetilde{D}_2\in \Ext^1(D,D)$ of the form
$$
0\to D \xrightarrow{i_1} \widetilde{D}_1 \xrightarrow{\pi_1} D \to 0, \quad 0\to D \xrightarrow{i_2} \widetilde{D}_2 \xrightarrow{\pi_2} D \to 0,
$$
the Baer sums satisfy 
$$
p_{i,\sigma}(\widetilde{D}_1+\widetilde{D}_2) = p_{i,\sigma}(\widetilde{D}_1)+p_{i,\sigma}(\widetilde{D}_2).
$$
Since both sides are $(\varphi,\Gamma_K)$-submodules of the Baer sum $\widetilde{D}_1+\widetilde{D}_2$ containing $t_\sigma(\widetilde{D}_1+\widetilde{D}_2)$, by Theorem \ref{thm:Wu_3.16}(i) it is enough to show that their Sen polynomials are equal.
For $\widetilde{D}\in \Ext^1(D,D)$, its Sen polynomial $P_{\Sen,\widetilde{D}}(T)\in (E[\vep]/\vep^2)[T]$ reduces modulo $\vep$ to the Sen polynomial $P_{\Sen, D}(T) = \prod_{\tau\in\Sigma_K}\prod_{j=1}^n (T-h_{j,\tau}) \in \prod_{\tau\in\Sigma_K} E[T] \cong (K\otimes_{\bbQ_p}E)[T]$ of $D$.
Since $E[\vep]/\vep^2$ is $(\vep)$-adically complete, Hensel's lemma applies to show that 
$$
P_{\Sen,\widetilde{D}}(T) = \prod_{\tau\in \Sigma_K}\prod_{j=1}^n (T-(h_{j,\tau}+a_{j,\tau}\vep))
$$
for some $a_{j,\tau}\in E$.
Suppose that for $l=1,2$, $P_{\Sen,\widetilde{D}_l}(T) = \prod_{\tau\in \Sigma_K}\prod_{j=1}^n (T-(h_{j,\tau}+a_{j,\tau}^{(l)}\vep))$ for some $a_{j,\tau}^{(l)}\in E$.
We show that $P_{\Sen, \widetilde{D}_1+\widetilde{D}_2}(T) = \prod_\tau\prod_j (T-(h_{j,\tau}+(a_{j,\tau}^{(1)}+a_{j,\tau}^{(2)})\vep))$, from which it would follow that $p_{i,\sigma}(\widetilde{D}_1+\widetilde{D}_2) = p_{i,\sigma}(\widetilde{D}_1)+p_{i,\sigma}(\widetilde{D}_2)$.
Since $D_{\Sen,\tau}$ is exact, for $l=1,2$, we can write $D_{\Sen,\tau}(\widetilde{D}_l) = \vep D_{\Sen,\tau}(D) \oplus D_{\Sen,\tau}(D)$ as free $K_\infty \otimes_{K,\tau} (E[\vep]/\vep^2)$-module of rank $n$. 
Let $\{e_{1,\tau},\dots, e_{n,\tau}\}$ be an eigenbasis of the $\tau$-Sen operator $\Theta_{\tau}$ on $D_{\Sen,\tau}(D)$.
Then, $\{e_{1,\tau},\dots, e_{n,\tau}\}$ is a $K_\infty \otimes_{K,\tau} (E[\vep]/\vep^2)$-basis of $D_{\Sen,\tau}(\widetilde{D}_l)$ on which the $\tau$-Sen operator $\widetilde{\Theta}_{\tau}$ acts by $\widetilde{\Theta}_{\tau}(e_{j,\tau})=(h_{j,\tau}+a_{j,\tau}^{(l)}\vep)e_{j,\tau}$.
The $\tau$-Sen module of $\widetilde{D}_1+\widetilde{D}_2$ is a Baer sum of $\tau$-Sen modules:
\begin{align*}
&D_{\Sen,\tau}(\widetilde{D}_1+\widetilde{D}_2) 
= D_{\Sen,\tau}(\widetilde{D}_1)+ D_{\Sen,\tau}(\widetilde{D}_2) 
\\&= 
\set{(x_1,x_2)\in D_{\Sen,\tau}(\widetilde{D}_1)\oplus D_{\Sen,\tau}(\widetilde{D}_2)}{\pi_1(x_1)=\pi_2(x_2)}/\set{\vep r,-\vep r}{r\in D_{\Sen,\tau}(D)},
\end{align*}
which has a $K_\infty \otimes_{K,\tau} (E[\vep]/\vep^2)$-basis represented by the classes of $\set{(e_{j,\tau},e_{j,\tau})}{1\leq j\leq n}$ on which the $\tau$-Sen operator $\widetilde{\Theta}_{\tau}$ acts by
\begin{align*}
\widetilde{\Theta}_{\Sen,\tau}([e_{j,\tau},e_{j,\tau}]) 
&=
[(h_{j,\tau}+a_{j,\tau}^{(1)}\vep)e_{j,\tau}, (h_{j,\tau}+a_{j,\tau}^{(2)}\vep)e_{j,\tau}]
\\&=
[(h_{j,\tau}+(a_{j,\tau}^{(1)}+a_{j,\tau}^{(2)})\vep) e_{j,\tau}, h_{j,\tau}e_{j,\tau}]
\\&=
(h_{j,\tau}+(a_{j,\tau}^{(1)}+a_{j,\tau}^{(2)})\vep)[e_{j,\tau},e_{j,\tau}].
\end{align*}
Hence, $P_{\Sen, \widetilde{D}_1+\widetilde{D}_2}(T) = \prod_\tau\prod_j (T-(h_{j,\tau}+(a_{j,\tau}^{(1)}+a_{j,\tau}^{(2)})\vep))$, and the isomorphism
$$
p_{i,\sigma}: \Ext^1(D,D) \xrightarrow{\sim} \Ext^1(p_{i,\sigma}(D),p_{i,\sigma}(D))
$$ 
is $E$-linear. 
For each trianguline deformation $\widetilde{D}\in \Ext^1_w(D,D)$, $p_{i,\sigma}(\widetilde{D})$ is given by a pullback using the unique lifting of $\Fil_w^\bullet(D)$ and belongs to $\Ext^1_w(p_{i,\sigma}(D),p_{i,\sigma}(D))$ by construction. 
Moreover, the resulting triangulation on $p_{i,\sigma}(\widetilde{D})$ has parameters $\delta_{w,j}(1+\psi_j\vep)$ for $j\leq n-i$ and $x_\sigma \delta_{w,j}(1+\psi_j\vep)$ for $j>n-i$, from which we see that $\kappa_w(p_{i,\sigma}(\widetilde{D}))= (\psi_1,\dots,\psi_n) = \kappa_w(\widetilde{D})$.
\end{proof}

Consequently, $p_{\bf k}$ preserves the subspace $\Ext^1_0(D,D)$ and induces an isomorphism 
$$
p_{\bf k}: \overline{\Ext}^1(D,D) \xrightarrow{\sim} \overline{\Ext}^1(p_{\bf k}(D), p_{\bf k}(D)).
$$
Thus, we have a natural commutative diagram 
\begin{equation}
\label{eqn:diagram_5.1.7}
\begin{tikzcd}
\bigoplus_{w\in S_n} \overline{\Ext}^1_w(D,D) \arrow[r, two heads] \arrow[d, "p_{\bf k}"] & \overline{\Ext}^1(D,D) \arrow[d, "p_{\bf k}"]\\
\bigoplus_{w\in S_n} \overline{\Ext}^1_w(p_{\bf k}(D),p_{\bf k}(D)) \arrow[r, two heads] & \overline{\Ext}^1(p_{\bf k}(D),p_{\bf k}(D)).
\end{tikzcd}
\end{equation}

\subsubsection{Translation} 

We fix notations for weights and translation functor on the automorphic side.
Recall that $G=\bfG(K) = \GL_n(K) \supset B = \bfB(K) \supset T = \bfT(K)$ and $B^{-}=\bfB^-(K)$.
Let 
$$
\rho := \left( \frac{n-1}{2}, \frac{n-3}{2},\dots, \frac{1-n}{2} \right)_{\sigma\in\Sigma_K} \in E^{\Sigma_K}
$$
be the half-sum of positive roots of $G$ with respect to $B$.
Then, the half-sum of the positive roots with respect to $B^-$ is $\rho^- := -\rho$.
The Weyl group of $\Res^K_{\bbQ_p}(\bfG)$ is $W\cong S_n^{\Sigma_K}$, whose element $w=(w_\sigma)_{\sigma}$ acts on $\frt_E^* = \prod_{\sigma\in \Sigma_K} \frt^* \otimes_{K,\sigma} E$ componentwise. 
Let
$
\theta := (0, -1, \dots, 1-n)_{\sigma\in \Sigma_K} \in \bbZ^{\Sigma_K},
$
which differs from $\rho$ by a central character $\frac{n-1}{2}(1,\dots,1)$.
Thus, we can compute the dot-action 
$w\cdot \lambda := w(\lambda+\rho)-\rho =w(\lambda+\theta)-\theta$ of the Weyl group on weights $\lambda\in \frt_E^*$ using either $\rho$ or $\theta$.

For locally $\bbQ_p$-analytic representations of $\GL_n(K)$ on compact type $E$-vector spaces (all representations later in this article are of compact type), \cite[Lemma 3.9]{Ding_wallcross} relates the translation functor on those $V$ that admit an infinitesimal character of $\frz_E$ with the translation functor on the strong dual $V'_b$:
the following diagram commutes for any integral weights $\mu,\lambda \in \frt_E^*$, 
\begin{equation}
\label{eqn:diagram_Ding3.9}
\begin{tikzcd}
D(G)\text{-mod}_{\frz_E = \chi_{\lambda^*}} \arrow[r, "T_{\lambda^*}^{\mu^*}"] & D(G)\text{-mod}_{\frz_E = \chi_{{\mu}^*}} \\
\Rep^{\la}(G)_{\frz_E=\chi_\lambda} \arrow[r, "T_{\lambda}^\mu"] \arrow[u, "(-)'_b"]& \Rep^{\la}(G)_{\frz_E=\chi_\mu} \arrow[u, "(-)'_b"]\,,
\end{tikzcd}
\end{equation}
where the vertical maps are taking the strong dual, which induce an anti-equivalence of categories from locally analytic representations of $G$ on compact type $E$-vector spaces to separately continuous $D(G)$-modules on nuclear Fr\'{e}chet spaces, cf. \cite[Cor. 3.3]{Schneider-Teitelbaum_GL2}, and $\lambda^* = -w_0\lambda$ is the dual highest weight with $w_0\in W$ the longest Weyl group element, cf. \cite[VIII, \S7, n\textsuperscript{o} 5, prop. 11]{Bourbaki_Lie78}.
We will use the same compatibility for objects with bounded generalized infinitesimal character. Namely, if a compact-type locally analytic representation $V$ is killed by a power $\mathfrak{m}_{\chi_\lambda}^N$ of the maximal ideal corresponding to $\chi_\lambda$, then its strong dual $V_b^{\prime}$ is killed by the corresponding power $\mathfrak{m}_{\chi_{-w_0\lambda}}^N$. Since strong duality and translation functors are exact on the relevant $\mathfrak{z}_E$-finite categories, \cite[Lemma 3.9]{Ding_wallcross} extends to this bounded generalized setting by dévissage on $N$. In the applications below, $N \leq 2$, because the principal series are deformations over $E[\varepsilon] / \varepsilon^2$.

\begin{lem}
\label{lem:T_equivalence}
Let $\lambda':={\bf h'}-\theta$ and $\lambda:={\bf h}-\theta$.
The translation functors $T_{\lambda^*}^{\lambda'^*}$ and  $T_{\lambda}^{\lambda'}$ in \eqref{eqn:diagram_Ding3.9} are equivalences of categories.
\end{lem}
\begin{proof}
It is enough to show that $T_{\lambda^*}^{\lambda'^*}:D(G)\text{-mod}_{|\lambda^*|}\longrightarrow D(G)\text{-mod}_{|\lambda'^*|} $ is an equivalence of categories, for which we use \cite[Theorem 1]{JLS_translation}.
Note that $T_{\lambda^*}^{\lambda'^*} = T_{w_0\cdot\lambda^*}^{w_0\cdot \lambda'^*}$ because $w_0\cdot \lambda'^*-w_0\cdot \lambda^* = w_0(\lambda'^*-\lambda^*)$.
Then, $w_0\cdot \lambda'^*-w_0\cdot \lambda^* = \lambda - \lambda'$ lifts to an algebraic character of $\bfT$, both $w_0\cdot \lambda'^*$ and $w_0\cdot \lambda^*$ are anti-dominant, and both of them have the same (trivial) stabilizer for the dot-action as ${\bf h},{\bf h'}$ are regular.
Hence, \cite[Theorem 1]{JLS_translation} applies to yield that 
$$
T_{w_0\cdot \lambda^*}^{w_0\cdot \lambda'^*} = T_{\lambda^*}^{\lambda'^*}:D(G)\text{-mod}_{|\lambda^*|}\longrightarrow D(G)\text{-mod}_{|\lambda'^*|} 
$$
is an equivalence of categories.
\end{proof}

\subsubsection{Principal series}
We recall the deformations to $E[\vep]/\vep^2$ on the locally analytic side considered in \cite[\S3.1.2]{Ding_crystabelline}.
For any $w\in S_n$, we define the locally analytic principal series
$$
\PS(w(\phi),{\bf h}) := \Ind_{B^-}^{G}(w(\phi)\eta z^\lambda)^\la,
$$
where $z^\lambda: T \to E^\times, (z_1,\dots,z_n)\mapsto \prod_{i=1}^n\prod_{\sigma\in \Sigma_K} \sigma(z_i)^{\lambda_{i,\sigma}}$ is the algebraic character of weight $\lambda$, and 
$
\eta:= 1 \boxtimes |\cdot|_K \boxtimes \dots \boxtimes |\cdot|_K^{n-1} = |\cdot|^{-1}_K \circ \theta
$
is a smooth character added for the purpose of normalization.
By \cite[Proposition 3.1(2)]{Ding_crystabelline}, the $\GL_n(K)$-socle of $\PS(w(\phi),{\bf h})$ is independent of the choice of $w\in S_n$, and the socle is locally algebraic, given by 
$$
\pi_\alg(\phi,{\bf h}) := \Ind_{B^-}^G(\phi\cdot\eta)^{\infty} \otimes_E L(\lambda),
$$
where $\Ind_{B^-}^G(-)^{\infty}$ denotes the smooth parabolic induction, and $L(\lambda)$ is the algebraic representation of $G$ of highest weight $\lambda$ with respect to $B$.
For $\psi\in \Hom(T,E)$, let $i_{w,{\bf h}}(\psi)$ be the locally analytic self-extension of $\PS(w(\phi),{\bf h})$ given by the locally analytic principal series
$$
i_{w,{\bf h}}(\psi):= \Ind_{B^-}^G(w(\phi)\eta z^\lambda(1+\psi\vep))^\la 
$$
with coefficients in $E[\vep]/\vep^2$, cf. \S\ref{subsection:translation_foundation}, whose strong dual is isomorphic to 
$$
D(G)\otimes_{D(B^-)} (E[\vep]/\vep^2)_{(w(\phi)\eta z^\lambda(1+\psi\vep))^{-1}}
$$
by Lemma \ref{lem:ST_dual} and Remark \ref{rem:dual_E[vep]}.
From this tensor product description, we deduce that
$$
i_{w,{\bf h}}: \Hom(T,E)\to \Ext^1_G(\PS(w(\phi),{\bf h}),\PS(w(\phi),{\bf h}))
$$
is an $E$-linear map.

\begin{prop}
\label{prop:translate_principal_series}
For any $w\in S_n$ and $\psi\in \Hom(T,E)$, we have
$
T_\lambda^{\lambda'}(i_{w,{\bf h}}(\psi)) \cong  i_{w,{\bf h'}}(\psi).
$
That is, the following diagram is commutative:
\begin{equation}
\label{eqn:diagram_5.2.1}
\begin{tikzcd}
\Hom(T,E) \arrow[r, "{i_{w,{\bf h}}}"]\arrow[d, Rightarrow, no head, "\id"] & \Ext^1_G(\PS(w(\phi),{\bf h}),\PS(w(\phi),{\bf h})) \arrow[d, "T_\lambda^{\lambda'}"] \\
\Hom(T,E) \arrow[r, "{i_{w,{\bf h'}}}"] & \Ext^1_G(\PS(w(\phi),{\bf h'}), \PS(w(\phi),{\bf h'})).
\end{tikzcd}
\end{equation}
\end{prop}
\begin{proof}
By \eqref{eqn:diagram_Ding3.9}, it suffices to prove the statement for the strong duals as $D(G)$-modules:
$$
T_{\lambda^*}^{\lambda'^*} \left(D(G)\otimes_{D(B^-)} (E[\vep]/\vep^2)_{(w(\phi)\eta z^\lambda(1+\psi\vep))^{-1}}\right) \cong D(G)\otimes_{D(B^{-})} (E[\vep]/\vep^2)_{(w(\phi)\eta z^{\lambda'}(1+\psi\vep))^{-1}},
$$
which we will reduce to Lemma \ref{lem:translate_D(G)_induction}.

For this purpose, it is necessary to translate the language from the upper triangular Borel $B$ to lower triangular Borel $B^-$.
For $w\in W$ and $\xi \in \frt_E^*$, the Weyl dot-action with respect to $B^-$ is given by $w \overline{\cdot} \xi = w(\xi + \rho^-)-\rho^-$.
By the independence of the twisted Harish-Chandra isomorphism from the choice of the system of positive roots \cite[p. 290]{KV_coh_induction} and the definition of the twisted Harish-Chandra isomorphism, one verifies (using the positive roots with respect to $B$ and the positive roots with respect to $B^-$) that for any weight $\xi\in \frt_E^*$, the eigencharacter $\chi_{\xi}$ of $\frz_E$ on the Verma module $U(\frg_E)\otimes_{U(\frb_E)} E_{\xi}$ coincides with the eigencharacter $\chi^-_{w_0\xi}$ of $\frz_E$ on the {\it opposite} Verma module $U(\frg_E)\otimes_{U(\frb_E^-)} E_{w_0\xi}$.
For any $\xi_1,\xi_2\in \frt_E^*$, let 
$\overline{T}_{\xi_1}^{\xi_2}$ denote the translation functor with respect to $B^-$.
Then, we have $\overline{T}_{\xi_1}^{\xi_2} = T_{w_0\xi_1}^{w_0\xi_2}$.

Therefore, we have
$
T_{\lambda^*}^{\lambda'^*} = T_{-w_0\lambda}^{-w_0\lambda'} = \overline{T}_{-\lambda}^{-\lambda'} = \overline{T}_{w_0\overline{\cdot}(-\lambda)}^{w_0\overline{\cdot}(-\lambda')}
$ 
and we need to show that
\begin{align*}
& \overline{T}_{w_0\overline{\cdot}(-\lambda)}^{w_0\overline{\cdot}(-\lambda')} \left( D(G)\otimes_{D(B^-)} (E[\vep]/\vep^2)_{w_0 \overline{\cdot} w_0 \overline{\cdot}(w(\phi)^{-1}\eta^{-1} z^{-\lambda}(1-\psi\vep))} \right)
\\\cong & \mbox{ }
D(G)\otimes_{D(B^-)} (E[\vep]/\vep^2)_{w_0 \overline{\cdot} w_0 \overline{\cdot}(w(\phi)^{-1}\eta^{-1} z^{-\lambda'}(1-\psi\vep))},
\end{align*}
which will follow from Lemma \ref{lem:translate_D(G)_induction} by taking $A = E[\vep]/\vep^2$, $\widetilde{\lambda}_A = w_0 \overline{\cdot}(w(\phi)^{-1}\eta^{-1} z^{-\lambda}(1-\psi\vep))$, $\widetilde{\mu}_A = w_0 \overline{\cdot}(w(\phi)^{-1}\eta^{-1} z^{-\lambda'}(1-\psi\vep))$, and $w= w_0 \in W_{[d\widetilde{\lambda}_A \text{ mod }\vep]} = W_{[d\widetilde{\mu}_A \text{ mod }\vep]} = W$.
It remains to check that the hypotheses of Lemma \ref{lem:translate_D(G)_induction} are satisfied by the above choices.
For (i), since ${\bf h}$ and ${\bf h'}$ are regular, both $w_0\overline{\cdot}(-\lambda)=-w_0({\bf h})+\theta$ and $w_0\overline{\cdot}(-\lambda')=-w_0({\bf h'})+\theta$ are anti-dominant with respect to $B^-$.
For (ii), the ratio 
\begin{align*}
\widetilde{\mu}_A / \widetilde{\lambda}_A 
&= 
\frac{w_0 \overline{\cdot}(w(\phi)^{-1}\eta^{-1} z^{-\lambda'}(1-\psi\vep))}{w_0 \overline{\cdot}(w(\phi)^{-1}\eta^{-1} z^{-\lambda}(1-\psi\vep))}
\\&=
w_0(z^{\lambda - \lambda'}) 
= z^{w_0(\lambda-\lambda')}
\end{align*}
is an $E^\times$-valued algebraic character of $\bfT(K)$.
For (iii), both $(w_0\overline{\cdot}(-\lambda))^\natural$ and $(w_0\overline{\cdot}(-\lambda'))^\natural$ lie in the same open $B^-$-anti-dominant chamber in the Euclidean space $\bbR\otimes_{\bbZ} \Phi$, for $\Phi$ the root system of $(\frg_E,\frt_E)$, cf. \cite[\S4.2.3]{JLS_translation}, because ${\bf h}$ and ${\bf h'}$ are regular.
Hence, Lemma \ref{lem:translate_D(G)_induction} applies.
\end{proof}

Recall from \cite[Proposition 3.1]{Ding_crystabelline} that the irreducible components of $\PS(w(\phi),{\bf h})$ are
$$
\cF_{B^-}^G(L^{-}(-u\cdot \lambda), w(\phi)\eta) = \left(D(G)\otimes_{D(\frg,B^-)} (L^{-}(-u\cdot \lambda)\otimes_E (w(\phi)\eta)^{-1}) \right)'_b\,,
$$
cf. \cite[Remark 4.1.11]{JLS_translation}, for $w\in S_n$ and $u=(u_\sigma)_\sigma \in W = S_n^{\Sigma_K}$, where $L^{-}(-u\cdot \lambda)$ is the unique simple quotient of $U(\frg_E)\otimes_{U(\frb_E^-)} E_{-u\cdot \lambda}$, \cite[\S3.1.1]{Ding_crystabelline}.
We denote these representations by $\scrC^{(\phi,{\bf h})}(w,u) := \cF_{B^-}^G(L^{-}(-u\cdot \lambda), w(\phi)\eta)$ to emphasize their dependence on $\phi$ and the weight ${\bf h}$.
By Lemma \ref{lem:T_equivalence} and \cite[Proposition 4.2.10]{JLS_translation}, we have 
\begin{equation}
\label{eqn:dagger_I}
T_\lambda^{\lambda'}(\scrC^{(\phi,{\bf h})}(w,u)) \cong \scrC^{(\phi,{\bf h'})}(w,u)
\end{equation}
since $L^{-}(-u\cdot \lambda)$ is the unique simple quotient of $U(\frg_E)\otimes_{U(\frb_E^-)} E_{-u\cdot \lambda}$ and $T_\lambda^{\lambda'}$ is an equivalence of categories.
Then, following \cite[(3.2)]{Ding_crystabelline}, we let $\PS_1(w(\phi),{\bf h})$ be the unique subrepresentation of $\PS(w(\phi),{\bf h})$ of socle $\pi_\alg(\phi,{\bf h}) = \scrC^{(\phi,{\bf h})}(w,\id)$ and cosocle $\bigoplus_{\substack{1\leq i\leq n-1 \\ \sigma\in\Sigma_K}} \scrC^{(\phi,{\bf h})}(w,s_{i,\sigma})$, where $s_{i,\sigma}$ denotes the element of $W=S_n^{\Sigma_K}$ that is the simple reflection $s_i$ at the $\sigma$-component and trivial elsewhere.
By its uniqueness, \eqref{eqn:dagger_I}, Lemma \ref{lem:T_equivalence} and Proposition \ref{prop:translate_principal_series}, we deduce
\begin{equation}
\label{eqn:dagger_II}
T_\lambda^{\lambda'}(\PS_1(w(\phi),{\bf h})) \cong \PS_1(w(\phi),{\bf h'}).
\end{equation}
By \cite[(3.3)]{Ding_crystabelline}, the amalgamated sum $\bigoplus_{\pi_{\alg}(\phi,{\bf h})}^{w\in S_n} \PS_1(w(\phi),{\bf h})$ has a unique subquotient of socle $\pi_{\alg}(\phi,{\bf h})$, and this subquotient is denoted $\pi_1(\phi,{\bf h})$.
By its uniqueness, \eqref{eqn:dagger_I}, \eqref{eqn:dagger_II} and Lemma \ref{lem:T_equivalence}, we deduce that
\begin{equation}
\label{eqn:dagger_III}
T_\lambda^{\lambda'}(\pi_1(\phi,{\bf h})) \cong \pi_1(\phi,{\bf h'}).
\end{equation}
After \cite[Remark 3.9]{Ding_crystabelline}, we denote by $\pi(\phi,{\bf h})$ the unique quotient of $\bigoplus_{\pi_{\alg}(\phi,{\bf h})}^{w\in S_n} \PS_1(w(\phi),{\bf h})$ of socle $\pi_{\alg}(\phi,{\bf h})$.
We similarly deduce that
\begin{equation}
\label{eqn:dagger_IV}
T_\lambda^{\lambda'}(\pi(\phi,{\bf h})) \cong \pi(\phi,{\bf h}').
\end{equation}

Following \cite[(3.11)]{Ding_crystabelline}, for each $w\in S_n$ we define 
$$
\zeta_{w,{\bf h},1}: \Hom(T,E) \to \Ext^1_G(\pi_{\alg}(\phi,{\bf h}), \pi_1(\phi, {\bf h})),
$$
firstly sending $\psi\in \Hom(T,E)$ to $i_{w,{\bf h}}(\psi) \in \Ext^1_G(\PS(w(\phi),{\bf h}), \PS(w(\phi),{\bf h}))$, and then pulling it back along the embedding $\pi_\alg(\phi,{\bf h}) \hookrightarrow \PS_1(w(\phi),{\bf h})$ to 
$$
\Ext^1_G(\pi_{\alg}(\phi,{\bf h}), \PS_1(w(\phi),{\bf h}))\cong \Ext^1_G(\pi_{\alg}(\phi,{\bf h}), \PS(w(\phi),{\bf h})),
$$ 
where the isomorphism is induced by pushout along $\PS_1(w(\phi),{\bf h})\hookrightarrow \PS(w(\phi),{\bf h})$, and finally pushing it forward to $\Ext^1_G(\pi_\alg(\phi,{\bf h}),\pi_1(\phi,{\bf h}))$ along $\PS_1(w(\phi),{\bf h}) \hookrightarrow\pi_1(\phi,{\bf h})$.
The same construction, using $\pi(\phi,{\bf h})$ instead of $\pi_1(\phi,{\bf h})$, gives the version
$$
\zeta_{w,{\bf h}}: \Hom(T,E) \to \Ext^1_G(\pi_{\alg}(\phi,{\bf h}),\pi(\phi,{\bf h})).
$$
By \cite[Remark 3.9]{Ding_crystabelline}, the pushforward along $\pi_1(\phi,{\bf h})\hookrightarrow \pi(\phi,{\bf h})$ is an isomorphism 
$$
\Ext^1_G(\pi_{\alg}(\phi,{\bf h}),\pi_1(\phi,{\bf h})) \cong \Ext^1_G(\pi_{\alg}(\phi,{\bf h}),\pi(\phi,{\bf h})).
$$
By \cite[(3.17)]{Ding_crystabelline}, both $\zeta_{w,{\bf h},1}$ and $\zeta_{w,{\bf h}}$ are $E$-linear injections, and we denote their images by
$$
\Ext^1_w(\pi_{\alg}(\phi,{\bf h}),\pi_1(\phi,{\bf h})) \cong \Ext^1_w(\pi_\alg(\phi,{\bf h}),\pi(\phi,{\bf h})).
$$
Since $T_\lambda^{\lambda'}$ is an equivalence of categories, it commutes with pullbacks and pushforwards.
Thus, together with Lemma \ref{lem:pullback_parameter} and Proposition \ref{prop:translate_principal_series}, we deduce the following corollary.
\begin{cor}
\label{cor:Ding_Ext-link}
\begin{enumerate}[(i)]
\item 
For any $w\in S_n$ and $\bullet=1$ or empty, we have a commutative diagram
$$
\begin{tikzcd}
\Hom(T,E) \arrow[d, Rightarrow, no head, "\id"] \arrow[r, hook, "{\zeta_{w,{\bf h},\bullet}}"] & \Ext^1_G(\pi_{\alg}(\phi,{\bf h}),\pi_\bullet(\phi,{\bf h})) \arrow[d, "{T_{\lambda}^{\lambda'}}"] \\
\Hom(T,E) \arrow[r, hook, "\zeta_{w,{\bf h'},\bullet}"] &  \Ext^1_G(\pi_{\alg}(\phi,{\bf h'}),\pi_\bullet(\phi,{\bf h'})).
\end{tikzcd}
$$
\item 
For any $w\in S_n$ and $\bullet=1$ or empty, we have a commutative diagram
$$
\begin{tikzcd}
\overline{\Ext}^1_w(D,D) \arrow[d, "p_{\bf k}"] \arrow[r, hook, "{\zeta_{w,{\bf h},\bullet}}\circ \kappa_w"] & \Ext^1_G(\pi_{\alg}(\phi,{\bf h}),\pi_\bullet(\phi,{\bf h})) \arrow[d, "{T_{\lambda}^{\lambda'}}"] \\
\overline{\Ext}^1_w(p_{\bf k}(D),p_{\bf k}(D)) \arrow[r, hook, "\zeta_{w,{\bf h'},\bullet}\circ \kappa_w"] &  \Ext^1_G(\pi_{\alg}(\phi,{\bf h'}),\pi_\bullet(\phi,{\bf h'})).
\end{tikzcd}
$$
\item 
For $\bullet=1$ or empty, we have a commutative diagram
$$
\begin{tikzcd}
\bigoplus_{w\in S_n} \overline{\Ext}^1_w(D,D) \arrow[r, "p_{\bf k}"] \arrow[d, "\zeta_{w,{\bf h},\bullet}\circ \kappa_w"] & \bigoplus_{w\in S_n} \overline{\Ext}^1_w(p_{\bf k}(D),p_{\bf k}(D)) \arrow[d, "\zeta_{w,{\bf h'},\bullet}\circ \kappa_w"] \\
\bigoplus_{w\in S_n} \Ext_w^1(\pi_{\alg}(\phi,{\bf h}),\pi_\bullet(\phi,{\bf h})) \arrow[r, "T_{\lambda}^{\lambda'}"] \arrow[d, "\mathrm{sum}", two heads]& \bigoplus_{w\in S_n} \Ext_w^1(\pi_{\alg}(\phi,{\bf h'}),\pi_\bullet(\phi,{\bf h'})) \arrow[d, "\mathrm{sum}", two heads] \\
\Ext_G^1(\pi_{\alg}(\phi,{\bf h}),\pi_\bullet(\phi,{\bf h})) \arrow[r, "T_\lambda^{\lambda'}"] & \Ext_G^1(\pi_{\alg}(\phi,{\bf h'}),\pi_\bullet(\phi,{\bf h'})),
\end{tikzcd}
$$
where the amalgamation maps ``sum" are surjective by \cite[Proposition 3.8(2)]{Ding_crystabelline}.
\end{enumerate}
\end{cor}

\subsubsection{Universal extension}
\label{univ_ext}
We discuss the notion of ``universal extension" in general.
Let $E$ be a field, and let $\cD$ be an $E$-linear abelian category with objects $A$ and $B$.
Then, to any finite dimensional $E$-linear subspace $W$ of $\Ext^1_{\cD}(A,B)$, one can attach a ``universal extension" $\cE^{W\text{-}\univ}$ of $A \otimes_E W$ ($\cong A^{\dim_E W}$) by $B$, satisfying the property that for any extension
$$
0 \to B \to \cE_e \to A \to 0
$$
with corresponding class $e\in W \subset \Ext^1_\cD(A,B)$, the map $\alpha_e: A \to A\otimes_E W, a\mapsto a\otimes e$ induces
$$
\alpha_e^*:\Ext^1_{\cD}(A\otimes_E W,B) \to \Ext^1_{\cD}(A,B)
$$
such that the pullback $\alpha_e^*(\cE^{W\text{-}\univ})$ is precisely $\cE_e$.

\begin{prop}
\label{prop:universal_ext_exists}
\begin{enumerate}[(i)]
\item 
Abstractly, this extension $\cE^{W\text{-}\univ}$ is the image of the inclusion map 
$$
i_W: W \hookrightarrow \Ext^1_{\cD}(A,B)
$$
under the canonical isomorphisms
$$
i_W\in 
\Hom_E(W,\Ext^1_{\cD}(A,B)) \xleftarrow[\rm can]{\sim} \Ext^1_{\cD}(A,B) \otimes_E W^\vee \xrightarrow[\rm can]{\sim} \Ext^1_{\cD}(A\otimes_E W,B) \ni \cE^{W\text{-}\univ},
$$
where the first map is a canonical isomorphism for finite dimensional $W$, and the second map is defined as follows: given $\cE_e \in \Ext^1_{\cD}(A,B)$ and a linear functional $f \in W^\vee$, we pullback $\cE_e$ via $\beta_f: A \otimes_E W \to A, a\otimes w \mapsto f(w)a$ to $\beta_f^*(\cE_e) \in \Ext^1_{\cD}(A,B)$, which defines
$$
\Ext^1_{\cD}(A,B) \otimes_E W^\vee \xrightarrow[\rm can]{\sim} \Ext^1_{\cD}(A\otimes_E W,B), \quad \cE_e \otimes f \mapsto \beta_f^*(\cE_e)
$$
with its inverse given by
$$
\Ext^1_{\cD}(A\otimes_E W,B) \xrightarrow[\rm can]{\sim} \Hom_E(W,\Ext^1_{\cD}(A,B)), \quad \cE \mapsto (e \mapsto \alpha_e^*(\cE)).
$$
\item 
Concretely, choose any $E$-basis $\{e_1,\dots,e_d\}$ of $W$, corresponding to extensions $\cE_1, \dots, \cE_d \in \Ext^1_{\cD}(A,B)$. Form the pushout of $\cE_1,\dots,\cE_d$ over the common object $B$, denoted $\bigoplus_B^{1\leq i\leq d} \cE_i$:
$$
\left(0 \to B \to \bigoplus_B^{1\leq i\leq d} \cE_i \to \bigoplus_{i=1}^d Ae_i \cong A\otimes_E W \to 0\right) \in \Ext^1_{\cD}(A\otimes_E W, B).
$$
Then, $\bigoplus_B^{1\leq i\leq d} \cE_i \cong \cE^{W\text{-}\univ}$ is the universal extension associated to $W$.
\end{enumerate}
\end{prop}
\begin{proof}
For any $e\in W$, the inclusion $i_e: E\hookrightarrow W, 1\mapsto e$ induces a commutative diagram:
$$
\begin{tikzcd}
\Hom_E(W, \Ext^1_\cD(A,B)) \arrow[r, "\sim"] \arrow[d, "i_e^*"] & \Ext^1_\cD(A,B)\otimes_E W^\vee \arrow[r, "\sim"] \arrow[d, "\id\otimes i_e^\vee"] & \Ext^1_\cD(A\otimes_E W, B) \arrow[d, "\alpha_e^*"] \\
\Hom_E(E, \Ext^1_\cD(A,B)) \arrow[r, "\sim"] & \Ext^1_\cD(A,B)\otimes_E E^\vee \arrow[r, "\sim"] & \Ext^1_\cD(A\otimes_E E, B).
\end{tikzcd}
$$
Keeping track of the element $i_W \in \Hom_E(W, \Ext^1_\cD(A,B))$ proves (i).

As for the second statement, we have a commutative diagram of $E$-linear isomorphisms
$$
\begin{tikzcd}
\Hom_E(W, \Ext^1_\cD(A,B)) \arrow[r, "\sim"] \arrow[d, "\bigoplus_{i=1}^d i_{e_i}^*", "\simeq"'] & \Ext^1_\cD(A,B)\otimes_E W^\vee \arrow[r, "\sim"] \arrow[d, "\bigoplus_{i=1}^d (\id\otimes i_{e_i}^\vee)", "\simeq"'] & \Ext^1_\cD(A\otimes_E W, B) \arrow[d, "\bigoplus_{i=1}^d \alpha_{e_i}^*", "\simeq"'] \\
\bigoplus_{i=1}^d \Hom_E(E, \Ext^1_\cD(A,B)) \arrow[r, "\sim"] & \bigoplus_{i=1}^d \Ext^1_\cD(A,B)\otimes_E E^\vee \arrow[r, "\sim"] & \bigoplus_{i=1}^d \Ext^1_\cD(A\otimes_E E, B)
\end{tikzcd}
$$
for the chosen basis $\{e_1,\dots,e_d\}$ of $W$.
Recall that the pushout $\bigoplus_B^{1\leq i\leq d} \cE_i$ sitting in 
$$
0 \to B \to \bigoplus_B^{1\leq i\leq d} \cE_i \to \bigoplus_{i=1}^d Ae_i \cong A\otimes_E W \to 0
$$
is by construction such that $\alpha_{e_j}^*\left(\bigoplus_B^{1\leq i\leq d} \cE_i\right) = \cE_j$ for each $1\leq j\leq d$. 
Hence, it corresponds to $i_W \in \Hom_E(W, \Ext^1_\cD(A,B))$, and it equals the universal extension attached to $W$ by (i).
\end{proof}

\subsubsection{Ding's correspondence and the proof of Theorem B}
We now recall Ding's construction of locally $\bbQ_p$-analytic $\GL_n(K)$-representations $\pi_{\min}(D)$ and $\pi_{\fs}(D)$ attached to $D\in \PGnc$, and then give a proof of Theorem B.
By \cite[Theorem 3.21]{Ding_crystabelline}, there is a unique $E$-linear map
$$
t_D: \Ext^1_G(\pi_\alg(\phi, {\bf h}), \pi_1(\phi, {\bf h})) \twoheadrightarrow \overline{\Ext}^1(D,D)
$$
characterized by the property that the following diagram is commutative
\begin{equation}
\label{eqn:diagram_tD}
\begin{tikzcd}
\bigoplus_{w\in S_n} \ovExt^1_{w}(D,D) \arrow[r, "\zeta_{w,{\bf h},1}\circ\kappa_w", "\sim"'] \arrow[d, "\text{sum}", two heads]& \bigoplus_{w\in S_n}\Ext^1_{w}(\pi_\alg(\phi,{\bf h}), \pi_1(\phi,{\bf h}))
\arrow[d, "\text{sum}", two heads] \\
\ovExt^1(D,D) & \arrow[l, "t_D", two heads] \Ext^1_{G}(\pi_\alg(\phi,{\bf h}), \pi_1(\phi,{\bf h})).
\end{tikzcd}
\end{equation}
The same map $t_D$ is used when $\pi_1(\phi,{\bf h})$ is replaced by $\pi(\phi,{\bf h})$, as in \cite[(3.16)]{Ding_crystabelline}.
In \cite[below Lemma 3.23]{Ding_crystabelline}, Ding defines $\pi_{\min}(D)$ as the universal extension $\cE^{\ker(t_D){\text{-univ}}}$ of $\ker(t_D)\otimes_E\pi_\alg(\phi,{\bf h})$ by $\pi_1(\phi,{\bf h})$ and defines $\pi_{\fs}(D)$ as the universal extension of $\ker(t_D)\otimes_E \pi_\alg(\phi,{\bf h})$ by $\pi(\phi,{\bf h})$ in the sense of \S\ref{univ_ext}.
Then, $\pi_{\fs}(D) \cong \pi_{\min}(D) \oplus_{\pi_1(\phi,{\bf h})} \pi(\phi,{\bf h})$ by \cite[(3.24)]{Ding_crystabelline}.

\begin{prop}
\label{prop:T_and_t}
For $\bullet = 1$ or empty, we have a commutative diagram
$$
\begin{tikzcd}
\overline{\Ext}^1(D,D) \arrow[d, "p_{\bf k}"] & \arrow[l, "t_D"'] \Ext^1_G(\pi_\alg(\phi,{\bf h}), \pi_\bullet(\phi,{\bf h})) \arrow[d, "T_\lambda^{\lambda'}"]\\ 
\overline{\Ext}^1(p_{\bf k}(D),p_{\bf k}(D)) & \arrow[l, "t_{p_{\bf k}(D)}"']  \Ext^1_G(\pi_\alg(\phi,{\bf h'}), \pi_\bullet(\phi,{\bf h'})).
\end{tikzcd}
$$
\end{prop}
\begin{proof}
Consider the following cube whose ``front" face is the square in the statement
$$
\begin{tikzcd}[scale cd=0.76, row sep=normal, column sep=tiny]
& {\bigoplus_{w\in S_n} \ovExt^1_{w}(D,D)} \arrow[ld, "\text{sum}", two heads] \arrow[rr, "\zeta_{w,{\bf h},\bullet} \circ \kappa_w", "\simeq"', pos=0.35] \arrow[dd, "{p_{\bf k}}", pos=0.28] & & {\bigoplus_{w\in S_n}\Ext^1_{w}(\pi_\alg(\phi,{\bf h}), \pi_\bullet(\phi,{\bf h}))} \arrow[ld, "\text{sum}", two heads] \arrow[dd,"{T_{\lambda}^{\lambda'}}", pos=0.3] \\
{\ovExt^1(D,D)} \arrow[dd, "p_{\bf k}", pos=0.4] & & {\Ext^1_{G}(\pi_\alg(\phi,{\bf h}), \pi_\bullet(\phi,{\bf h}))} \arrow[ll, "t_D", pos=0.58, two heads, crossing over] \arrow[dd, "T_{\lambda}^{\lambda'}"{xshift=0.8ex}, pos=0.3] & \\
 & {\bigoplus_{w\in S_n} \ovExt^1_{w}(p_{\bf k}(D),p_{\bf k}(D))} \arrow[ld, "\text{sum}", two heads] \arrow[rr, "\zeta_{w,{\bf h'},\bullet} \circ \kappa_w", "\simeq"', pos=0.3] & & {\bigoplus_{w\in S_n}\Ext^1_{w}(\pi_\alg(\phi,{\bf h}'), \pi_\bullet(\phi,{\bf h}'))} \arrow[ld, "\text{sum}", two heads]. \\
{\ovExt^1(p_{\bf k}(D),p_{\bf k}(D))} & & {\Ext^1_{G}(\pi_\alg(\phi,{\bf h}'), \pi_\bullet(\phi,{\bf h}'))} \arrow[ll, "t_{p_{\bf k}(D)}", pos=0.6, two heads] \arrow[from=uu, crossing over] &                       
\end{tikzcd}
$$
All other five faces commute either by Corollary \ref{cor:Ding_Ext-link}(iii) or by definition.
Since the map 
$$
\mathrm{sum}:\bigoplus_{w\in S_n} \Ext_w^1(\pi_{\alg}(\phi,{\bf h}), \pi_\bullet(\phi,{\bf h})) \twoheadrightarrow \Ext^1_G(\pi_{\alg}(\phi,{\bf h}), \pi_\bullet(\phi,{\bf h}))
$$ 
is surjective by \cite[Proposition 3.8(2)]{Ding_crystabelline}, it follows that all faces are commutative.
\end{proof}
Since $T_\lambda^{\lambda'}$ is an exact equivalence of categories, it induces an isomorphism 
$$
\Ext^1_G(\pi_{\alg}(\phi,{\bf h}),\pi_\bullet(\phi,{\bf h})) \xrightarrow{\sim} \Ext^1_G(\pi_{\alg}(\phi,{\bf h'}),\pi_\bullet(\phi,{\bf h'}))
$$
and sends the universal extension associated to a subspace $W \subset \Ext^1_G(\pi_{\alg}(\phi,{\bf h}),\pi_\bullet(\phi,{\bf h}))$ to the universal extension associated to $T_\lambda^{\lambda'}(W)$. 
By Proposition \ref{prop:T_and_t}, we have 
$$
T_\lambda^{\lambda'}(\pi_{\min}(D))=\pi_{\min}(p_{\bf k}(D)) \quad \text{ and } \quad T_\lambda^{\lambda'}(\pi_{\fs}(D))=\pi_{\fs}(p_{\bf k}(D)).
$$ 
This completes the proof of Theorem B.

\subsubsection*{Data availability}
We do not analyse or generate any datasets, because our work proceeds within a theoretical and mathematical approach. One can obtain the relevant materials from the references below.

\subsubsection*{Conflict of interest}
On behalf of all authors, the corresponding author states that there is no conflict of interest.

\printbibliography

@book {ABV-RealLanglands,
    AUTHOR = {Adams, Jeffrey and Barbasch, Dan and Vogan, Jr., David A.},
     TITLE = {The {L}anglands classification and irreducible characters for
              real reductive groups},
    SERIES = {Progress in Mathematics},
    VOLUME = {104},
 PUBLISHER = {Boston, MA etc.: Birkh{\"a}user},
      YEAR = {1992},
     PAGES = {xii+318},
      ISBN = {0-8176-3634-X},
   MRCLASS = {22-02 (22E47)},
  MRNUMBER = {1162533},
}

@article {Agrawal-Strauch,
    AUTHOR = {Agrawal, Shishir and Strauch, Matthias},
     TITLE = {From category {$\mathcal{O}^\infty$} to locally analytic
              representations},
   JOURNAL = {J. Algebra},
  FJOURNAL = {Journal of Algebra},
    VOLUME = {592},
      YEAR = {2022},
     PAGES = {169--232},
      ISSN = {0021-8693,1090-266X},
   MRCLASS = {22E50 (20G25)},
  MRNUMBER = {4342507},
MRREVIEWER = {Yixin\ Bao},
       DOI = {10.1016/j.jalgebra.2021.09.032},
       URL = {https://doi.org/10.1016/j.jalgebra.2021.09.032},
}

@article {BeCh_families,
    AUTHOR = {Bella\"iche, Jo\"el and Chenevier, Ga\"etan},
     TITLE = {Families of {G}alois representations and {S}elmer groups},
   JOURNAL = {Ast\'erisque},
  FJOURNAL = {Ast\'erisque},
    NUMBER = {324},
      YEAR = {2009},
     PAGES = {xii+314},
      ISSN = {0303-1179,2492-5926},
      ISBN = {978-2-85629-264-8},
   MRCLASS = {11F80 (11F70 11F85 11G40 11S25 14F30 14G22)},
  MRNUMBER = {2656025},
MRREVIEWER = {Adolfo\ Quir\'os},
}

@article {Bergdall,
    AUTHOR = {Bergdall, John},
     TITLE = {Paraboline variation over {$p$}-adic families of
              {$(\phi,\Gamma)$}-modules},
   JOURNAL = {Compos. Math.},
  FJOURNAL = {Compositio Mathematica},
    VOLUME = {153},
      YEAR = {2017},
    NUMBER = {1},
     PAGES = {132--174},
      ISSN = {0010-437X,1570-5846},
   MRCLASS = {11F80 (11F33 11F55 11F85 14F30)},
  MRNUMBER = {3622874},
MRREVIEWER = {Ivan\ Mati\'c},
       DOI = {10.1112/S0010437X16007831},
       URL = {https://doi.org/10.1112/S0010437X16007831},
}

@article {Berger_B-pairs,
    AUTHOR = {Berger, Laurent},
     TITLE = {Construction de {$(\phi,\Gamma)$}-modules: repr\'esentations
              {$p$}-adiques et {$B$}-paires},
   JOURNAL = {Algebra Number Theory},
  FJOURNAL = {Algebra \& Number Theory},
    VOLUME = {2},
      YEAR = {2008},
    NUMBER = {1},
     PAGES = {91--120},
      ISSN = {1937-0652,1944-7833},
   MRCLASS = {14F30 (11F80 11S15 11S20 11S25 14G20)},
  MRNUMBER = {2377364},
MRREVIEWER = {Adolfo\ Quir\'os},
       DOI = {10.2140/ant.2008.2.91},
       URL = {https://doi.org/10.2140/ant.2008.2.91},
}

@incollection {Berger_filtres,
    AUTHOR = {Berger, Laurent},
     TITLE = {\'Equations diff\'erentielles {$p$}-adiques et
              {$(\phi,N)$}-modules filtr\'es},
      NOTE = {Repr\'esentations $p$-adiques de groupes $p$-adiques. I.
              Repr\'esentations galoisiennes et $(\phi,\Gamma)$-modules},
   JOURNAL = {Ast\'erisque},
  FJOURNAL = {Ast\'erisque},
    NUMBER = {319},
      YEAR = {2008},
     PAGES = {13--38},
      ISSN = {0303-1179,2492-5926},
      ISBN = {978-2-85629-256-3},
   MRCLASS = {11F80 (12H25 14F30)},
  MRNUMBER = {2493215},
MRREVIEWER = {Michael\ M.\ Schein},
}

@book{Bosch_lectures,
 author = {Bosch, Siegfried},
 title = {Lectures on formal and rigid geometry},
 fseries = {Lecture Notes in Mathematics},
 series = {Lect. Notes Math.},
 issn = {0075-8434},
 volume = {2105},
 isbn = {978-3-319-04416-3; 978-3-319-04417-0},
 year = {2014},
 publisher = {Cham: Springer},
 language = {English},
 doi = {10.1007/978-3-319-04417-0},
 zbMATH = {6255263},
 Zbl = {1314.14002}
}

@misc{Bourbaki_Lie78,
 author = {Bourbaki, N.},
 title = {{\'E}l{\'e}ments de math{\'e}matique. {Fasc}. {XXXVIII}: {Groupes} et alg{\`e}bres de {Lie}. {Chap}. {VII}: {Sous}-alg{\`e}bres de {Cartan}, {\'e}l{\'e}ments r{\'e}guliers. {Chap}. {VIII}: {Alg{\`e}bres} de {Lie} semi-simples d{\'e}ploy{\'e}es},
 year = {1975},
 language = {French},
 howpublished = {Actualit{\'e}s {Scientifiques} et {Industrielles}, 1364. {Paris}: {Hermann}. 271 p. {F} 78.00 (1975).},
 zbMATH = {3515658},
 Zbl = {0329.17002}
}

@article{Breuil_socle1,
 author = {Breuil, Christophe},
 title = {Locally analytic socle. {I}.},
 fjournal = {Annales de l'Institut Fourier},
 journal = {Ann. Inst. Fourier},
 issn = {0373-0956},
 volume = {66},
 number = {2},
 pages = {633--685},
 year = {2016},
 language = {French},
 doi = {10.5802/aif.3021},
 zbMATH = {6644876},
 Zbl = {1378.11098}
}

@article {BM-multiplicity,
    AUTHOR = {Breuil, Christophe and M\'ezard, Ariane},
     TITLE = {Multiplicit\'es modulaires et repr\'esentations de {${\rm
              GL}_2({\bf Z}_p)$} et de {${\rm Gal}(\overline{\bf Q}_p/{\bf
              Q}_p)$} en {$l=p$}},
      NOTE = {With an appendix by Guy Henniart},
   JOURNAL = {Duke Math. J.},
  FJOURNAL = {Duke Mathematical Journal},
    VOLUME = {115},
      YEAR = {2002},
    NUMBER = {2},
     PAGES = {205--310},
      ISSN = {0012-7094,1547-7398},
   MRCLASS = {11F80 (11F85 11S23)},
  MRNUMBER = {1944572},
MRREVIEWER = {Jacques\ Tilouine},
       DOI = {10.1215/S0012-7094-02-11522-1},
       URL = {https://doi.org/10.1215/S0012-7094-02-11522-1},
}

@article {BS_functoriality,
    AUTHOR = {Breuil, Christophe and Schneider, Peter},
     TITLE = {First steps towards {$p$}-adic {L}anglands functoriality},
   JOURNAL = {J. Reine Angew. Math.},
  FJOURNAL = {Journal f\"ur die Reine und Angewandte Mathematik. [Crelle's
              Journal]},
    VOLUME = {610},
      YEAR = {2007},
     PAGES = {149--180},
      ISSN = {0075-4102,1435-5345},
   MRCLASS = {11S37 (22E50)},
  MRNUMBER = {2359853},
MRREVIEWER = {Jacques\ Tilouine},
       DOI = {10.1515/CRELLE.2007.070},
       URL = {https://doi.org/10.1515/CRELLE.2007.070},
}

@article {BHS3,
    AUTHOR = {Breuil, Christophe and Hellmann, Eugen and Schraen, Benjamin},
     TITLE = {A local model for the trianguline variety and applications},
   JOURNAL = {Publ. Math. Inst. Hautes \'Etudes Sci.},
  FJOURNAL = {Publications Math\'ematiques. Institut de Hautes \'Etudes
              Scientifiques},
    VOLUME = {130},
      YEAR = {2019},
     PAGES = {299--412},
      ISSN = {0073-8301,1618-1913},
   MRCLASS = {11F70 (14D24 22E60)},
  MRNUMBER = {4028517},
MRREVIEWER = {Yiwen\ Ding},
       DOI = {10.1007/s10240-019-00111-y},
       URL = {https://doi.org/10.1007/s10240-019-00111-y},
}

@incollection{Buzzard_eigenvar,
 author = {Buzzard, Kevin},
 title = {Eigenvarieties},
 booktitle = {\(L\)-functions and Galois representations. Based on the symposium, Durham, UK, July 19--30, 2004},
 isbn = {978-0-521-69415-5},
 pages = {59--120},
 year = {2007},
 publisher = {Cambridge: Cambridge University Press},
 language = {English},
 zbMATH = {5282812},
 Zbl = {1230.11054}
}

@article {Colmez_poids,
    AUTHOR = {Colmez, Pierre},
     TITLE = {Correspondance de {L}anglands locale {$p$}-adique et
              changement de poids},
   JOURNAL = {J. Eur. Math. Soc. (JEMS)},
  FJOURNAL = {Journal of the European Mathematical Society (JEMS)},
    VOLUME = {21},
      YEAR = {2019},
    NUMBER = {3},
     PAGES = {797--838},
      ISSN = {1435-9855,1435-9863},
   MRCLASS = {11S37 (22E50)},
  MRNUMBER = {3908766},
MRREVIEWER = {Geo\ Kam-Fai\ Tam},
       DOI = {10.4171/JEMS/851},
       URL = {https://doi.org/10.4171/JEMS/851},
}

@article {Colmez_localement,
    AUTHOR = {Colmez, Pierre},
     TITLE = {Repr\'esentations localement analytiques de {${\bf GL}_2({\bf
              Q}_p)$} et {$(\varphi,\Gamma)$}-modules},
   JOURNAL = {Represent. Theory},
  FJOURNAL = {Representation Theory. An Electronic Journal of the American
              Mathematical Society},
    VOLUME = {20},
      YEAR = {2016},
     PAGES = {187--248},
      ISSN = {1088-4165},
   MRCLASS = {22E50 (11F70 14D24)},
  MRNUMBER = {3522263},
MRREVIEWER = {Nadir\ Matringe},
       DOI = {10.1090/ert/484},
       URL = {https://doi.org/10.1090/ert/484},
}

@incollection {Colmez_trianguline,
    AUTHOR = {Colmez, Pierre},
     TITLE = {Repr\'esentations triangulines de dimension 2},
      NOTE = {Repr\'esentations $p$-adiques de groupes $p$-adiques. I.
              Repr\'esentations galoisiennes et $(\phi,\Gamma)$-modules},
   JOURNAL = {Ast\'erisque},
  FJOURNAL = {Ast\'erisque},
    NUMBER = {319},
      YEAR = {2008},
     PAGES = {213--258},
      ISSN = {0303-1179,2492-5926},
      ISBN = {978-2-85629-256-3},
   MRCLASS = {11S37 (11S23)},
  MRNUMBER = {2493219},
MRREVIEWER = {Jan\ Nekov\'a\v r},
}

@unpublished{Ding_wallcross,
    author = {Ding, Yiwen},
     title = {Towards the wall-crossing of locally $\mathbb{Q}_p$-analytic representations of $\mathrm{GL}_{n}(K)$ for a p-adic field $K$},
      note = {Preprint, v1, 2024, \url{https://arxiv.org/abs/2404.06315}}
}

@unpublished{Ding_weight,
    author = {Ding, Yiwen},
     title = {Change of weights for locally analytic representations of $\mathrm{GL}_2(\mathbb{Q}_p)$},
      note = {Preprint, v1, 2023, \url{https://arxiv.org/abs/2307.04332}}
}

@article {Ding_crystabelline,
    AUTHOR = {Ding, Yiwen},
     TITLE = {{$p$}-adic {H}odge parameters in the crystabelline
              representations of {${\rm GL}_{n}$}},
   JOURNAL = {Publ. Math. Inst. Hautes \'Etudes Sci.},
  FJOURNAL = {Publications Math\'ematiques. Institut de Hautes \'Etudes
              Scientifiques},
    VOLUME = {142},
      YEAR = {2025},
     PAGES = {1--74},
      ISSN = {0073-8301,1618-1913},
   MRCLASS = {22E50 (11F70 11S37)},
  MRNUMBER = {4992849},
       DOI = {10.1007/s10240-025-00156-2},
       URL = {https://doi.org/10.1007/s10240-025-00156-2},
}

@unpublished{EGH_CatLLC,
    author = {Emerton, Matthew and Gee, Toby and Hellmann, Eugen},
     title = {An introduction to the categorical $p$-adic Langlands program},
      note = {Preprint, v5, 2025, \url{https://arxiv.org/abs/2210.01404v5}}
}

@book {EG_modulistack,
    AUTHOR = {Emerton, Matthew and Gee, Toby},
     TITLE = {Moduli stacks of \'etale ({$\varphi, \Gamma$})-modules and the
              existence of crystalline lifts},
    SERIES = {Annals of Mathematics Studies},
    VOLUME = {215},
 PUBLISHER = {Princeton University Press, Princeton, NJ},
      YEAR = {[2023] \copyright 2023},
     PAGES = {ix+298},
      ISBN = {978-0-691-24134-0; 978-0-691-24135-7; 978-0-691-24136-4},
   MRCLASS = {14D23 (11F80 11S37 14F30)},
  MRNUMBER = {4529886},
MRREVIEWER = {Eran\ Assaf},
       DOI = {10.1515/9780691241364},
       URL = {https://doi.org/10.1515/9780691241364},
}

@incollection {Fontaine_arithmetique,
    AUTHOR = {Fontaine, Jean-Marc},
     TITLE = {Arithm\'etique des repr\'esentations galoisiennes
              {$p$}-adiques},
      NOTE = {Cohomologies $p$-adiques et applications arithm\'etiques. III},
   JOURNAL = {Ast\'erisque},
  FJOURNAL = {Ast\'erisque},
    NUMBER = {295},
      YEAR = {2004},
     PAGES = {xi, 1--115},
      ISSN = {0303-1179,2492-5926},
   MRCLASS = {11F80 (11F85 11S15 11S20 11S25 13K05 14F30)},
  MRNUMBER = {2104360},
MRREVIEWER = {Laurent\ N.\ Berger},
}

@incollection {Fontaine_p-adiques_semi-stables,
    AUTHOR = {Fontaine, Jean-Marc},
     TITLE = {Repr\'esentations {$p$}-adiques semi-stables},
      NOTE = {With an appendix by Pierre Colmez,
              P\'eriodes $p$-adiques (Bures-sur-Yvette, 1988)},
   JOURNAL = {Ast\'erisque},
  FJOURNAL = {Ast\'erisque},
    NUMBER = {223},
      YEAR = {1994},
     PAGES = {113--184},
      ISSN = {0303-1179,2492-5926},
   MRCLASS = {14F30 (14F40)},
  MRNUMBER = {1293972},
MRREVIEWER = {Burt\ Totaro},
}

@article {HS_density,
    AUTHOR = {Hellmann, Eugen and Schraen, Benjamin},
     TITLE = {Density of potentially crystalline representations of fixed
              weight},
   JOURNAL = {Compos. Math.},
  FJOURNAL = {Compositio Mathematica},
    VOLUME = {152},
      YEAR = {2016},
    NUMBER = {8},
     PAGES = {1609--1647},
      ISSN = {0010-437X,1570-5846},
   MRCLASS = {11S20 (11F80 11F85)},
  MRNUMBER = {3542488},
MRREVIEWER = {Alan\ Koch},
       DOI = {10.1112/S0010437X16007363},
       URL = {https://doi.org/10.1112/S0010437X16007363},
}

@book {Humphreys_CategoryO,
    AUTHOR = {Humphreys, James E.},
     TITLE = {Representations of semisimple {L}ie algebras in the {BGG} category {$\mathscr{O}$}},
    SERIES = {Graduate Studies in Mathematics},
    VOLUME = {94},
 PUBLISHER = {American Mathematical Society, Providence, RI},
      YEAR = {2008},
     PAGES = {xvi+289},
      ISBN = {978-0-8218-4678-0},
   MRCLASS = {17B10},
  MRNUMBER = {2428237},
MRREVIEWER = {Serge\ M.\ Skryabin},
       DOI = {10.1090/gsm/094},
       URL = {https://doi.org/10.1090/gsm/094},
}

@article {JLS_translation,
    AUTHOR = {Jena, Akash and Lahiri, Aranya and Strauch, Matthias},
     TITLE = {Translation functors for locally analytic representations},
   JOURNAL = {Int. J. Number Theory},
  FJOURNAL = {International Journal of Number Theory},
    VOLUME = {20},
      YEAR = {2024},
    NUMBER = {10},
     PAGES = {2611--2640},
      ISSN = {1793-0421,1793-7310},
   MRCLASS = {17B10 (11F85 11S37 22E50)},
  MRNUMBER = {4833414},
       DOI = {10.1142/S1793042124501252},
       URL = {https://doi.org/10.1142/S1793042124501252},
}

@book {KV_coh_induction,
    AUTHOR = {Knapp, Anthony W. and Vogan, Jr., David A.},
     TITLE = {Cohomological induction and unitary representations},
    SERIES = {Princeton Mathematical Series},
    VOLUME = {45},
 PUBLISHER = {Princeton University Press, Princeton, NJ},
      YEAR = {1995},
     PAGES = {xx+948},
      ISBN = {0-691-03756-6},
   MRCLASS = {22E46 (22-02 22E47)},
  MRNUMBER = {1330919},
MRREVIEWER = {William\ M.\ McGovern},
       DOI = {10.1515/9781400883936},
       URL = {https://doi.org/10.1515/9781400883936},
}

@article {KPX_cohomology,
    AUTHOR = {Kedlaya, Kiran S. and Pottharst, Jonathan and Xiao, Liang},
     TITLE = {Cohomology of arithmetic families of
              {$(\varphi,\Gamma)$}-modules},
   JOURNAL = {J. Amer. Math. Soc.},
  FJOURNAL = {Journal of the American Mathematical Society},
    VOLUME = {27},
      YEAR = {2014},
    NUMBER = {4},
     PAGES = {1043--1115},
      ISSN = {0894-0347,1088-6834},
   MRCLASS = {11F33 (11R23 11S25 11S31 13D09)},
  MRNUMBER = {3230818},
MRREVIEWER = {Th\cfac ong\ Nguy\cftil en-Quang-\Dbar\cftil o},
       DOI = {10.1090/S0894-0347-2014-00794-3},
       URL = {https://doi.org/10.1090/S0894-0347-2014-00794-3},
}

@article {Liu_cohomology,
    AUTHOR = {Liu, Ruochuan},
     TITLE = {Cohomology and duality for {$(\phi,\Gamma)$}-modules over the
              {R}obba ring},
   JOURNAL = {Int. Math. Res. Not. IMRN},
  FJOURNAL = {International Mathematics Research Notices. IMRN},
      YEAR = {2008},
    NUMBER = {3},
     PAGES = {Art. ID rnm150, 32},
      ISSN = {1073-7928,1687-0247},
   MRCLASS = {11S25 (11F80 11S20 12G05)},
  MRNUMBER = {2416996},
MRREVIEWER = {Laurent\ N.\ Berger},
       DOI = {10.1093/imrn/rnm150},
       URL = {https://doi.org/10.1093/imrn/rnm150},
}

@article {Nakamura_classification,
    AUTHOR = {Nakamura, Kentaro},
     TITLE = {Classification of two-dimensional split trianguline
              representations of {$p$}-adic fields},
   JOURNAL = {Compos. Math.},
  FJOURNAL = {Compositio Mathematica},
    VOLUME = {145},
      YEAR = {2009},
    NUMBER = {4},
     PAGES = {865--914},
      ISSN = {0010-437X,1570-5846},
   MRCLASS = {11S25 (11F80 11F85)},
  MRNUMBER = {2521248},
MRREVIEWER = {David\ L.\ Savitt},
       DOI = {10.1112/S0010437X09004059},
       URL = {https://doi.org/10.1112/S0010437X09004059},
}

@article {Nakamura_B-pairs,
    AUTHOR = {Nakamura, Kentaro},
     TITLE = {Deformations of trianguline {$B$}-pairs and {Z}ariski density
              of two dimensional crystalline representations},
   JOURNAL = {J. Math. Sci. Univ. Tokyo},
  FJOURNAL = {The University of Tokyo. Journal of Mathematical Sciences},
    VOLUME = {20},
      YEAR = {2013},
    NUMBER = {4},
     PAGES = {461--568},
      ISSN = {1340-5705},
   MRCLASS = {11F80 (11F85 11S25)},
  MRNUMBER = {3185293},
MRREVIEWER = {Stefano\ Vigni},
}

@article {Schneider-Teitelbaum_duality,
    AUTHOR = {Schneider, Peter and Teitelbaum, Jeremy},
     TITLE = {Duality for admissible locally analytic representations},
   JOURNAL = {Represent. Theory},
  FJOURNAL = {Representation Theory. An Electronic Journal of the American
              Mathematical Society},
    VOLUME = {9},
      YEAR = {2005},
     PAGES = {297--326},
      ISSN = {1088-4165},
   MRCLASS = {22E50},
  MRNUMBER = {2133762},
MRREVIEWER = {Anne-Marie\ H.\ Aubert},
       DOI = {10.1090/S1088-4165-05-00277-3},
       URL = {https://doi.org/10.1090/S1088-4165-05-00277-3},
}

@article{Schneider-Teitelbaum_GL2,
 author = {Schneider, Peter and Teitelbaum, Jeremy},
 title = {Locally analytic distributions and {{\(p\)}}-adic representation theory, with applications to {{\(\text{GL}_{2}\)}}},
 fjournal = {Journal of the American Mathematical Society},
 journal = {J. Am. Math. Soc.},
 issn = {0894-0347},
 volume = {15},
 number = {2},
 pages = {443--468},
 year = {2002},
 language = {English},
 doi = {10.1090/S0894-0347-01-00377-0},
 zbMATH = {1720953},
 Zbl = {1028.11071}
}

@misc{stacks-project,
    shorthand    = {Stacks},
    author       = {The {Stacks Project Authors}},
    title        = {\textit{Stacks Project}},
    howpublished = {\url{https://stacks.math.columbia.edu}},
    year         = {2018},
  }

@book {Weibel_HA,
    AUTHOR = {Weibel, Charles A.},
     TITLE = {An introduction to homological algebra},
    SERIES = {Cambridge Studies in Advanced Mathematics},
    VOLUME = {38},
 PUBLISHER = {Cambridge University Press, Cambridge},
      YEAR = {1994},
     PAGES = {xiv+450},
      ISBN = {0-521-43500-5; 0-521-55987-1},
   MRCLASS = {18-01 (16-01 17-01 20-01 55Uxx)},
  MRNUMBER = {1269324},
MRREVIEWER = {Kenneth\ A.\ Brown},
       DOI = {10.1017/CBO9781139644136},
       URL = {https://doi.org/10.1017/CBO9781139644136},
}

@phdthesis{Wu_thesis,
  TITLE = {{Trianguline variety and eigenvariety at points with non-regular Hodge-Tate weights}},
  AUTHOR = {Wu, Zhixiang},
  URL = {https://theses.hal.science/tel-04116114},
  NUMBER = {2022UPASM015},
  SCHOOL = {{Universit{\'e} Paris-Saclay}},
  YEAR = {2022},
  MONTH = Jul,
  TYPE = {Theses},
  HAL_ID = {tel-04116114},
  HAL_VERSION = {v1},
}

@unpublished{Wu_translation,
    author = {Wu, Zhixiang},
     title = {Geometric translations of $(\varphi,\Gamma)$-modules for $\mathrm{GL}_2(\mathbb{Q}_p)$},
      note = {Preprint, v2, 2025, \url{https://arxiv.org/abs/2405.16637v2}}
}

@incollection {Zhu_affineGrassmann,
    AUTHOR = {Zhu, Xinwen},
     TITLE = {An introduction to affine {G}rassmannians and the geometric
              {S}atake equivalence},
 BOOKTITLE = {Geometry of moduli spaces and representation theory},
    SERIES = {IAS/Park City Math. Ser.},
    VOLUME = {24},
     PAGES = {59--154},
 PUBLISHER = {Amer. Math. Soc., Providence, RI},
      YEAR = {2017},
      ISBN = {978-1-4704-3574-5},
   MRCLASS = {14M15 (14D24 20F65 22E57)},
  MRNUMBER = {3752460},
MRREVIEWER = {Felipe\ Zald\'ivar},
}

\end{document}